\documentclass{article}
\usepackage[utf8]{inputenc}
\usepackage[a4paper, total={7.1in, 8in}]{geometry}
\pdfoutput=1

\usepackage{amsmath, bm}
\usepackage{amsfonts}
\usepackage{amsthm}
\usepackage{amssymb}
\usepackage{bbm}
\usepackage{ragged2e}
\usepackage{lineno}
\usepackage{comment}
\usepackage[ruled]{algorithm2e}
\usepackage{graphicx,subfig}
\usepackage{pgfplots}
\usetikzlibrary{spy}
\usepackage{hyperref}
\usepackage{todonotes}

\newtheorem{thm}{Theorem}[section]
\newtheorem{lem}[thm]{Lemma}

\theoremstyle{remark}
\newtheorem*{rem}{\rm\bf{Remarks}}

\let\oldproofname=\proofname
\renewcommand{\proofname}{\rm\bf{\oldproofname}}

\newcommand{\grad}{\nabla}

\newcommand{\R}{\mathbb{R}}
\newcommand{\N}{\mathbb{N}}

\newcommand{\p}{\mathbf{p}}

\newcommand{\argmin}[1]{\underset{#1}{\text{{\rm argmin} }}}
\renewcommand*{\sup}[1]{\underset{#1}{\text{{\rm sup} }}}
\renewcommand*{\inf}[1]{\underset{#1}{\text{{\rm inf} }}}
\renewcommand*{\liminf}[1]{\underset{#1}{\text{{\rm liminf} }}}

\newcommand{\dom}{\mathfrak{D}}

\newcommand{\lam}{\lambda}

\newcommand{\lap}{\Delta}

\renewcommand{\th}{{{\bf{\theta}}}}
\renewcommand{\u}{\mathbf{u}}
\newcommand{\Ra}{\Rightarrow}

\newcommand{\eps}{\epsilon}

\newcommand{\weakly}{\rightharpoonup}

\newcommand{\Loo}{L^2(\dom)}
\newcommand{\Ho}{H^1(\dom)}
\newcommand{\Hoo}{H^2(\dom)}
\newcommand{\HooN}{H^2_N(\dom)}

\newcommand{\Htp}{H^2_p(\dom)}

\newcommand{\Linf}{L^{\infty}(\dom)}
\newcommand{\LinfT}{L^{\infty}_T(\dom)}
\newcommand{\Li}{L^{\infty}}

\newcommand{\intdom}{\int_{\dom}}
\newcommand{\intdomk}{\int_{\dom_{u_k}}}

\newcommand{\domT}{\dom_T}
\newcommand{\domti}{\dom_{t_i}}
\newcommand{\domto}{\dom_{t_1}}

\newcommand{\CLTi}[2]{C([t_{i-1},t_i];L_2(\dom)}
\newcommand{\CHTi}[2]{C([t_{i-1},t_i];H^1(\dom)}

\newcommand{\chti}{C_{H^1_{t_i}}}
\newcommand{\chto}{C_{H^1_{t_1}}}
\newcommand{\hot}{H^1(\dom)}
\newcommand{\hott}{H^2(\dom)}
\newcommand{\LtHsp}[1]{L^2(I;{#1})}
\newcommand{\Z}{\mathbf{Z}}
\newcommand{\LtHspsq}{L^2(I;\Z)}

\newcommand{\X}{\mathbf{X}}

\newcommand{\const}{{{\large{\rm k}}}}

\newcommand{\lbd}{{\rm m}}
\newcommand{\ubd}{{\rm M}}

\newcommand{\Th}{\Theta}
\newcommand{\THad}{\mathbf{\Xi}}

\title{Data Driven Modeling of Pseudopalisade Pattern Formation}

\author{Sandesh Athni Hiremath\footnote{ {Mechanical and Process Engineering}, {TU Kaiserslautern}, Gottlieb-Daimler-Straße 42, {Kaiserslautern}, {67663}, Germany } \:(sandesh.hiremath@mv.uni-kl.de) \\ and \\ Christina Surulescu\footnote{{Felix-Klein-Zentrum für Mathematik}, {TU Kaiserslautern}, Paul-Ehrlich-Str. 31, {Kaiserslautern}, {67663}, Germany} \:(surulescu@mathematik.uni-kl.de)}

\date{}

\begin{document}
\maketitle

\abstract{Pseudopalisading is an interesting phenomenon where cancer cells arrange themselves to form a dense garland-like pattern. Unlike the palisade structure, a similar type of pattern first observed in schwannomas, by pathologist J.J. Verocay \cite{Wippold2006}, pseudopalisades are less organized and associated with a necrotic region at their core. These structures are mainly found in glioblastoma (GBM), a grade IV brain tumor, and provide a way to assess the aggressiveness of the tumor. Identification of the exact bio-mechanism responsible for the formation of pseudopalisades is a difficult task, mainly because pseudopalisades seem to be a consequence of complex nonlinear dynamics within the tumor. In this paper we propose a data-driven methodology to gain insight into the formation of different types of pseudopalisade structures. To this end, we start from a state of the art macroscopic model for the dynamics of GBM, that is coupled with the dynamics of extracellular pH, and formulate a terminal value optimal control problem. Thus, given a specific, observed pseudopalisade pattern, we determine the evolution of parameters (bio-mechanisms) that are responsible for its emergence. Random histological images exhibiting pseudopalisade-like structures are chosen to serve as target pattern. Having identified the optimal model parameters that generate the specified target pattern, we then formulate two different types of pattern counteracting ansatzes in order to determine possible ways to impair or obstruct the process of pseudopalisade formation. This provides the basis for designing active or live control of malignant GBM. Furthermore, we also provide a simple, yet insightful, mechanism to synthesize new pseudopalisade patterns by  linearly combining the optimal model parameters responsible for generating different known target patterns. This particularly provides a hint that complex pseudopalisade patterns could be synthesized by a linear combination of parameters responsible for generating simple patterns. Going even further, we ask ourselves if complex therapy approaches can be conceived by linearly combining such that are able to reverse or disrupt simple pseudopalisade patterns, which is then positively answered with the help of numerical simulations.}


\section{Introduction}
Biological phenomena produce some of the visually most appealing patterns, but unfortunately not all of them can be associated with a beneficial outcome. For example, skin and tissue patterns of animals  could actually indicate the onset/progression of a harmful process. Pseudopalisades belong to such category,  wherein the microscopic cellular arrangement, although both visually and dynamically quite intriguing, 
indicates the most advanced stage of glioblastoma multiforme (GBM), a type of brain tumor, which in most cases is lethal.
Such patterns are actually used to pathologically characterize the aggressiveness or malignancy of the tumor \cite{Wippold2006}. Unlike the highly regular palisade structure observed in schwannoma cells, pseudopalisade structures are less organized and more irregular in appearance. The initiation of such pathological structures is not very clear, but is mainly hypothesized (see \cite{Brat2004}) to be a complex interaction of different biophysical processes such as: (i) rapidly proliferating neoplastic cells, (ii) cells being highly resistant to apoptosis and (iii)  cells migrating away from the toxic debris formed by cellular necrosis, as a consequence of hypoxia and acidosis.  \\
In contrast, the microenvironment surrounding pseudopalisades is fairly better understood. According to the studies  \cite{Rong2006,Wippold2006,Brat2004,Gonzalez2012} the cells made of such structures are mainly hypoxic and have less proliferating capabilities. These cells, however, show increased vascular endothelial growth factor (VEGF) expression, that results in development of microvascular structures \cite{Zagzag2000,Plate1992}. Due to the dense structure of the brain tissue, this additional vascular growth is very irregular and even results in the formation of glomeruloid bodies. The area enclosed by pseudopalisades is composed of mainly dead cells and other cellular debris forming the necrotic core. Moreover, due to the vascular aberrations there exist anisotropic oxygen gradients, with the center being hypoxic.\\ 
As hypoxia is closely associated with acidosis \cite{Horn2007,Chiche2010,Jing2019}, it is well possible that the migratory cells have switched to a glycolytic pathway,  which in turn exacerbates the micro-acidity and promotes migratory behavior of the cells \cite{Estrella2013,Piasentin2020}. Because GBM is the most dominant type of malignant brain tumors \cite{Dolecek2012} and since the detection of pseudopalisades indicates a worsening condition of a glioma patient\cite{Brat2003,Kleihues1995}, it is highly important to get insight into the formation and behavior of these structures. 
For this purpose, mathematical models have proven to be highly effective, especially for understanding and validating the dynamics of biological processes. In the context of GBM, various types of models have been proposed. Nice reviews  on the chronological evolution of such models can be found in \cite{Hatzikirou2005, Harpold2007, Martirosyan2015, Alfonso2017}. Broadly speaking, there are mainly two classes of models: discrete and continuous. The former mainly comprise rule-based computational models (e.g. \cite{Sander2002, Khain2011, Bottger2012, Kim2013} that try to identify the self-organization behavior of the cells. These models take advantage of the computational power to explore different rule-configurations that could explain the phenomena.
On the other hand, continuous models are based on continuous abstraction of the evolution of physical processes. The most simple, yet effective continuous models are ODE based. They are employed not only to study the proliferating capabilities of glioma \cite{Sturrock2015}, but also for assessing the effects of radio- and chemotherapy \cite{Yu2021}. However, when one is interested in studying the invasive behavior of GBMs, space becomes important, thus spatial dynamics needs to be taken into account. Most of the cancer invasion models, inspired by the early works of \cite{Murray2002}, are based on reaction-diffusion equations, see e.g.,  \cite{Jbabdi2005,Swanson2011,Hatzikirou2012,Gonzalez2012, Kim2013, Alfonso2016}. They  only consider movement based on random motion with very limited ability to incorporate direction/orientation information from the microenvironment, e.g. only via some anisotropic diffusion coefficient. These models were generalized in \cite{Hinow2008,Kim2009,Colombo2015}, where advection/taxis terms were introduced to incorporate the relevant microenvironment information such as tissue structure, vasculature etc.\\
Because cancer growth and spread is a complex multiscale process, mere macroscopic models fail to illicit the outcomes of cross-scale interactions. 
Many cellular motion models originate at the subcellular or cellular scale by first considering the dynamics of individual cells followed by modeling the interactions with other cells and physical/chemical components of the environment. This is then upscaled to the tissue level, where experimental observations are possible. Such multiscale framework has been considered in \cite{Painter2013,Engwer2015,Engwer2016,Engwer2016a,Hunt2016, Corbin2018, Conte2020, Conte2021, Corbin2021, Dietrich2020, Conte2022} to study the invasive patterns of glioma. There, mainly the parabolic scaling is used to obtain a corresponding macroscopic PDE which consequently involves diffusion and reaction coefficients that are coupled with the dynamics of the lower scales. In contrast to the deterministic models, authors in \cite{Hiremath2015,Hiremath2016,Hiremath2017,Hiremath2018} have used stochastic multiscale settings to illicit transient invasive patterns of cancer. \\
The kind of models considered so far are in some sense phenomenological descriptions that try to explain or justify \emph{in vivo} or \emph{in vitro} observations. This is a bottom-up approach, wherein theoretical reasoning is used to explain the observed data. In contrast, one could resort to statistical techniques to infer relevant properties of the dynamics directly from the data without considering any biophysical model. Alternatively, one could simply complement the phenomenological model by statistically incorporating the observed data.  Such models are called data-driven, where the data and the (bio)physics are coupled through an optimization formulation. This type of inverse problem formulation has been used e.g., in  \cite{Hogea2008,Konukoglu2010,Gholami2016} to estimate patient specific model parameters that can subsequently be used for making predictions. Similar to these approaches, in this paper we formulate an optimal control problem with the aim of gaining insight into the dynamical processes responsible for generating a specific type of pseudopalisade patterns. Unlike previous studies \cite{Kim2009,Caiazzo2015,Gonzalez2012,Kumar2021,KSZ21} where a more or less phenomenological approach was employed, we adopt a data-driven approach where all model unknowns along with the involved model parameters are estimated from the data itself. Starting from a macroscopic model \cite{Kumar2021}, which in turn is obtained using a multiscale modeling technique, given some arbitrary initial condition and a target pseudopalisade pattern, we compute the optimal model parameters such as growth rate, diffusion coefficient, taxis direction, such that the initial tumor density optimally evolves to the final pseudopalisade pattern. The advantages of this approach are the following:
\begin{enumerate}
    \item given some fixed arbitrary initial condition $\u_0$ and different target (final)  patterns $(O_k)_{k\in\N}$, we can compute corresponding optimal model parameters $\th_{O_k}$ and solutions $\u_{O_k}$ (see \eqref{eq:pseudOCP}) which are able to directly explain the data. By analyzing the qualitative properties of the obtained parameters one can gain insight into their interactions that eventually result in the formation of the observed structures. Furthermore, this approach also provides a way to directly compare the differences in the model parameters, thereby also the resulting dynamics. This eventually results in different end patterns. 
    
    \item by reversing the initial and target (final) conditions, we can compute what are the optimal parameters that can reverse or undo the developed pattern. A typical application of this would be in developing strategies to renormalize or neutralize the tumor micro-environment with the aim of reducing the malignancy of the developed tumor. In cases where direct intervention on the tumor is conceivable, appropriate additional parameterized equations (based on the type of intervention) can be introduced with the aim of stopping further progression of the pattern. 
    
    \item by combining the parameters $\th_{O_k}$ which are responsible for generating simple target patterns $O_k$, to obtain a new parameter vector $\th_{O'}$, thus synthesizing new unseen patterns $O'$. E.g. we can define $\th_{O'} := \sum_{k=1}^N c_k \th_{O_k}$, for $ N \in \N, c_k \in \R$,  and simulate the dynamics to synthesize a new pattern $O'$. The main application of this would be that, if one can design strategies leading to simpler pseudopalisade patterns, then a similar linear combination of interventions could putatively also work in neutralizing complex pseudopalisade patterns. Whether this is useful for conceiving new therapy approaches remains arguable, however it can help understanding the histological patterns.
\end{enumerate}
Based on the above discussion, the rest of the document is organized in the following manner. In Section \ref{sec:modeling} we present a multiscale mathematical model for the dynamics involved in the formation of general pseudopalisade patterns. In Section \ref{sec:modelwellposedness} we establish the wellposedness of the  model. Following that, in Section \ref{sec:ocpwellposedness} we formulate a terminal optimal control problem (TOCP) and establish the existence of its solution. Subsequently, in Section \ref{sec:numerics} we present the numerical results. Following this, in Section \ref{sec:thsynth} we not only discuss the application of TOCP to therapy problems but also for the synthesis of new unseen patterns and its plausible value for interventions. Finally in Section \ref{sec:conclude} we discuss the results and draw conclusions. 

\section{Modeling} \label{sec:modeling}
The goal of this section is to set up a system of equations that is able to mimick the complex interactions of cancer cells with their host tissue. Because we are interested in analyzing the influence of tissue acidity on the type of glioma patterns that emerge, we restrict the description to mainly the interactions between cancer cells and protons in the extracellular region. The latter is modeled by accounting for the dynamics of extracellular proton concentration $H$. Since protons are much smaller than cancer cells, their dynamics is much faster, thus can be  considered directly at the macroscopic level. Following this reasoning, the acid evolution is described by the following reaction-diffusion equation:
\begin{align}
\label{eq:pHe}
    \partial_t H &= D_s \lap H - \alpha H + \beta f_2(C,H),
\end{align}
where $D_s$ is the effective proton diffusion coefficient, $\alpha$ is the effective acid removal rate  by vasculature, and $\beta$ represents the effective expulsion rate  of protons by cancer cells mainly as a byproduct of glycolytic energy cycle \cite{Gatenby03}. Here $f_2(C,H) := {CH \over (1 + C^2 + H^2)^2}$ models the efflux of protons by the membrane transporters of the cell such as MCT, NHE \cite{Webb99,Webb04}. For the sake of generality, the terms $\alpha$ and $\beta$ are considered to be functions of space and time. 

On the other hand, the dynamics of cancer cells is much slower and therefore can be modelled not only by considering intracellular events, but also transmembrane and extracellular interactions. This basically results in multiscale modeling of tumor evolution which for the case of GBM has been done previously by several authors \cite{Painter2013,Engwer2015,Engwer2016,Hunt2016,Swan2018}. For our study, we refer to the more recent works \cite{Kumar2021,KS2020}, particularly to the former, where the activity of proton-specific transmembrane units was considered to deduce a kinetic transport equation for the evolution of tumor density which was subsequently upscaled. The parabolic scaling procedure resulted in a myopic-diffusion-based PDE which not only translates the averaged random fluctuations at the microscopic level to the macroscopic one, but also correctly incorporates the advection term representing directed movement of cells. The resulting parabolic PDE in the non-dimensionalized form reads:
\newcommand{\DT}{\mathbb{D}}
\begin{align}
\label{eq:can}
    \partial_t C &=  \nabla \cdot (\nabla \cdot (\DT C)) + \nabla \cdot ( \delta(H) C \DT \nabla H) + \mu C (1-C)(1-H), 
\end{align}
where $\DT$ is the anisotropic diffusion tensor, $\delta(H)$ is the pH taxis coefficient, and $\mu$ is the proliferation rate. For our study, we consider a slightly modified version of the macroscopic equation, given as:
\begin{align}
\label{eq:modCan}
    \partial_t C &=  \nabla \cdot (\sigma(t,x) \nabla C + C \nabla \kappa(t,x)) + \nabla \cdot ( \delta(t,x) C \nabla H) + \mu(t,x) f_1(C,H), 
\end{align}
where we have reduced the diffusion tensor $\DT$ and pH taxis coefficient $\delta(H)$ to space-time functions $\sigma$ and $\delta$, respectively. The growth term is modified to a bounded function $f_1(C,H) := {C(1-C)(1-H) \over (1 + C^2 + H^2)^2}$ and the rate constant $\mu$ is taken be a space-time function. Additionally, we have introduced the advection term $C \nabla \kappa$ to  model haptotactic movement described by a tissue-dependent time-varying  function $\kappa$ (to be estimated). Let $T > 0$ and $I = (0,T] \subset \mathbb{R}_+ $ be a finite time interval. Let $\dom \subset \mathbb{R}^2$ be an open bounded spatial domain with sufficiently smooth boundary. The resulting coupled PDE system is given by the following initial boundary value problem (IBVP):
\begin{subequations}
\label{eq:appPseudo}
\begin{align}
\partial_t C &= \nabla \cdot (\sigma \nabla C + C \nabla \kappa ) + \nabla \cdot (\delta \: C  \nabla H) + \mu  f_1(C,H)  \quad \text{ in } (0,T]\times \dom \label{eq:Can}\\ 
\partial_t H &= \lap H - \alpha H + \beta f_2(C,H) \quad \text{ in } (0,T]\times \dom \label{eq:He}\\ 
C(0) &= C_0, \quad H(0) = H_0 \quad \text{ in } \dom \notag \\
0 &= (\sigma \nabla C + C \nabla \kappa) \cdot \hat n \quad \text{ on } [0,T] \times \partial \dom \notag \\
0 &= \nabla H \cdot \hat n \quad \text{ on } [0,T] \times \partial \dom . \notag
\end{align}
\end{subequations}
The non-dimensionalized PDE system \eqref{eq:appPseudo} serves as the abstract macroscopic model for the underlying dynamics for the evolution of pseudopalisades in GBM under the influence of acidity. Due to the nonlinear coupling of the reaction terms and interplay between different taxis terms, the resulting dynamics of GBM can be very complex and most importantly very much dependent on the qualitative and quantitative properties of the model coefficient functions (model parameters). As a result, accurately determining the model parameters for a specific observed dynamics can be very challenging. The usual way is to look for stationary solutions  by means of linear stability and bifurcation analysis. This process, although very effective during the modeling phase to gain analytic insight, is, however, often unable to explain real and interesting experimental observations which are usually very complicated. In order to explain each observation accurately, one is usually required to formulate an inverse problem, which is in fact the paradigm of this paper. To this end, starting from \eqref{eq:appPseudo} we formulate a minimization problem whose goal is to determine the optimal parameters for the model such that the final state of the tumor closely matches the real observations. This is realized by devising an optimal control problem (OCP) for which the objective function is based on the final spatial distribution of the tumor, hence it is termed as the terminal optimal control problem (TOCP). In the following section we shall first establish wellposedness of the dynamical model and then present the corresponding TOCP for which we prove the existence of a minimizer which then paves the way for performing data based numerics. 
\newcommand{\SW}{\mathcal{S}}
\newcommand{\Y}{\mathbf{Y}}
\newcommand{\U}{\mathbf{U}}
\renewcommand{\TH}{\mathbf{\Th}}
\newcommand{\LtHspsqstar}{{L^2(I;\Z^*)}}
\newcommand{\B}{{\bf B}}
\newcommand{\bvsig}{{\bf \varsigma}}
\newcommand{\bbar}{\Big\|}
\newcommand{\Bbar}{\Bigg\|}
\newcommand{\Zloc}{\mathcal{Z}}
\newcommand{\Zzero}{Z^0}
\newcommand{\Zstar}{Z^*}
\newcommand{\Kloc}{\mathcal{K}}
\newcommand{\Aop}{\mathbf{A}}
\newcommand{\Rop}{\mathbf{R}}
\newcommand{\Fop}{\mathbf{F}}
\newcommand{\rop}{\mathbf{r}}
\newcommand{\aop}{\mathbf{a}}
\newcommand{\fop}{\mathbf{f}}
\newcommand{\Gu}{G_{\u}}
\newcommand{\iGu}{G^{-1}_{\u}}
\newcommand{\Guad}{G^*_{\u}}
\newcommand{\iGuad}{(G^*_{\u})^{-1}}
\newcommand{\bphi}{{\bm{\varphi}}}
\newcommand{\bthet}{{\bm{\theta}}}
\newcommand{\bvarth}{{\bm{\vartheta}}}
\renewcommand{\v}{\mathbf{v}}
\newcommand{\w}{\mathbf{w}}
\newcommand{\V}{\mathbf{V}}
\newcommand{\W}{\mathbf{W}}
\newcommand{\g}{\mathbf{g}}
\newcommand{\thetabf}{{\boldsymbol \theta}}
\newcommand{\q}{\mathbf{q}}
\newcommand{\Csp}{\mathbf{\mathcal{C}}}
\newcommand{\fotil}{\tilde f_1}
\newcommand{\fttil}{\tilde f_2}
\newcommand{\blam}{{\bm{\lambda}}}
\newcommand{\ubdk}{\ubd_{\kappa}}
\newcommand{\ubdd}{\ubd_{\delta}}
\newcommand{\ubdmu}{\ubd_{\mu}}
\newcommand{\ubdsig}{\ubd_{\sigma}}
\newcommand{\ubdalp}{\ubd_{\alpha}}
\newcommand{\lbdsig}{\lbd_{\sigma}}
\newcommand{\lbdalp}{\lbd_{\alpha}}
\newcommand{\lbdmu}{\lbd_{\mu}}
\newcommand{\lbdd}{\lbd_{\delta}}
\newcommand{\CTH}{C([0,T]; H^s_p(\dom))}

\section{Analysis}
\subsection{Assumptions and prerequisites:}
Let $I = (0,T] \subset \mathbb{R}_+ $ be a finite time interval and $\dom \subset \mathbb{R}^2$ be an open bounded spatial domain with sufficiently smooth boundary. 
Letting $ H^2_{N} := \big\{ u \in \Htp: { \partial u \over \partial \nu } = 0\big\},$ we use the following notations for the common Lebesgue and Sobolev spaces:
\begin{align}
\label{eq:funSp}
\begin{array}{lll}
    L^p := L^p(\dom,\|\cdot\|_{L^p}),& L^p_T := L^p([0,T] \times \dom, \|\cdot\|_{L^p_T}), & L^p_T \equiv L^p([0,T];L^p(\dom)),\\
    V :=  \Big( L^2(\dom), (\cdot, \cdot)\Big),& Z^2 := \HooN,\ Z := Z^1 := \Hoo,& W := \Ho,\ \W := W \times W,\\
    Z^2(T) := \LtHsp{Z^2},\qquad & \V := V \times V,\ \Z := Z^1 \times Z^2,& Z(T) := Z^1(T) := \LtHsp{Z},\\
     W(T) := \LtHsp{W},& \V(T) := \LtHsp{\V},\Z(T) := \LtHspsq,& \W(T) := \LtHsp{\W},\\
     \Y := L^{\infty}(I\times \dom),& \Y(T) := L^{\infty}(I;\Z),& \Csp(T) := C(\bar I;\Z),\\
\end{array}\notag\\
X := \{ u \in \LtHsp{W} , \: u' \in \LtHsp{W'} \}, \quad \X := \{ \u \in \LtHsp{\W}: \u' \in \LtHsp{\W'} \},
\end{align} 
with $W'$ denoting the dual space of $W$. Also, we denote the space of linear operators from $\U$ to $\U' \times \Z$ by $L(\U ; \U' \times \Z)$.

\noindent
Finally, we define the solution space $\U$ and the parameter space $\TH$ as 
\begin{align}
\label{eq:solSp}
    \U := \X \cap \Y(T) \cap \Z(T), \quad \text{and} \quad \TH := \Big(\LtHsp{ Z^{\times 6}}, (\cdot,\cdot)_{\TH}\Big). 
\end{align} 



\paragraph{Formulation of the data-driven model}
The data-driven model comprises two main components: firstly, the model for the dynamics of the system and secondly, the objective/cost functional that couples the observation data with the system dynamics. The former is described by \eqref{eq:appPseudo} which when coupled with a terminal type objective functional results in a terminal valued optimal control problem (TOCP). We shall now reformulate \eqref{eq:appPseudo} in an abstract form so that it enables TOCP to be represented in a way that is conducive for mathematical analysis.
To this end, let $\u = (u_1, u_2)^{\top}$ represent the cancer density $C$ and extracellular acidity $H$, then equation \eqref{eq:appPseudo} can be rewritten as
\begin{subequations}
\label{eq:absPseudo}
\begin{align}
\partial_t u_1 - \nabla \cdot (\sigma \nabla u_1 + u_1 \nabla \kappa ) &= \nabla \cdot (\delta \: u_1  \nabla u_2)  + \mu  f_1(u_1, u_2) \quad \text{ in } (0,T]\times \dom \label{eq:Can}\\ 
\partial_t u_2 - \lap u_2 + \alpha u_2 &=  \beta f_2(u_1,u_2)  \quad \text{ in } (0,T]\times \dom \label{eq:He}\\ 
u_1(0) &= u_{1_0}, \quad u_2(0) = u_{2_0} \quad \text{ in } \dom \notag \\
(\sigma \nabla u_1 + u_1 \nabla \kappa) \cdot \hat n &= 0,  \quad \nabla u_2 \cdot \hat n = 0 \quad \text{ on } [0,T] \times \partial \dom, \notag
\end{align}
\end{subequations}
where $f_1(u_1,u_2) := {u_1 (1 - u_1) (1 - u_2) \over (1 + u_1^2 + u_2^2)^2}$,  $f_2(u_1,u_2) := {u_1 u_2 \over (1 + u_1^2 + u_2^2)^2}$. \\
Letting $\hat f_1 := f_1+ \mu u_1$, the weak formulation of \eqref{eq:absPseudo} is given as:
\begin{subequations}
\label{eq:weakPseudo}
\begin{align}
    (\partial_t u_1, \varphi) + (\sigma \nabla u_1, \nabla \varphi) + (\mu u_1, \varphi) +  (\delta u_1 \nabla u_2, \nabla \varphi) &= (\mu \hat f_1, \varphi) - (u_1 \nabla \kappa, \nabla \varphi), \\
    (\partial_t u_2, \psi) + (\nabla u_2 , \nabla \psi) + (\alpha u_2, \psi) &= (\beta f_2, \psi), \\
    u_1(0) = u_{1,0}, \quad u_2(0) &= u_{2,0}. \notag
\end{align}
\end{subequations}
$ \forall \varphi, \psi \in W$, and $t \in (0,T]$.
\noindent Let $\Aop: \W \to \mathcal{L}(\W; \W')$ and $\rop: \W \to \W'$ be linear and nonlinear operators, respectively, which for $\u, \v, \w \in \W$ are defined as
\begin{align}
\label{eq:linopA}
\begin{array}{l l}
{\Aop}(\w) := \begin{bmatrix} 
	 A_1(w_1) & A_2(w_1)  \\
	 0 & A_3(w_2)
\end{bmatrix}, & \rop(\u) := \begin{bmatrix} 
	 r_1(\u) + r_2(\u)  \\
	 r_3(\u)
\end{bmatrix} \\ 
\big(A_1(w_1) u_1, v_1\big) := (\sigma \grad u_1 , \grad v_1) + (\mu u_1, v_1), & (r_1(\u),\v) := (\mu f_1(u_1,u_2) + \mu u_1, v_1), \\
\big(A_2(w_1) u_2, v_1\big) := (\delta w_1 \nabla u_2, \grad v_1 ), & (r_2(\u),\v) := (-u_1 \grad \kappa , \grad v_1),\\
\big(A_3(w_2)u_2,v_2\big) := (\grad u_2, \grad v_2) + (\alpha u_2, v_2), & (r_3(\u),\v) := (\beta f_2(u_1,u_2), v_2).
\end{array}
\end{align}


\noindent
Then, $ \forall \bphi := (\varphi_1, \varphi_2)^{\top} \in \W$, the above weak formulation can be compactly rewritten as
\begin{subequations}
\label{eq:matPseudo}
\begin{align}
(\partial_t u_1, \varphi_1) + (A_1(u_1)u_1 ,\varphi_1) + (A_2(u_1)u_2 ,\varphi_1)  &= (r_1(\u) + r_2(\u),\varphi_1) \quad  t \in (0,T]\\
(\partial_t u_2, \varphi_2) + (A_3(u_2)u_2,\varphi_2)  &= (r_3(\u),\varphi_2) \quad t \in (0,T]\\
\u(0) &= \u_0. \notag
\end{align}
\end{subequations}
Based on the definition of $\Aop$ the corresponding trilinear form $\aop: \W \times \W \times \W \to \R$ can be defined 
in the following way:
\begin{align}
\label{eq:bilinopA}
    \aop(\w)[\u,\v] := \langle A(\w) \u, \v \rangle_{\W', \W} &= \big(A_1(w_1) u_1, v_1\big) + \big(A_2(w_1) u_2, v_1\big) + \big(A_3(w_2)u_2,v_2\big) \notag\\
    &= a_1(w_1)[u_1,v_1] + a_2(w_1)[u_2,v_1] + a_3(w_2)[u_2,v_2] \\
    &= a_1(\w)[\u,\v] + a_2(\w)[\u,\v] + a_3(\w)[\u,\v]. \notag
\end{align}
\noindent
Let the parameters appearing in \eqref{eq:absPseudo} be represented by the vector function $\bthet$ defined as $\bthet := (\th_1, \dots, \th_6)^{\top} = (\sigma, \kappa, \delta, \alpha, \beta, \mu)^{\top}. $ Based on this, we can now formulate the TOCP. First, let the objective (or cost) functional $J$ be defined as
\begin{align}
\label{eq:cost}
    J(\u,\bthet) = {1 \over 2} \|u_1(T) - O\|^2_{L^2} + {\lam  \over 2} \|\bthet\|^2_{V(T)},
\end{align}
 where $\lambda \in \R$ and $O \in \Ho$ is the image data of the observed pseudopalisade pattern. The aim of the TOCP is to find $\u \in \U$ and $\bthet \in \TH$ such that $\u, \bthet$ satisfy the state equation ${\bf G}(\u,\bthet) = 0$ while minimizing the functional $J$.
Altogether, it results in the following minimization problem: find $(\u^*, \bthet^*)$ with $\u^* := \u(\bthet^*)$ and
\begin{align}
\label{eq:pseudOCP}
    \begin{split}
        \u^*, \bthet^* &= \argmin{\u \in \U, \bthet \in \TH} J(\u,\bthet) \quad {\rm{ s.t. }} \quad {\bf G}(\u; \bthet) = 0, \text{ with } \u \in \U \text{ and } \bthet \in \THad,
    \end{split} 
\end{align}
where the equality constraint ${\bf G}(\u; \bthet) = 0$ represents the system dynamics specified by \eqref{eq:weakPseudo}. The mapping ${\bf G}: \U \times \THad \to \U' \times \Z$ reads as
\begin{align}
\label{eq:opConstraint}
    {\bf G}(\u, \bthet) &= \begin{pmatrix}(\partial_t u_1, \cdot) + (\sigma \nabla u_1, \nabla \cdot) + (u_1 \nabla \kappa, \nabla \cdot)  +  (\delta u_1 \nabla u_2, \nabla \cdot)   - (\mu f_1, \cdot)\\
    (\partial_t u_2, \cdot) + (\nabla u_2 , \nabla \cdot) + (\alpha u_2, \cdot) - (\beta f_2, \cdot) \\
    \u(0) - \u_0.
\end{pmatrix}
\end{align}
\noindent
Letting $\| \cdot \|_{\THad} := \| \cdot \|_{{C(\bar I; Z^{\times 6})}}$, the subspace $\THad \subset \TH$ is the set of admissible parameters defined as 
\begin{align*}
    \THad &:= \Big\{ \bthet \in \TH \cap C(\bar I; Z^{\times 6}) :\ \bthet \in C(\bar I; Z^{\times 6}),\  \|\bthet\|_{\THad} \le \ubd_{\bthet} \Big\}. 
\end{align*}
\begin{rem}
It is actually sufficient to define $\THad$ as $$ \THad := \Big \{ \bthet \in \TH \cap C(\bar I; Z^{\times 6}): \th_{i} \in C(\bar I; Z), \|\th_{i}\|_{C(\bar I; Z)} \le \ubd_{\th_i} \le \ubd_{\bthet}, \:\: i \in \{1,2,3\} \Big \}.$$ For the sake of notational simplicity we shall avoid it here, since it does not bring any major difference in the analysis.
\end{rem}

\subsection{Model wellposedness} \label{sec:modelwellposedness} 
In this section we look at the wellposedness of the system and investigate the existence of the optimal solution. For this purpose we first introduce the following general assumptions.
 

\paragraph{Assumptions on $f_1$, $f_2$ and $\th$:}

\begin{enumerate}
\item $f_1, f_2 \in C^{\infty}$, and $f_1(u,v), f_2(u,v) \in W \cap \Linf $ whenever $u,v \in W$, $t\in I$. 
\item $\partial_{u,v} f_1, \partial_{u,v} f_2$  are elements of $W \cap \Linf $ for $u,v \in W$, $t\in I$. Thus $\|f_j\|_{C(\bar I; W)} \le \ubd_{f_j}$ for $j \in \{1, 2\}$. 
\item In particular, for our application we have that $f_1(u,v):= {u(1-u)(1-v) \over (1 + u^2 + v^2)^2}$ and $f_2(u,v) := {uv \over (1 + u^2 + v^2)^2}$, which satisfy the above conditions.
\item The model parameter functions $\th_i \in \{ \sigma, \kappa, \delta, \alpha, \beta, \mu \} $ are such that:\\ $0 < \lbd_{\th_i} \le \th_i(t,x) \le \ubd_{\th_i} < \infty $, for all $t \in I$, $x \in \dom$. Additionally, it is also assumed that $\|\th_i\|_{C(\bar I; Z)} \le \ubd_{\th_i}$. 
\end{enumerate}

\paragraph{Energy estimates of solutions}
\begin{lem} \label{lem:aprEst} Let $\u = (u_1,u_2) \in \Z $ satisfy the equation \eqref{eq:absPseudo}. Then its components fulfill the following energy estimates:
\begin{align}
\label{eq:eEst}
    \|u_1\|^2_{W(T)} \le \const_{u_1} \|u_{1,0}\|^2_{W} \quad \text{ and } \quad \|u_2\|^2_{W(T)} \le \const_{u_2} \|u_{2,0}\|^2_{W},
\end{align}
with $k_u(T)$ and $k_v(T)$ appropriate constants.
\end{lem}

\begin{proof}
Multiplying by $u_2$ both sides of \eqref{eq:He} 
we get
\begin{align*}
{1\over2}{d \over dt} \intdom u_2^2 + \intdom \vert \nabla u_2 \vert^2  {+\alpha \intdom u_2^2}&\le \ubd_{\beta} \ubd_{f_2} \|u_2\|^2_V \Ra \|u_2\|^2_{W(T)} \le \|u_{2,0}\|^2_V \exp(\ubd_{\beta} \ubd_{f_2}T)\\
\Ra \|u_2\|^2_{W(T)} &\le \const_{u_2}(T) \|u_{2,0}\|^2_W \le \ubd_{u_2}(T),
\end{align*}
with sufficiently large $\ubd_{u_2}(T)$.

\noindent
Similarly, multiplying by $u_2'$ both sides of \eqref{eq:He} 
we get
\begin{align}
\label{eq:u2gradBd}
2\|u_2'\|^2_V + {1 \over 2} {d \over dt} \intdom \vert \nabla u_2 \vert^2 {+\alpha \frac{d}{dt}\|u_2\|_V^2}&\le \ubd_{\beta} \ubd_{f_2} {d \over dt} \|u_2\|^2_V  \le \ubd_{\beta} \ubd_{f_2} \const_{u_2}\|u_{2,0}\|^2_{W} \notag \\
\Ra \|u_2'\|^2_{V(T)} + \sup{t\in [0,T]} \|\grad u_2\|^2_{V} &\le (1 + \ubd_{\beta} \ubd_{f_2} \const_{u_2} T)\|u_{2,0}\|^2_{W} \le \ubd_{u_2}(T). 
\end{align}
Lastly, the second derivative can be bounded from above as follows:
\begin{align}
\label{eq:u2lapBd}
    \|\lap u_2\|_V &\le (\ubd_{\alpha} + \ubd_{\beta} \ubd_{f_2} ) \|u_2\|_V + \|u'_2\|_V \notag\\
    &\le  (\ubd_{\alpha} + \ubd_{\beta} \ubd_{f_2} ) [ \|u_2\|_V + \|u'_2\|_V ] \notag\\
    &\le (\ubd_{\alpha} + \ubd_{\beta} \ubd_{f_2} )( 2 \const_{u_2} + 1 + \ubd_{\beta} \ubd_{f_2} \const_{u_2} T) \|u_{2,0}\|_W \notag \\
    &\le \ubd_{u_2}(T).
\end{align}
Altogether, we get that $u_2 \in L^{\infty}([0,T];Z)$ and $u'_2 \in L^2(I;V)$. {{Additionally, differentiating \eqref{eq:He} and multiplying with $u_2'$ and using \eqref{eq:u2energy_ineq} below and boundedness of $\partial_{u_1} f_2, \partial_{u_2} f_2$ we can get that $u_2 \in H^1([0,T];W)$.}}


\noindent Now let us consider equation \eqref{eq:Can}. Like above, multiplying by $u_1$  both sides of the equation 
we get: 
\begin{align*}
{1 \over 2} {d \over dt} \intdom u_1^2 + \ {\lbd_{\sigma}} \intdom \vert \nabla u_1 \vert^2 &\le  (\ubd_{\kappa} + \ubd_{\delta} \|\nabla u_2\|_{\Li}) \intdom u_1 1 \cdot \nabla u_1  + 2 \ubd_{\mu} \ubd_{f_1} \intdom {1\over 2 } u_1^2 \\
&\le {(\ubd_{\kappa} + \ubd_{\delta} \|\nabla u_2\|_{\Li})^2 \over 2 \lbd_{\sigma}} \|u_1\|^2 + {\lbd_{\sigma} \over 2} \|\nabla u_1\|^2  + 2 \ubd_{\mu} \ubd_{f_1} {1 \over 2} \|u_1\|^2 \\
\Ra {1 \over 2} {d \over dt} \intdom u_1^2 + \ {\lbd_{\sigma} \over 2} \intdom \vert \nabla u_1 \vert^2 &\le \Big[{(\ubd_{\kappa} + \ubd_{\delta} \ubd_{u_2})^2 \over 2 \lbd_{\sigma}} + 2 \ubd_{\mu} \ubd_{f_1} \Big] {1 \over 2} \|u_1\|^2  \\ 
\Ra \|u_1\|^2_{V(T)} &\le \exp(\const_{u_1}T) \|u_{1,0}\|^2_{W} \text{ and } \|\nabla u_1\|^2_{L^2} \le \exp(\const_{u_1}T) \|u_{1,0}\|^2_{W} \\
\Ra \|u_1\|^2_{W(T)} &\le \exp(\const_{u_1}T) \|u_{1,0}\|^2_{W}.
\end{align*}
Now multiplying by $u_1'$ both sides of \eqref{eq:Can} we get
\begin{align}
\label{eq:u2energy_ineq}
\|u_1'\|^2 &+ {\lbd_{\sigma} \over 2} {d \over dt} \|\nabla u_1\|^2 \le (\ubdk + \ubdd \ubd_{u_2})  (u_1 1,\grad u_1')  + 2\ubdmu \ubd_{f_1} {1 \over 2} {d \over dt} \|u_1\|^2 + \ubd_{\sigma} \|\nabla u_1\|^2 \notag \\
&\le (\ubdk + \ubdd \ubd_{u_2}) (u_1 1,\grad u_1')  + (2\ubdmu \ubd_{f_1} + \ubd_{\sigma}) \exp(\const_{u_1}T) \|u_{1,0}\|^2_{W}
\end{align}
Using integration by parts for the time derivative i.e. \begin{align*}
\int_0^T \int_{\dom} u_1 1 \cdot \grad \u'_1 &= \int_{\dom} \int_0^T  u_1 1 \cdot \grad \u'_1 = \int_{\dom} (u_1 1 \cdot \grad u_1) \Big \vert_0^T - \int_0^T \int_{\dom} \u'_1 1 \cdot \grad u_1,
\end{align*} 
we get that $(\ubdk + \ubdd \ubd_{u_2}) \int_0^T(u_1 ,\grad u_1')$ satisifes the following inequalities.
\begin{align*}
 &\le (\ubdk + \ubdd \ubd_{u_2}) \Big( \eps \|u_1'\|^2_{V(T)} + {1 \over 4\eps} \|\grad u_1\|^2_{V(T)} + {1 \over 2} (\|u_{1,0}\|^2_{V} + \|\grad u_{1,0}\|^2_{V} + \|u_{1,T}\|^2_{V} + \|\grad u_{1,T}\|^2_{V} ) \Big) \\
&\le {1 \over 2} \|u_1'\|^2_{V(T)} +  {(\ubdk + \ubdd \ubd_{u_2})^2 \over 2} \Big(  \|\grad u_1\|^2_{V(T)} +  {1 \over 2} (\|u_{1,0}\|^2_{V} + \|\grad u_{1,0}\|^2_{V} + \|u_{1,T}\|^2_{V} + \|\grad u_{1,T}\|^2_{V} ) \Big) \\
&\le {1 \over 2} \|u_1'\|^2_{V(T)} +  {(\ubdk + \ubdd \ubd_{u_2})^2 \over 2} \Big(  \|u_1\|^2_{W(T)} +  {1 \over 2} (\|u_{1,0}\|^2_{W} + \|u_{1,T}\|^2_{W} ) \Big).
\end{align*}
Now plugging this in \eqref{eq:u2energy_ineq} we get that
\begin{align*}
    \|u_1'\|^2_{V(T)} + {\lbd_{\sigma}} \sup{t\in [0,T]}  \|\grad u_1\|^2_{V} &\le {(\ubdk + \ubdd \ubd_{u_2})^2} \Big(  \|u_1\|^2_{W(T)} +  {1 \over 2} (\|u_{1,0}\|^2_{W} + \|u_{1,T}\|^2_{W} ) \Big) \\
    &\hspace*{3cm} + (4\ubdmu \ubd_{f_1} + 2\ubd_{\sigma}) \exp(\const_{u_1}T) \|u_{1,0}\|^2_W\\
    &\le \Big(2 (\ubdk + \ubdd \ubd_{u_2})^2 + 4\ubdmu \ubd_{f_1} + 2\ubd_{\sigma} \Big) \exp(\const_{u_1}T) \|u_{1,0}\|^2_W.
\end{align*}
Additionally, if $u_1 \in L^{\infty}([0,T];V)$, we also have that
\begin{align*}
    \lbdsig \|\lap u_1\|_V &\le \ubdsig \|\grad u_1\|_V + \ubdk \|\grad u_1\|_V + \|u_1\|_V \ubdk +  \ubdmu \ubd_{f_1} \|u_1\|_V + \|u_1'\|_V \\
    &\hspace*{2cm} + \ubdd ( \ubd_{u_2} \|u_1\|_V + \ubd_{u_2} \|\grad u_1\|_V + \|u_1\|_{\Li} \|\lap u_2\|_V) \\
    &\le (\ubdsig + \ubdk + \ubdmu \ubd_{f_1} + \ubdd \ubd_{u_2}) \|u_1\|_W +  \ubdd \|u_1\|_{\Li} \|\lap u_2\|_V + \|u_1'\|_V\\
    &\le (\ubdsig + \ubdk + \ubdmu \ubd_{f_1} + \ubdd \ubd_{u_2}) \exp(\const_{u_1}T)\|u_{1,0}\|_V + \ubdd \ubd_{u_2} \|u_1\|_{\Li} + \|u_1'\|_V.
\end{align*}
\end{proof}

\paragraph{Non-negativity of solutions}
\begin{lem}\label{lem:nonneg} Let $\u = (u_1,u_2) \in \Z $ satisfy system  \eqref{eq:absPseudo}. Then $\u(t) \ge 0$ for all $t \in I$ if $\u_0 \ge 0$. 
\end{lem}
\begin{proof}
Let $f_1 := u_1 \fotil$, with $\fotil:= {(1-u_1)(1-u_2) \over {(1 + u_1^2 + u_2^2)^2}}$, 
$q(u_1):= \begin{cases} {1 \over 2} u_1^2 \quad \text{if } u_1 \in (-\infty, 0) \\ 0 \quad \text{else } \end{cases}$. \\
\noindent
Then the function $Q(t) = \int_{\dom} q(u_1(t)) dx$ is continuously differentiable. 
Its derivative (using \eqref{eq:Can}) is given by:
\begin{align*}
    Q'(t) &= \intdom q'(u_1) \nabla \cdot (\sigma \nabla u_1 + u_1 \grad \kappa) + \intdom q'(u_1) \nabla \cdot (\delta \: u_1  \nabla u_2) + \intdom q'(u_1) \mu  u_1 \fotil\\
    &= - \intdom \nabla q'(u_1) \cdot  \sigma \nabla u_1 - \intdom \grad \kappa \cdot u_1 \nabla q'(u_1)   - \intdom \nabla q'(u_1) \cdot  \delta u_1 \nabla u_2 + \mu \intdom q'(u_1) u_1 \fotil \\
    &\le - \intdom \sigma \nabla q'(u_1) \cdot  \nabla u_1 + \ubd_{\kappa} \intdom \vert \nabla q'(u_1) \cdot \mathbf{1} u_1  \vert   +   \ubd_{\delta} \intdom \vert \nabla q'(u_1) \cdot  u_1 \nabla u_2 \vert  + \ubd_{\mu} \ubd_{f_1} \intdom q'(u_1) u_1   \\
    &{\le} - \lbd_{\sigma} \intdom \vert\nabla q'(u_1) \vert  ^2 + \ubd_{\kappa} \intdom \eps_1 \vert \nabla q'(u_1) \vert^2 + {\ubd_{\kappa}\over 4 \eps_1} {\intdom} u_1^2 + \ubd_{\delta} \intdom \eps_2 \vert \nabla q'(u_1) \vert^2\\
    &+ {\ubd_{\delta}\over 4 \eps_2} \int_{\dom} u_1^2 \vert \nabla u_2 \vert^2  
   + \ubd_{\mu} \ubd_{f_1} \intdom q'(u_1) u_1  \quad \text{(using Young's inequality)}\\ 
    & \le {\ubd^2_{\kappa} \over \lbd_{\sigma}}\intdom {1\over 2} u_1^2 + { \ubd^2_{\delta} \ubd^2_{u_2} \over \lbd_{\sigma}} \intdom {1\over 2 } u_1^2 +  \ubd_{\mu} \ubd_{f_1} \intdom {1\over 2} q'(u_1)u_1 \\ 
    &{\le} \Big[{ \ubd^2_{\kappa} \over \lbd_{\sigma}} + { \ubd^2_{\delta} \ubd^2_{u_2} \over \lbd_{\sigma}} +  \ubd_{\mu} \ubd_{f_1} \Big] \intdom  q(u_1)   :=   \const(t) \: Q(t) \\
    { \Ra }\ Q(t) &= 0 \quad \forall t \ge 0 \quad \text{using Gronwall's inequality.}
\end{align*}
Similarly, letting $q(u_2):= \begin{cases} {1 \over 2} u_2^2 \quad \text{if } u_2 \in (-\infty, 0) \\ 0 \quad \text{else } \end{cases}$ and $Q(t) = \int_{\dom} q(u_2(t)) dx$, and using the above result that $u_1 \ge 0$, we get that $Q'(t) \le \ubd_{\alpha} Q(t)$. This in turn (due to Gronwall) implies that $ Q(t) = 0$ for all $t \ge 0$.
    \\
\end{proof}
\begin{lem} 
\label{lem:bdedPseudo}
Let $\u$ be a solution to equation \eqref{eq:absPseudo}. Then $\u \in L^{\infty}(I\times\dom)$.
\end{lem} \begin{proof} In order to prove uniform boundedness of $\u$ in $\domT := [0,T]\times \dom $, following the approach of \cite{Finotti12}, we partition $\domT$ along the time axis and show that the increment of the magnitude of $\u$ over these different partitions tends to zero. 
To this end, let the finite time sequence $(t_i)_{i}$, with $i \in \{0,\dots,K\}$ and $K\in \N$ finite, represent the partition of $[0,T]$. Correspondingly, let $\domti := [t_{i-1},t_i]\times \dom $ represent the $i^{th}$ partition of the time-partitioned space-time cylinder. Let $\u := (u,v)$ and let the norm for the time-continuous $H^1(\dom)$ valued functions be denoted as 
$$\|w\|_{\chti (\mathcal D)} := \sup {t\in[t_{i-1},t_i]} \|w(t)\|_{\hot}.$$
Based on the energy estimate \eqref{eq:eEst} we have that $\|u\|_{\chti(\mathcal D)} \le \const_{u_1}$ and $\|v\|_{\chti (\mathcal D)} \le \const_{u_2}$ for every $i \in \N$.
Now let $u_k = \max\{u-k,0\}$ for $k > \hat k := \max\{1+\eps,\|u_0\|_{\Li}\}$. Correspondingly, the supporting sets of $u_k$ are denoted by $\dom_{u_k}(t) := \{x \in \dom: u(t,x) > k\}$ and $\domti(k) := \{(t,x) \in \domti: u(t,x) > k\}$ for each $i \in \N$.
Now, testing \eqref{eq:Can} with $u_k$ and letting $f_1 := u \fotil$ we get
\begin{align*}
{1 \over 2} {d \over dt} \intdom u_k^2 + \lbd_{\sigma}  \intdom \vert\nabla u_k\vert^2 \le - \intdom u \nabla \kappa \cdot \nabla u_k    &- \intdom \nabla u_k \cdot \delta u \nabla v + \ubd_{\mu} \intdom u_k u \fotil \\
\le \intdom  \vert\nabla u_k \cdot u \grad \kappa \vert &+  \intdom \vert\nabla u_k \cdot \delta u \nabla v\vert  + \ubd_{\mu} \intdom u_k u \fotil\\
= \intdomk  \vert\nabla u_k \cdot \grad \kappa u\vert &+  \intdomk \vert\nabla u_k \cdot \delta u \nabla v\vert  + \ubd_{\mu}  \intdomk u_k u \fotil\\
\hspace*{-2cm}\le \intdomk  \eps_1 \vert\nabla u_k\vert^2  + \intdomk {\ubd_{\kappa}^2 \over 4 \eps_1} u^2 &+ \intdomk \eps_2 \vert\nabla u_k\vert^2 + \intdomk {\ubd_{\delta}^2  \over 4 \eps_2} u^2 \vert\nabla v\vert^2  \\
&\hspace*{1cm} + \ubd_{\mu}\ubd_{f_1} \intdomk u_k u \\
\Ra \ {1 \over 2} {d \over dt} \intdom u_k^2 + {\lbd_{\sigma} \over 2}  \intdom \lvert \nabla u_k \lvert^2 &\le \Big[{\ubd_{\kappa}^2 \over \lbd_{\sigma}} + {\ubd^2_{\delta}\over \lbd_{\sigma}} \ubd_{u_2}^2 \Big]\intdomk u^2 + \ubd_{\mu} \ubd_{f_1} \intdomk u_k u \\
{1 \over 2} {d \over dt} \intdom u_k^2 + {\lbd_{\sigma} \over 2}  \intdom \vert\nabla u_k\vert^2 &\le \Big[{\ubd_{\kappa}^2 \over \lbd_{\sigma}} + {\ubd^2_{\delta}\over \lbd_{\sigma}} \ubd_{u_2} ^2\Big]\intdomk u^2 + \ubd_{\mu} \ubd_{f_1} \intdomk u_k u \\
&\hspace*{-2cm}\le \Big[{\ubd_{\kappa}^2 \over \lbd_{\sigma}} + {\ubd^2_{\delta}\over \lbd_{\sigma}} \ubd_{u_2}^2  + \ubd_{\mu} \ubd_{f_1} \Big] \intdomk (u^2 + u u_k) \\
&\hspace*{-2cm}\le 4\Big[{\ubd_{\kappa}^2 + \ubd^2_{\delta} \ubd_{u_2}^2 \over \lbd_{\sigma}} + \ubd_{\mu} \ubd_{f_1} \Big] \intdomk ((u_k)^2 + k^2) \\
&\hspace*{-2cm}= \const(\lbd_{\sigma},\ubd_{\kappa},\ubd_{u_2},\ubd_{\delta},\ubd_{\mu},\ubd_{f_1})(t) \|u_k\|^2_{L^2} + k^2 \lvert\dom_{u_k}(t)\rvert.
\end{align*}
Integrating w.r.t. $t \in [0,t_1]$, with $t_1 > 0$ small enough such that $$t_1 
\sup{t\in[0,t_1]}\const(\lbd_{\sigma},\ubd_{\kappa},\ubd_{u_2},\ubd_{\delta},\ubd_{\mu},\ubd_{f_1})(t) < 1/2,$$ 
we get that
\begin{align}
\label{eq:ukchtoIneq}
\|u_k\|^2_{\chto} &\le 2 \const(\lbd_{\sigma},\ubd_{\kappa},\ubd_{u_2},\ubd_{\delta},\ubd_{\mu},\ubd_{f_1}) k^2 \eta_k, \text{ with } \eta_k = \lvert \domto \lvert = \int_0^{t_1}  \lvert\dom_k(t) \lvert dt \notag \\
\Ra \|u_k\|_{\chto} &\le \vartheta_1 \: k \: \eta_k^{1\over 2}.
\end{align}
Let $N_0 := n_0 \hat k$ for some $n_0 > 1$, let $k_i = N_0 ( 2 - 2^{-i})$ for $i \in \N_0$, then $\eta_k$ fulfills the following inequality: 
\begin{align}
\label{eq:delIneq}
    k_{i+1} \eta^{1\over r}_{k_{i+1}} &= k_{i+1} \Bigg(\int_0^{t_1} \vert\dom_{k_{i+1}}(t)\vert dt\Bigg)^{1 \over r} = k_{i+1} \Bigg(\int_0^{t_1} \int_{\dom_{k_{i+1}}(t)} dx \: dt \Bigg)^{1 \over r} \notag \\
    &= [k_i + (k_{i+1}-k_i)] \Bigg(\int_0^{t_1} \int_{\dom_{k_{i+1}}(t)} dx \: dt \Bigg)^{1 \over r} \\
    &\le [k_i + (k_{i+1}-k_i)] \Bigg(\int_0^{t_1} \int_{\dom_{k_i}(t)} dx \: dt \Bigg)^{1 \over r} \notag \\
\Ra (k_{i+1} - {k_i}) \eta^{1 \over r}_{k_{i+1}}   &\le \Bigg(\int_0^{t_1} \int_{\dom_{k_i}(t)} (k_{i+1}-k_i)^r dx \: dt \Bigg)^{1 \over r} \quad (\text{since } k_i > \hat k > 1 \text{ are constants}) \notag\\
    &= \Bigg(\int_0^{t_1} \int_{\dom_{k_i}(t)} (u_{k_i}-u_{k_{i+1}})^r dx \: dt  \Bigg)^{1 \over r}  \notag\\
    \le \Bigg(\int_0^{t_1} \Big(\int_{\dom_{k_i}(t)} &(u_{k_i}-u_{k_{i+1}})^r dx  + \int_{\dom_{k_i}^c(t)} u_{k_{i}}^r dx \Big) \: dt  \Bigg)^{1 \over r}  \quad \text{ with $\dom^c_{k_i}(t) := \dom - \dom_{k_i}(t)$}\notag \\
    &\le \Bigg(\int_0^{t_1} \Big(\int_{\dom_{k_i}(t)} u_{k_i}^r dx  + \int_{\dom_{k_i}^c(t)} u_{k_{i}}^r dx \Big) \: dt  \Bigg)^{1 \over r}  = \Bigg(\int_0^{t_1} \int_{\dom} u_{k_i}^r dx \: dt  \Bigg)^{1 \over r}  \notag\\
\Ra (k_{i+1} - {k_i}) \eta^{1 \over r}_{k_{i+1}} &\le \|u_{k_i}\|_{L_r(\domto)}.
\end{align}
Using the Sobolev embedding inequality for $2 (\nu + 1) =: r \in [2, {2d \over d-2}]$ with $\nu > 0$ we have that
\begin{align*}
    \|u_{k_i}\|_{L_r(\domto)} \le \gamma_0 \|u_{k_i}\|_{\hot} &\le \gamma_0 \vartheta_1 k_i \eta_{k_i}^{1\over2}, \quad \forall i \in \N_0, \gamma_0 > 1\\
\Ra \eta^{1 \over r}_{k_{i+1}} \le {\gamma_0 \vartheta_1 k_i \over (k^{i+1} - {k^i})} \eta^{1 \over r}_{k_{i+1}} &\le 4 \gamma_0 \vartheta_1 2^i \eta_{k_i}^{1 \over 2} \le 4 \gamma_0 \vartheta_1 2^i (\eta_{k_i}^{1 \over r})^{1+\nu}, \quad \forall i \in \N_0,\ (\text{since }  {\nu + 1 \over r} = {1 \over 2}).
\end{align*}
In particular, by taking $k_i := \hat k$ and $k_{i+1} := N_0$ and defining $\eta_{_{_{_0}}} := \eta_{_{_{k_0}}}$ we have that $\eta_{0}^{1\over r}   \le {\gamma_0 \vartheta_1 \over n_0 - 1} (T\vert\dom\vert)^{1\over2}$. Thus for $n_0 > 1 + \gamma_0 \vartheta_1 (T\vert\dom\vert)^{1\over2} (4\gamma_0\vartheta_1)^{1\over\nu} 2^{1\over\nu^2}$ we get that $\eta_{0}^{1\over r} \le (4\gamma_0\vartheta_1)^{-1\over\nu} 2^{-1\over\nu^2}.$ Thus invoking Theorem 2.4.1 of \cite{Zacher2010} we get that
$(\eta_{k_i})^{1\over r} \to 0 $ as $i \to \infty$. In particular, we get that $\eta_{_{_{2N_0}}} := \eta_{k_{\infty}} = 0$. Consequently, we get that $u \le c_1 := k_{\infty} = 2 n_0 \hat k$ on $\dom_{t_1}$. Now taking $N \in \N$ such that $\cup_{i=1}^N [t_{i-1},t_i] = [0,T]$ we get that $u \le \sum_{i=1}^N c_i < \infty$ on $\domT$. 
For the $v$ component of $\u$, by repeating the above steps, we also get that $v \le \sum_{i=1}^N c_{v,i}$ in $\domT$. Consequently, we get that $\u$ is bounded in $\domT$.
\end{proof}

\paragraph{Properties of the operators $\aop$ and $\rop$}
\begin{lem}
\label{lem:propAopFop}
Let $\w, \u, \v \in \W$ be non-negative, $\bthet \in \THad \subset \TH$. Then there exist $\ubd_{\aop}, \ubd_{\rop}, \lbd_{\aop} \in (0,\infty)$ such that operators $\aop$ and $\rop$ appearing in \eqref{eq:linopA} satisfy the following inequalities:
\begin{align*}
    \vert {\aop}[\w](\u,\v)\vert &\le \ubd_{\aop} \|\w\|_W \|\u\|_W \|\v\|_W  \\
    \aop(\w)[\u,\u] &\ge \lbd_{\aop}\|\u\|^2_{\W} \\
    \vert \aop(\w^1)[\u,\v] - \aop(\w^2)[\u,\v] \vert &\le \ubd_{\aop}  \| \w^1 - \w^2 \|_{V} \|u_2\|_{W} \|v_2\|_{W} \\ 
    \|\rop(\u)\|_{\W'} &\le \ubd_{\rop} \|\u\|_{\W} \\
    \|\rop(\u^1) - \rop(\u^2)\|_{\W'} &\le \ubd_{\rop} \|\u^1 - \u^2\|_{W}.
\end{align*}
\end{lem}
\begin{proof} First let us recall the definition of $\aop$ and $\rop$ (see \eqref{eq:linopA}, \eqref{eq:bilinopA}). For any $\u, \v, \w \in \W$ they are defined as \\[-5ex]
\begin{align*}
\begin{array}{l l}
{\aop}(\w)[\u,\v] := a_1(\w)[\u,\v] + a_2(\w)[\u,\v] + a_3(\w)[\u,\v], & \quad \rop(\u) := \begin{bmatrix} 
	 r_1 + r_2  \\
	 r_3
\end{bmatrix}, \\
a_1(\w)[\u,\v] := (\sigma \grad u_1 , \grad v_1) + (\mu u_1, v_1), & r_1(\u)[\v] := (\mu f_1(u_1,u_2) + \mu u_1, v_1), \\
a_2(\w)[\u,\v] := (\delta w_1 \nabla u_2, \grad v_1 ), & r_2(\u)[\v] := (-u_1 \grad \kappa , \grad v_1),\\
a_3(\w)[\u,\v] := (\grad u_2, \grad v_2) + (\alpha u_2, v_2), & r_3(\u)[\v] := (\beta f_2(u_1,u_2), v_2)
\end{array}
\end{align*}
Based on this, now let us verify the properties of the operator $\aop: \W \times \W \times \W \to \R$. For the boundedness property we have that: \\[-2ex]
\begin{minipage}{.5\linewidth}
\begin{align*}
    a_1(\w)[\u,\v] &= (\sigma \grad u_1 , \grad v_1) + (\mu u_1, v_1) \\
    &\le \ubdsig \|\grad u_1\| \|\grad v_1\| + \ubdmu \|u_1\|\|v_1\| \\
    &\le \const_{a_1}(\sigma,\mu) \|u_1\|_W  \|v_1\|_W
\end{align*}
\end{minipage}
\begin{minipage}{.5\linewidth}
\begin{align*}
    a_3(\w)[\u,\v] &= (\grad u_2 , \grad v_2) + (\alpha u_1, v_1) \\
    &\le \|\grad u_2\| \|\grad v_2\| + \ubdalp \|u_2\|\|v_2\| \\
    &\le \const_{a_3}(\alpha) \|u_2\|_W \|v_2\|_W
\end{align*}
\end{minipage}
\begin{minipage}{.5\linewidth}
\begin{align*}
    a_2(\w)[\u,\v] &= (\delta w_1 \grad u_2 , \grad v_2)  \\ 
    &\le \ubdd \|w_1\|_{L^6} \| \grad u_2\|_{L^3} \|\grad v_2\|_{L^2} \\ 
    &\le \const_{a_2}(\delta)  \|w_1\|_W \| u_2\|_W \|v_2\|_W
\end{align*}
\end{minipage}
\begin{minipage}{.5\linewidth}
\begin{align*}
   \Ra \vert \aop(\w)[\u,\v] \vert &\le \const_{a_1} \|u_1\|_W \|v_1\|_W + \const_{a_2}\|u_2\|_W \|v_2\|_W \\
   &\hspace*{.5cm} + \const_{a_3}\|u_2\|_W \|v_2\|_W \\
   &\le \ubd_{\aop}  \|\w\|_{\W} \| \u\|_{\W} \|\v\|_{\W}
\end{align*}
\end{minipage}\\[2ex]
with $0 < \ubd_{\aop} < \infty$ being a common big enough constant independent of any $\u,\v,\w \in \W$. Next, for the coercivity property let $\w = (w_1, w_2)^{\top}$ be such that $w_1 \ge 0$, then for any $\u \in \W$ we have that: \\[-6ex]
\begin{align*}
    \aop(\w)[\u,\u] &= a_2(\w)[\u,\u] + a_2(\w)[\u,\u] + a_3(\w)[\u,\u] \\
    &= (\sigma \grad u_1 , \grad u_1) + (\mu u_1, u_1) +  (\grad u_2 , \grad u_2) + (\alpha u_1, u_1) + (\delta w_1 \grad u_2 , \grad u_2) \\
    &= (\sigma, |\grad u_1|^2) + \|\grad u_2\|^2 +  \|\mu u_1\|^2 +  \|\alpha u_2\|^2 + (\delta w_1, |\grad u_2|^2) \\
    &\ge \lbdsig (1, |\grad u_1|^2) + \|\grad u_2\|^2 +  \lbdmu \| u_1\|^2 +  \lbdalp \| u_2\|^2 + \lbdd(w_1, |\grad u_2|^2) \\
    &\ge \const(\sigma,\alpha,\mu, \delta) \Big(\|\grad u_1\|^2 + \|\grad u_2\|^2 + \|u_1\|^2 +  \|u_2\|^2 + (w_1, |\grad u_2|^2) \Big)\\
    &\ge \const(\sigma,\alpha,\mu, \delta) (\|u_1\|^2_W + \|u_2\|^2_W)  \quad (\text{since } w_1 \ge 0) \\
\Ra a(\w)[\u,\u] &\ge \lbd_{\aop} \|\u\|^2_{\W}
\end{align*}
with $0 < \lbd_{\aop} < \infty$ is a small enough constant independent of any $\u, \w \in \W$. Next, for the Lipschitz continuity property let $\w^1, \w^2 \in \W$, then we have that: 
\begin{align*}
\vert\aop(\w^1)[\u,\v] - \aop(\w^2)[\u,\v]\vert &= \vert a_2(\w^1)[\u,\v] - a_2(\w^2)[\u,\v]\vert = \vert (\w^1 - \w^2  \grad u_2, \grad v_2 ) \vert \\
&\le \ubd_{\delta} \| \w^1 - \w^2 \|_{L^2} \|\grad u_2\|_{L^3} \|\grad v_2\|_{L^6} \\
&\le \ubd_{\delta}  \| \w^1 - \w^2 \|_{\V} \|u_2\|_{\W} \|v_2\|_{\W} \\
&\le \ubd_{\aop}  \| \w^1 - \w^2 \|_{\V} \|u_2\|_{\W} \|v_2\|_{\W} \\
\Ra \vert\aop(\w^1)[\u,\v] - \aop(\w^2)[\u,\v]\vert &\le \ubd_{\aop} \|\w^1 - \w^2\|_{\W}
\end{align*}
with $0 < \ubd_{\aop} < \infty$ beign a large enough generic constant independent of any $\u,\v,\w \in \W$. \\
\noindent
Now, for  the nonlinear operator $\rop: \W \to \W'$ let $\u, \v \in \W$. Then we have that
\begin{align*}
    \vert\rop(\u)[\v]\vert  &\le \vert r_1(\u)[\v]\vert + \vert r_3(\u)[\v] \vert + \vert r_2(\u)[\v]\vert  \le \vert(\mu f_1(u_1,u_2) + \mu u_1, v_1) \vert  + \vert\beta (f_2(u_1,u_2), v_2) \vert \\
    &\hspace*{2cm} + \vert(-u_1 \grad \kappa , \grad v_1) \vert \\
    &\le \ubd_{\mu} (1 + \ubd_{f_1}) \|u_1\|_{W} \|v_1\|_{W}+ \ubd_{\beta} \ubd_{f_2} \|u_2\|_{W} \|v_2\|_{W} + \ubd_{\kappa} \|u_1\|_{W} \|v_1\|_{W} \\
    &\le \const(\kappa, \mu, \alpha, \beta, f_1, f_2) \: \|\u\|_{\W} \: \|\v\|_{\W} \\
\Ra  \|\rop(\u)\|_{\W'}  &\le \ubd_{\rop} \|\u\|_{\W}
\end{align*}
with $0 < \ubd_{\rop} < \infty$ begin a big enough constant independent of any $\u,\v \in \W$.\\
Next for the Lipschitz continuity, let $\u^1, \u^2 \in \W$ then we have that
\begin{align*}
     \vert(\rop(\u^1) &- \rop(\u^2)) [\v]\vert  \le \vert (r_1(\u^1) - r_1(\u^2)) [ v_1]\vert + \vert (r_3(\u^1) - r_3(\u^2)) [v_2] \vert + \vert (r_2(\u) - r_2(\v)) [v_1]\vert \\
    &\le (\ubd_{\mu} \|\nabla f_1\|_{\Li} + \ubd_{\beta} \|\nabla f_2\|_{\Li})  \|\u^1 - \u^2\|_{W} \|\v\|_{\W}  + \ubd_{\mu} \|\u^1_1 - \u^2_1\|_{W} \|v_1\|_{W}  \\
    &\hspace*{2cm} + \ubd_{\kappa} \|\u^2_2 - \u^1_2\|_{W} \|v_1\|_{W} \\
    &\le \const(\kappa, \mu, \alpha, \beta, f_1, f_2) \|\u^1 - \u^2\|_{\W} \|\v\|_{\W} \\
\Ra  \|\rop(\u^1) &- \rop(\u^2)\|_{\W'}  \le \ubd_{\rop} \|\u^1 - \u^2\|_{\W}
\end{align*}
with $0 < \ubd_{\rop} < \infty$ being a large enough generic constant independent of any $\u^1,\u^2 \in \W$.
\end{proof}

\paragraph{Existence of a solution}

\begin{thm} 
\label{thm:existPseudo}
For every $0 \le \u_0  \in \Z$ and $\bthet \in \TH$, the pseudopalisade system \eqref{eq:absPseudo} possesses a unique non-negative weak-solution
\begin{align*}
    \u \in \U \cap \Y \cap \Csp(T), \quad \u \ge 0,  
\end{align*}
where $T > 0$ is dependent on $\|\u_0\|_{\Z} \le \ubd_{\u_0}$ and $\|\bthet\|_{\THad} \le \ubd_{\bthet}$. Moreover, $\u$ satisfies the following inequality:
\begin{align*}
\|\u\|_{\Y} + \|\u\|_{\Csp(T)} + \|\u\|_{\X} \le \const_{\u}({\ubd_{\u_0},\ubd_{\bthet}}).
\end{align*}
\end{thm}
\begin{proof}
This is a direct consequence of Lemmas \ref{lem:aprEst}, \ref{lem:nonneg}, \ref{lem:bdedPseudo}, \ref{lem:propAopFop} and Theorem 5.10 in \cite{Yagi09}.
\end{proof}

\subsection{Existence of optimal parameter functions} \label{sec:ocpwellposedness}
Let us recall that the state equation $G(\u, \bthet) = 0$ (the weak form of equation \eqref{eq:absPseudo}) is given as:  
\begin{align*}
    G(\u, \bthet) &= \begin{pmatrix}(\partial_t u_1, \cdot) + (\sigma \nabla u_1, \nabla \cdot) + (u_1 \nabla \kappa, \nabla \cdot)  +  (\delta u_1 \nabla u_2, \nabla \cdot)   - (\mu f_1, \cdot)\\
    (\partial_t u_2, \cdot) + (\nabla u_2 , \nabla \cdot) + (\alpha u_2, \cdot) - (\beta f_2, \cdot) \\
    \u(0) - \u_0
\end{pmatrix},\ \bthet \in \THad.
\end{align*}
\begin{lem}
\label{lem:OpDiffx}
For a fixed $\bthet \in \THad$, the operator $G: \U \times \THad \to \U' \times \Z$ is infinitely Fr\'echet differentiable. The partial derivative $G_{\u} := \partial_{\u} G$ of $G$ with respect to $\u$ at point $(\u,\th)$ is represented as a mapping $(\u, \bthet) \mapsto G_{\u}[\u,\bthet]$, $G_{\u}: \U \times \THad \to L(\U ; \U' \times \Z)$. Its exact form is given by
\begin{align}
\label{eq:diffU}
&G_{\u}[\u,\bthet](v_1,v_2) =  \notag \\
&\hspace*{-1cm}\begin{pmatrix} (\partial_t v_1, \cdot) + (\sigma \nabla v_1, \nabla \cdot) +  (v_1 \nabla \kappa, \nabla \cdot)  +  (\delta v_1 \nabla u_2, \nabla \cdot)   - (\mu \partial_{_{u_1}} f_1 v_1, \cdot) + (\delta u_1 \nabla v_2, \nabla \cdot)   - (\mu \partial_{_{u_2}}  f_1 v_2, \cdot)\\
(\partial_t v_2 , \cdot) + (\nabla v_2 , \nabla \cdot) + (\alpha v_2, \cdot) - (\partial_{_{u_1}}  f_2 v_1, \cdot) - (\beta \partial_{_{u_2}} f_2 v_2, \cdot) \\
v_1(0) \\
v_2(0)
\end{pmatrix}.
\end{align}
Similarly, the partial derivative $G_{\th} := \partial_{\th} G$ of $G$ with respect to $\bthet$ at the point $(\u,\bthet)$ is represented as a mapping $(\u, \bthet) \mapsto G_{\bthet}[\u,\bthet]$, $G_{\bthet}: \U \times \THad \to L(\THad; \U' \times \Z)$. Its exact form is given as
\begin{align}
\label{eq:diffTh}
G_{\bthet}[\u,\bthet](\phi_1, \dots, \phi_6) = \begin{pmatrix} (\phi_1 \nabla u_1, \nabla \cdot ) + (u_1 \nabla \phi_2, \nabla \cdot ) + (\phi_3 \: u_1  \nabla u_2, \nabla \cdot) + (\phi_4  f_1(\u). \cdot)\\
(\phi_5 u_2, \cdot) + (-\phi_6 f_2(\u), \cdot) \\
0 \\
0 \\
\end{pmatrix}.
\end{align}
\end{lem}
\begin{proof} By following the technique of Lemma 1.17 of \cite{Hinze08} and by using Equation (1.95) of \cite{Yagi09} (due to Lemma \ref{lem:bdedPseudo}) we can obtain the first partial Fr\'echet derivatives $G_{\u}$ and $G_{\th}$ as given in \eqref{eq:diffU} and \eqref{eq:diffTh}, respectively. For the latter, it is clear that the higher-order derivatives $\partial_{\bthet}^{(k)} G$, $k > 1$ are equal to 0.
However, the second derivative $\partial_{\u}^{(2)} G$ reads
\begin{align*}
&G_{\u\u}[\u,\th](v_1,v_2)(v_3,v_4) = \\
&\hspace*{-1cm}\begin{pmatrix}  (\delta v_1 \nabla v_4, \nabla \cdot)   + (\delta v_3 \nabla v_2, \nabla \cdot) - (\mu \partial^2_{_{u_1,u_1}} f_1 v_1 v_3, \cdot) - (\mu \partial^2_{_{u_2,u_1}} f_1 [v_1 v_4+v_2v_3], \cdot)    - (\mu \partial_{_{u_2,u_2}}  f_1 v_2 v_4, \cdot)\\
- (\beta \partial^2_{_{u_1,u_1}} f_2 v_1 v_3, \cdot) - (\beta \partial^2_{_{u_2,u_1}} f_2 [v_1 v_4 +v_2v_3], \cdot) - (\beta \partial_{_{u_2,u_2}}  f_2 v_2 v_4, \cdot) \\
0 \\
0
\end{pmatrix}.
\end{align*}
Now, further derivatives of $G$ with respect to $\u$ are just multiples of the  partial derivatives $\partial_{\u}^{(k)} f_1$ and $\partial_{\u}^{(k)} f_2$. Since $f_1$ and $f_2$ are $C^{\infty}$ functions, we get that $G$ is infinitely Fr\'echet differentiable. 
\end{proof}

\begin{lem}
\label{lem:invGx}
The operator $G_{\u} \in L(\U \times \TH; \U' \times \Z)$ has a bounded inverse.
\end{lem}
\begin{proof} Let $\U' \ni \g = (g_1, g_2)$ and $\Z \ni \v_0$, let $G_{\u}$ be the operator $G_{\u}$ linearized at $\u \in \U$, $\bthet \in \TH$. Then  $G_{\u}$ is said to have an inverse if there exists a unique $\v = (v_1,v_2)$ which satisfies the equation $G_{\u}(\v) = \g$. This is to say that $\v$ satisfies the following PDE in weak form:
\begin{align*}
    (\partial_t v_1, \cdot) + (\sigma \nabla v_1 + \grad \kappa v_1, \nabla \cdot ) + (\delta \: v_1  \nabla u_2, \nabla \cdot) + (\delta \: u_1  \nabla v_2, \nabla \cdot) &= (g_1 + \mu  \partial_{u_1} f_1 v_1 + \mu  \partial_{u_2} f_1 v_2, \cdot)\\
    (\partial_t v_2, \cdot) + (\grad v_2, \grad \cdot) + (\alpha v_2, \cdot) &= (g_2 + \beta \partial_{u_2} f_2 v_2 + \partial_{u_1} f_2 v_1, \cdot) \\
    \v(0) &= \v_0,
\end{align*}
which can be rewritten as 
\begin{align*}
    (\partial_t v_1, \cdot) + (\sigma \nabla v_1 + v_1[\grad \kappa + \delta \grad u_2], \nabla \cdot ) + (\delta \: u_1  \nabla v_2, \nabla \cdot) &= (g_1 + \mu  \partial_{u_1} f_1 v_1 + \mu  \partial_{u_2} f_1 v_2, \cdot)\\
    (\partial_t v_2, \cdot) + (\grad v_2, \grad \cdot) + (\alpha v_2, \cdot) &= (g_2 + \beta \partial_{u_2} f_2 v_2 + \partial_{u_1} f_2 v_1, \cdot) \\
    \v(0) &= \v_0.
\end{align*}
This can be further simplified to
\begin{align*}
    (\partial_t v_1, \cdot) + (\sigma \nabla v_1 + v_1 \tilde \kappa, \nabla \cdot ) + (\delta \: u_1  \nabla v_2, \nabla \cdot) &= (\tilde g_1(\u,\v,\bthet), \cdot)\\
    (\partial_t v_2, \cdot) + (\grad v_2, \grad \cdot) + (\alpha v_2, \cdot) &= (\tilde g_2(\u,\v,\bthet), \cdot) \\
    \v(0) &= \v_0,
\end{align*}
where $\tilde \kappa := \grad \kappa + \delta \grad u_2 $, $\tilde g_1(\u,\v,\bthet) := g_1 +\mu  \partial_{u_1} f_1 v_1 + \mu  \partial_{u_2} f_1 v_2$ and $\tilde g_2(\u,\v,\bthet) := g_2 + \beta \partial_{u_2} f_2 v_2 + \beta \partial_{u_1} f_2 v_1$. Since $\u \in \U$ (cf. Theorem \ref{thm:existPseudo}) and since $\fop$, $\grad \fop$ are elements of $\Z \cap \Y$ for all $t \in I$, Lemmas \ref{lem:propAopFop}, \ref{lem:bdedPseudo} can be used to invoke Theorem 4.7 of \cite{Yagi09} for obtaining  the existence of a unique $\v$ that satisfies the equation $G_{\u}[\u,\bthet](\v) = \g$. Specifically, for $\g \in \U' \cap V(T)$, performing computations similar to Lemma \ref{lem:aprEst}, it follows that the solution $\v \in \X \cap \Y$ satisfies the following inequality:
\begin{align}
\label{eq:ineq_Guinv}
    \|\v\|_{W(T)} &\le \const_{G_{\u}}(\u,\bthet) (\|\v_0\|_{V} + \|\g\|_{V(T)}) \notag\\
    &\le \const_{G_{\u}} (\|\v_0\|_{Z} + \|\g\|_{V(T)}) \notag \\
\Ra    \|G_{\u}^{-1} \|_{L(\U',\U)} &\le \const_{G_{\u}}(\u,\bthet).
\end{align}

\end{proof}

\noindent As a consequence of this, we can apply the implicit function theorem to consider $\u$ as a function of $\bthet$ via the mapping $\bthet \mapsto \u(\bthet)$. Moreover, it also follows that $\u(\bthet)$ is infinitely differentiable, and thus also Lipschitz, with respect to $\bthet$. 
Now we are ready to establish the existence of an optimal parameter $\bthet$. 
\begin{thm}
Let $T \in (0,\infty)$,  $\u_0 \in \Z$ and $\u_0$ be non-negative. Then there exists an optimal parameter vector $\bthet \in \THad$ minimizing the functional $J(\bthet)$ from \eqref{eq:cost}.
\end{thm}
\begin{proof}
Clearly, by definition of $J(\bthet)$, we have that $J \ge 0$ thus $\inf{\bthet} J(\bthet)$ exists. Due to continuity of the mapping $\bthet \mapsto J(\bthet)$ and closedness of $\TH$, there exists a sequence $(\bthet_n)_{n \in \N} \subset \THad$ such that  $J(\bthet_n) \to J(\bthet^*) = \inf{\bthet} J(\bthet)$. Since $\THad \subset \TH$ is bounded, the sequence $(\bthet_n)_{n \in \N}$ is also bounded in $\TH$. Hence, we can extract a weakly convergent sub-sequence $(\bthet_j)_{j \in \N}$  such that $\bthet_j \weakly \bthet^*$ in $\TH$. Now, since the mapping $\THad \ni \bthet \mapsto \u(\bthet) \in \X$ is a bounded operator (due to Lemma \ref{lem:bdedPseudo} and Theorem  \ref{thm:existPseudo}) we can let $(\u_j)_{j\in\N}$ be a sequence of solutions corresponding to the parameter sequence $(\bthet_j)_{j\in\N}$. Since, due to Lemma \ref{lem:bdedPseudo}, the sequence $(\u_j)_{j\in\N}$ is bounded in $\X$, there exists a sub-sequence $(\u_{j_k})_{k\in\N}$ such that $\u_{j_k} \rightharpoonup  \u^*$ in $W(T)$. Moreover, we also have that $(\u'_{j_k})_{k \in \N} \in W(T)'$ is a bounded sequence. Now, applying the Lions-Aubin compactness theorem (see e.g., \cite{Necas19}) we get that:
\begin{align*}
    \u_{j_k} \to \u^*  \text{ in } V(T),& \quad \nabla \u_{j_k} \rightharpoonup \nabla \u^* \text{ in } V(T), \text{ and } \u'_{j_k} \rightharpoonup {\u'}^* \text{ in } W(T)' .
\end{align*}
Due to uniqueness, we get that $\u^*$ is the solution to \eqref{eq:absPseudo} corresponding to $\bthet^*$, i.e. $\bthet^* \mapsto \u^*(\bthet^*)$. Finally, due to the weak lower semicontinuity of the norm and the weak convergence of $\u(\bthet_{j_k})$ to $\u^*(\bthet^*)$ in $V(T)$, we have that $$ J(\u^*(\bthet^*),\bthet^*) = J(\u^*,\bthet^*) \le \liminf{k\to \infty} J(\u_{j_k},\bthet_{j_k}) = \inf{\bthet} J(\u(\bthet),\bthet) \le J(\u^*(\bthet^*),\bthet^*).$$ 
Thus $J(\u(\bthet^*),\bthet^*)$  is indeed equal to $\inf{\bthet} J(\u(\bthet),\bthet).$
\end{proof}

\noindent Next we shall construct the minimizing sequence by deducing the adjoint equation and the necessary optimality condition. 

\newcommand{\proj}{\mathcal{P}}
\newcommand{\projthad}{\proj_{\THad}}

\begin{thm} Let $(\bar \u,\bar \bthet)$ be an optimal solution to problem $\eqref{eq:pseudOCP}$. Then there exists an adjoint state $\bar \p \in \U \subset \U''$ s.t. the following optimality condition holds
\begin{subequations}
\label{eq:optSys}
\begin{align}
    G(\bar \u,\bar \bthet) &= 0  \label{eq:optSysState}\\
    G^*_{\u}(\bar \u,\bar \bthet) \bar \p &= - J_{\u}(\bar \u,\bar \bthet)  \label{eq:optSysAdjoint}\\
    (J_{\bthet}(\bar \u,\bar \bthet) + G^*_{\bthet}(\bar \u,\bar \bthet) \bar \p, \bthet - \bar \bthet)_{\TH} &\ge 0 \label{eq:optSysCond}
\end{align}
\end{subequations}
\end{thm}
\begin{proof}
Due to Theorem \ref{thm:existPseudo} and Lemma \ref{lem:invGx} we can invoke Theorem 1.48 of \cite{Hinze08} for the reduced cost function $\hat J(\bthet) := J(\u(\bthet),\bthet)$, which gives us that the local solution $\bar \bthet \in \THad$ satisfies the following variational inequality:
\[ \langle \hat J'(\bar \bthet), \bthet - \bar \bthet \rangle_{\Th',\Th} \ge 0  \quad \forall \bthet \in \THad.\]
For $\bvarth := \bthet - \bar \bthet, \bthet \in \THad$ we have that 
\begin{align*}
    \langle \hat J'(\bar \bthet), \bvarth \rangle_{\Th',\Th} &= \langle J_{\u}(\bar \u, \bar \bthet), \bar \u' \bvarth \rangle_{\U',\U} + \langle J_{\bthet}(\bar \u, \bar \bthet), \bvarth \rangle_{\Th',\Th}\\
    &= \langle (\bar \u')^* J_{\u}(\bar \u, \bar \bthet),  \bvarth \rangle_{\Th',\Th} + \langle J_{\bthet}(\bar \u, \bar \bthet), \bvarth \rangle_{\Th',\Th}
\end{align*}
Based on the state equation $\hat G(\bthet) := G(\u(\bthet),\bthet) = 0$, we have that $G_{\u} \u'(\bthet) + G_{\bthet} = 0$. This implies $\u'(\bthet) = -G_{\u}^{-1} G_{\bthet}$, Consequently, we get that $\u'(\bthet)^* J_{\u} = -G_{\bthet}^{*} (G_{\u}^*)^{-1} J_{\u}$. Defining $\p := -(G_{\u}^*)^{-1} J_{\u}$ we get the following adjoint equation: 
\begin{align}
\label{eq:preAdjEqn}
    G_{\u}^* \p = - J_{\u}  \text{ in } \U'.
\end{align}
Explicitly, the equation can be written as 
\hspace*{-1cm}
\begin{align}
\label{eq:adjoint}
     \begin{pmatrix} -\partial_t - \nabla \cdot (\mathbf{\sigma} \nabla ) + (\nabla \kappa   + \delta \nabla u_2) \cdot \nabla   &  -\nabla \cdot (\delta u_1 \nabla )\\
      0 & -\partial_t -\lap + \alpha 
    \end{pmatrix} \begin{pmatrix} p_1 \\ p_2 \end{pmatrix} &=   \begin{pmatrix} \mu \partial_{u_1} f_1(\u) &  \mu \partial_{u_2} f_1(\u)  \\ \beta  \partial_{u_1} f_2(\u) &  \beta \partial_{u_2} f_2(\u) \end{pmatrix} \begin{pmatrix} p_1 \\ p_2 \end{pmatrix} \notag\\
    \nabla_n p_1 = 0, \quad   \nabla_n p_2 = 0, \quad p_1(T) = u_1(T) - O,& \quad p_2(T) = 0.
\end{align}
\begin{align}
\label{eq:weakAdjoint}
    \hspace*{-1cm}-(\partial_t p_1, \varphi) + (\sigma \nabla p_1, \nabla \varphi) + ((\grad \kappa + \delta \nabla u_1) \cdot \nabla p_1, \varphi) + (\delta u_1 \nabla p_2, \varphi) &= -(\mu \partial_{u_1} f_1 p_1, \varphi) - (\nu \partial_{u_2} f_1 p_2, \varphi), \notag\\
    \hspace*{-1cm}-(\partial_t p_2, \psi) + (\nabla p_2 , \nabla \psi) + (\alpha p_2, \psi) &= -(\beta \partial_{u_2} f_2 p_2, \psi) - (\partial_{u_1} f_2 p_1, \psi), \notag\\
    \hspace*{-1cm}p_1(T) = u_1(T) - O, \quad &p_2(0) = 0.
\end{align}

\noindent In light of \eqref{eq:preAdjEqn} and with the Riesz isomorphism associated to the corresponding duality pairing, the optimality condition reads as:
\begin{align}
\label{eq:optCond}
(G_{\bthet}^* \p + J_{\bthet},  \bthet - \bar \bthet )_{\TH} \ge 0, \quad \forall \bthet \in \THad.
\end{align}
\end{proof}

\noindent We now establish the stability result for the adjoint equation
\begin{lem}
\label{lem:stabAdj}
Let $\p^1, \p^2$ be two solutions to equation \eqref{eq:adjoint} generated by two $\u^1,\bthet^1$ and $\u^2, \bthet^2$ respectively. Then inequality \eqref{eq:stabAdjIneq} holds.
\end{lem}
\begin{proof} 

Let $\q := \p^1 - \p^2$, $\v := \u^1 - \u^2$ and $\bvarth := \bthet^1 - \bthet^2$
\begin{align}
\label{eq:stabAdj}
    \begin{pmatrix} -\partial_t - \nabla \cdot (\sigma^1 \nabla ) + (\nabla \kappa^1   + \delta^1 \nabla u^1_2 ) \cdot \nabla + \mu^1 g_{11}  &  -\nabla \cdot (\delta^1 u^1_1 \nabla ) + \mu^1 g_{12}\\
      \beta^1g_{21}  & -\partial_t -\lap + \alpha^1  + \beta^1g_{22}
    \end{pmatrix} &\begin{pmatrix} q_1 \\ q_2 \end{pmatrix} = 
    \begin{pmatrix} R_1   \\  R_2\end{pmatrix} \notag\\
    \nabla_n q_1 = 0, \quad   \nabla_n q_2 = 0, \quad q_1(T) = v_1(T), \quad q_2(T) &= 0. 
\end{align}
\noindent where $R_1 := \nabla \cdot (\vartheta_1 \nabla p^2_1) - (\nabla \vartheta_2   + \delta^1 \nabla v_2  + \vartheta_3 \nabla u^2_2) \cdot \nabla p^2_1 + \nabla \cdot (\delta^1 v^1_1 \nabla p^2_2 + \vartheta_3 u^2_1 \nabla p^2_2) + g^3$
\noindent
\[
\begin{array}{l l l}
R_2 := \vartheta_4 p^2_2  + g^4 , & g^3(\v,\bvarth) := g^{3,1}(\v) + g^{3,2}(\bvarth), & g^4(\v,\bvarth) := g^{4,1}(\v) + g^{4,2}(\bvarth), \\
g_{11} := \partial_{u_1} f_1(\u^1), &  g_{12} := \partial_{u_2} f_1(\u^1), &  g_{21} :=  \partial_{u_1} f_2(\u^1), \\
g_{22} :=  \partial_{u_2} f_2(\u^1), & g^3_{11} := \mu^1 \partial^2_{u_1 u_1} f_1(\u^1) p^2_1 v_1,& g^3_{21} := \mu^1 \partial^2_{u_2 u_1} f_1(\u^1) p^2_1 v_2,\\ g^3_{12} := \mu^1 \partial^2_{u_1 u_2} f_1(\u^1) p^2_2 v_1, & g^3_{22} := \mu^1 \partial^2_{u_2 u_2} f_1(\u^1) p^2_2 v_2,& g^4_{11} := \beta^1 \partial^2_{u_1 u_1} f_2(\u) p^2_1 v_1,\\
g^4_{21} := \beta^1\partial^2_{u_2 u_1} f_2(\u^1) p^2_1 v_2, & g^4_{12} := \beta^1\partial^2_{u_1 u_2} f_2(\u^1) p^2_1 v_1,& g^4_{22} := \beta^1 \partial^2_{u_2 u_1} f_2(\u^1) p^2_2 v_2, \\
g^{3,1} := g^3_{11} + g^3_{21} + g^3_{12} + g^3_{22} & g^{4,1} := g^4_{11} + g^4_{21} + g^4_{12} + g^4_{22} & g^{4,2} := \vartheta_6 g_{21} p^2_1 + \vartheta_5 g_{22} p^2_2 \\
g^{3,2} := \vartheta_5 g_{11} p^2_1 + \vartheta_5 g_{12} p^2_2 & &
 \end{array}\]
Since $\p^{1,2} \in \U \subset \U''$, $\u^{1,2} \in \U$ and $\bthet^{1,2} \in \TH$, the RHS terms $R_1$ and $R_2$ are elements of $\Z(T)$. Thus, letting $R := (R_1, R_2)$ the above equation \eqref{eq:stabAdj} can be abstractly written as:
\begin{align*}
    \Guad \q = R \text{ $\:$ in $\:$ } \U' 
\end{align*}
Due to the invertibility of $\Gu$ we get the existence of a unique solution to \eqref{eq:stabAdj}. Consequently, we have that
\begin{align}
\label{eq:stabAdjIneq}
\|\q\|_{\U} &\le \|\iGuad R\|_{\U} \le \|\iGuad\|_{L(\U',\U)} \|R \|_{\U'}  \notag\\
&\le \|\iGu \|_{L(\U',\U)} \|R\|_{\U'} \notag\\
&\le \const_{\Gu}(\u,\bthet) \|R\|_{\Z}.
\end{align}
\end{proof}



\noindent  Now we shall provide the smoothness (in terms of $\bthet$) result for the cost functional $J$.

\begin{thm}
\label{thm:smoothJac}
Let $\u \in \Y$, $\p \in \Y$ and $\bthet \in \Th$, then the cost functional $\bthet \mapsto \hat J(\bthet) := J(\u(\bthet),\bthet)$ is infinitely Frechet differentiable. Moreover, the mapping $\bthet \mapsto \nabla \hat J(\bthet)$, $\nabla \hat J: \Y \times \Y \times \THad \to \THad$ is Lipschitz continuous. 
\end{thm}
\begin{proof}
The infinite differentiability of $\hat J$ follows from the facts that $\bthet \mapsto \u(\bthet)$ is a smooth mapping (due to Lemma \ref{lem:invGx}) and $J$ is a quadratic functional of $\u$ and $\bthet$. Moreover, as already mentioned above, $\nabla \hat J$ takes the following form:
\begin{align*}
    \grad \hat J = \lambda \bthet + F(\bthet), \quad F(\bthet) := (G^*_{\bthet} \p)(\bthet) = -(G^*_{\bthet} \iGuad J_{\u})(\bthet).
\end{align*}
Moreover, based on \eqref{eq:diffTh}, $\grad \hat J$ can be explicitly written as
\begin{align}
\label{eq:Jacobian}
\nabla \hat J := \begin{pmatrix}
\nabla u_1 \cdot \nabla p_1 + \lam \: \th_1\\
\grad \cdot (u_1 \nabla p_1) + \lam \: \th_2\\
u_1 \nabla u_2 \cdot \nabla p_1 + \lam \: \th_3\\
\partial_{u_1} f_1(\u) p_1 + \lam \: \th_4\\
u_2 p_2 + \lam \: \th_5\\
\partial_v f_2(\u) p_2 + \lam \: \th_6\\
\end{pmatrix}.
\end{align}
Since $\u, \p \in \Y$, we get that $F(\bthet) \in \TH$. Due to the stability result of the adjoint $\p$ (see Lemma \ref{lem:stabAdj}) and the linear structure of $J_{\bthet}$,  for the Lipschitz continuity of $\grad J$ it suffices to only consider the operator $G_{\bthet}^* \p$ and establish its stability with respect to $\u, \p$ and $\bthet$. To this end we shall consider each component of the Jacobian vector function \eqref{eq:Jacobian} \\
\noindent
\textbf{Component 1:}
\begin{align*}
\|\nabla u_1 \nabla p_1 - \nabla u_2 \nabla p_2 \|_V &= \|\nabla u_1 (\nabla p_1 - \nabla p_2) + \nabla p_2 (\nabla u_1 - \nabla u_2) \|_V \\
&\le \|\nabla u_1\|_{L^4} \|(\nabla p_1 - \nabla p_2)\|_{L^4} + \|\nabla p_2\|_{L^4} \|(\nabla u_1 - \nabla u_2)\|_{L^4}  \\
&\le \|u_1\|_{Z} \|p_1 - p_2\|_{Z} + \|p_2\|_{Z} \|u_1 - u_2\|_{Z} \\
\Ra \|\nabla u_1 \nabla p_1 - \nabla u_2 \nabla p_2 \|_{V(T)} &\le \|u_1\|_{Y(T)} \|p_1 - p_2\|_{Z(T)} + \|p_2\|_{Y(T)} \|u_1 - u_2\|_{Z(T)}.
\end{align*}
\textbf{Component 2:}
\begin{align*}
\|\grad \cdot (u_1 \nabla p_1 - u_2 \nabla p_2 ) \|_V &= \|\grad \cdot \Big(u_1 (\nabla p_1 - \nabla p_2) + p_2 (\nabla u_1 - \nabla u_2)\Big) \|_V \\
&\hspace*{-3cm}\le \|\grad \cdot \Big(u_1 (\nabla p_1 - \nabla p_2) + p_2 (\nabla u_1 - \nabla u_2)\Big)\|_V  \\
&\hspace*{-3cm}\le \|\grad \cdot \Big(u_1 (\nabla p_1 - \nabla p_2) \Big)\| + \|\grad \cdot \Big(p_2 (\nabla u_1 - \nabla u_2)\Big)\|_V  \\
&\hspace*{-3cm}\le \|\grad u_1 \cdot (\nabla p_1 - \nabla p_2) + u_1 (\lap p_1 - \lap p_2)\| + \|\grad p_2 \cdot (\nabla u_1 - \nabla u_2) + p_2 (\lap u_1 - \lap u_2)\|_V  \\
&\hspace*{-3cm}\le \|\grad u_1\|_{L^4} \|\nabla p_1 - \nabla p_2\|_{L^4} + \|u_1\|_{\Li} \|\lap p_1 - \lap p_2\|_V + \|\grad p_2\|_{L^4} \|\nabla u_1 - \nabla u_2\|_{L^4} \\
&\hspace*{2cm} + \|p_2\|_{\Li} \|\lap u_1 - \lap u_2\|_V \\
&\hspace*{-3cm}\le \|u_1\|_Z \|p_1 - p_2\|_Z + \|u_1\|_{\Li} \|p_1 - p_2\|_Z + \|p_2\|_Z \|u_1 - u_2\|_Z + \|p_2\|_{\Li} \|u_1 - u_2\|_Z  \\
&\hspace*{-3cm}\le (\|u_1\|_Z + \|u_1\|_{\Li}) \|p_1 - p_2\|_Z  + (\|p_2\|_Z + \|p_2\|_{\Li}) \|u_1 - u_2\|_Z   \\
&\hspace*{-3cm}\le \|u_1\|_Z \|p_1 - p_2\|_Z + \|u_1\|_{\Li} \|p_1 - p_2\|_Z + \|p_2\|_Z \|u_1 - u_2\|_Z + \|p_2\|_{Z} \|u_1 - u_2\|_Z  \\
&\hspace*{-3cm}\le 2\|u_1\|_Z \|p_1 - p_2\|_Z  + 2\|p_2\|_Z \|u_1 - u_2\|_Z  \\
\Ra \|\grad \cdot (u_1 \nabla p_1 - u_2 \nabla p_2 ) \|_{V(T)} &\le 2\|u_1\|_{Y(T)} \|p_1 - p_2\|_{Z(T)}  + 2\|p_2\|_{Y(T)} \|u_1 - u_2\|_{Z(T)}.  
\end{align*}
\textbf{Component 3:}
\begin{align*}
\|u_1 \nabla u_1 \nabla p_1 - u_2 \nabla u_2 \nabla p_2 \|_V &= \|u_1 \nabla u_1 (\nabla p_1 - \nabla p_2) + u_1 \nabla p_2 (\nabla u_1 - \nabla u_2)\\
&\hspace*{2cm} + \nabla u_2 \nabla p_2 (u_1 - u_2)  \|_V \\
&\hspace*{-3cm}\le \|u_1 \nabla u_1\|_{L^4} \| (\nabla p_1 - \nabla p_2)\|_{L^4} + \|u_1 \nabla p_2\|_{L^4} \|(\nabla u_1 - \nabla u_2)\|_{L^4} \\
&\hspace*{2cm} + \|\nabla u_2 \nabla p_2\|_V \|u_1 - u_2\|_{\Li}   \\
&\hspace*{-3cm}\le \|u_1 \nabla u_1\|_{L^4} \|p_1 - p_2\|_{W} + \|u_1 \nabla p_2\|_{L^4} \|\nabla u_1 - \nabla u_2\|_{W} + \|\nabla u_2 \nabla p_2\|_V \|u_1 - u_2\|_{\Li}   \\
&\hspace*{-3cm}\le \|u_1\|_{\Li} \|\nabla u_1\|_{L^4} \|p_1 - p_2\|_{W} + \|u_1\|_{\Li} \|\nabla p_2\|_{L^4} \|\nabla u_1 - \nabla u_2\|_{W} \\
&\hspace{2cm}+ \|\nabla u_2\|_{L^4} \nabla p_2\|_{L^4} \|u_1 - u_2\|_{\Li}   \\
&\hspace*{-3cm}\le \|u_1\|_{\Li} \|u_1\|_Z \|p_1 - p_2\|_{Z} + \|u_1\|_{\Li} \|p_2\|_Z \| u_1 - u_2\|_{Z} + \|u_2\|_Z \|p_2\|_Z \|u_1 - u_2\|_{\Li}   \\
&\hspace*{-3cm}\le \|u_1\|_{Z} \|u_1\|_Z \|p_1 - p_2\|_{Z} + \|u_1\|_{Z} \|p_2\|_Z \| u_1 - u_2\|_{Z} + \|u_2\|_Z \|p_2\|_Z \|u_1 - u_2\|_{Z}  \\
\Ra \|u_1 \nabla u_1 \nabla p_1 - u_2 \nabla u_2 \nabla p_2 \|_{V(T)} &\le \|\u^2_1\|_{Y(T)} \|p_1 - p_2\|_{Z(T)} + \|u_1 p_2\|_{Y(T)} \| u_1 - u_2\|_{Z(T)} \\
&\hspace*{2cm} + \|u_2 p_2\|_{Y(T)} \|u_1 - u_2\|_{Z(T)}.
\end{align*}
\textbf{Component 4:}
\begin{align*}
\|\partial_{u_1} f^1_1 p_1 - \partial_{u_1} f^2_1 p_2 \|_V &\le \|\partial_{u_1} f^1_1 p_1 - \partial_{u_1} f^1_1 p_2 + \partial_{u_1} f^1_1 p_2 - \partial_{u_1} f^2_1 p_2  \| \\
&\le \|\partial_{u_1} f^1_1 p_1 - \partial_{u_1} f^1_1 p_2\|_V + \|\partial_{u_1} f^1_1 p_2 - \partial_{u_1} f^2_1 p_2\|_V  \\
&\le \|\partial_{u_1} f^1_1\|_{W} \|p_1 - p_2\|_W + \|p_2\|_W \|\partial_{u_1} f^1_1 - \partial_{u_1} f^2_1 \|_{W} \\
&\le \|\partial_{u_1} f^1_1\|_{W} \|p_1 - p_2\|_W + \|p_2\|_W \|\partial^2_{\u} f^1_1\|_{\Li} \|\u^1 - \u^2\|_Z \\
\Ra \|\partial_{u_1} f^1_1 p_1 - \partial_{u_1} f^2_1 p_2 \|_{V(T)} &\le \ubd_{f_1} \|p_1 - p_2\|_{Z(T)} + \|p_2\|_{Y(T)} \ubd_{f_1} \|\u^1 - \u^2\|_{Z(T)}
\end{align*}
\textbf{Component 5:}
\begin{align*}
\|u_1  p_1 - u_2  p_2 \|_V &= \|u_1 ( p_1 -  p_2) + p_2 ( u_1 -  u_2) \|_V \\
&\le \| u_1\|_{L^4} \|( p_1 -  p_2)\|_{L^4} + \|p_2\|_{L^4} \|u_1 -  u_2\|_{L^4}  \\
&\le \|u_1\|_{W} \|p_1 - p_2\|_{W} + \|p_2\|_{W} \|u_1 - u_2\|_{W} \\
\Ra \|u_1  p_1 - u_2  p_2 \|_V &\le \|u_1\|_{Y(T)} \|p_1 - p_2\|_{Z(T)} + \|p_2\|_{Y(T)} \|u_1 - u_2\|_{Z(T)}. 
\end{align*}
\textbf{Component 6:}
\begin{align*}
\|\partial_{u_2} f^1_2 p_1 - \partial_{u_2} f^2_2 p_2 \| &\le \|\partial_{u_2} f^1_2 p_1 - \partial_{u_2} f^1_2 p_2 + \partial_{u_2} f^1_2 p_2 - \partial_{u_2} f^2_2 p_2 \| \\
&\le \|\partial_{u_2} f^1_2 p_1 - \partial_{u_2} f^1_2 p_2\|_V + \|\partial_{u_2} f^1_2 p_2 - \partial_{u_2} f^2_2 p_2\|_V  \\
&\le \|\partial_{u_2} f^1_2\|_{W} \|p_1 - p_2\|_W + \|p_2\|_W \|\partial_{u_2} f^1_2 - \partial_{u_2} f^2_2 \|_W  \\
&\le \|\partial_{u_2} f^1_2\|_{W} \|p_1 - p_2\|_W + \|p_2\|_W \|\partial^2_{\u} f^1_1\|_{\Li} \|\u^1 - \u^2\|_Z  \\
\Ra \|\partial_{u_2} f^1_2 p_1 - \partial_{u_2} f^2_2 p_2 \| &\le \ubd_{f_2} \|p_1 - p_2\|_W + \|p_2\|_{Y(T)} \ubd_{f_2} \|\u^1 - \u^2\|_{Z(T)}  
\end{align*}
\end{proof}

\noindent Finally, we {need} the following result for the numerical solution of the minimization problem.

\begin{thm}
\label{thm:projGradDesc}
Let $\bthet^* \in \THad$ be the solution of system \eqref{eq:optSys}. Then Algorithm \ref{algo:pgd}., i.e. the projected gradient descent method, generates a minimizing sequence $(\bthet_n)_{n \in \N} \in \THad$ that converges to $\bthet^*$ in $\TH$.
\end{thm}
\begin{proof}
\noindent Since $\TH$ is a Hilbert space and $\THad \subset \TH$ is a closed convex set, the optimality condition can be written as 
\begin{align}
\label{eq:projOptCond}
    \bthet^* = \projthad(\bthet^* - \gamma \nabla \hat J(\bthet^*)), 
\end{align}
where, $\gamma > 0$ is some arbitrary fixed constant, $\projthad(\bthet) = \argmin{\hat \bthet \in \THad} \|\hat \bthet - \bthet \|_{\TH}$ being the projection operator onto the convex subset $\THad$. First we notice that any arbitrary $\bthet$ obtained via the equation \eqref{eq:projOptCond} is an element of $\THad$. Without loss of generality, letting $\gamma := \lam^{-1}$  we have that
\begin{align*}
    \bthet^* &= \projthad(\bthet^* - \lam^{-1} \nabla \hat J(\bthet^*)) = \projthad( - \lam^{-1} F(\bthet^*) ) \\
\Ra \|\bthet^*\| &\le \| \projthad(0) + \bthet^*  - \projthad(0) \| \\
&\le \|\projthad(0)\| + \|\projthad(-\lam^{-1} F(\bthet^*)) - \projthad(0)\|\\
&\le \|\projthad(0)\| + \lam^{-1}\|F(\bthet^*)\|
\end{align*}
The last inequality follows from the nonexpansivenss property of the projection operator. Since $F(\cdot) \in \TH$ we get that $\bthet^* \in \TH$ and the projection operator $\projthad$ ensures $\bthet^* \in \THad$. Finally, due to Lipschitz continuity of $F$, 
we can invoke Theorem 2.4 from \cite{Hinze08} to conclude that the projected gradient descent generates a minimizing sequence $(\bthet)_{k} \subset \THad$. 
\end{proof}
\noindent This completes the analysis section. Next we perform numerical simulations and discuss the obtained results. 
\newcommand{\bxi}{\bm{\xi}}
\newcommand{\bpsi}{\bm{\psi}}
\renewcommand{\a}{\mathbf{a}}
\newcommand{\f}{\mathbf{f}}
\newcommand{\blp}{\big(}
\newcommand{\brp}{\big)}
\newcommand{\img}{\mathcal{I}}
\newcommand{\obs}{O}
\newcommand{\proc}{\mathfrak{P}}
\newcommand{\opt}{\mathfrak{E}}
\newcommand{\Gauss}{\mathfrak{G}}
\newcommand{\DS}{\mathfrak{D}}
\newcommand{\I}{\mathcal{I}}
\newcommand{\xstar}{\u^*}
\newcommand{\thstar}{\bthet^*}
\newcommand{\xhat}{\hat\u}
\newcommand{\thhat}{\hat\bthet}
\newcommand{\uhtau}{\u_{h,\tau}}
\newcommand{\phtau}{\p_{h,\tau}}
\newcommand{\thhtau}{\bthet_{h,\tau}}
\section{Numerical simulations}\label{sec:numerics}
In this section we present a numerical method to compute the optimal parameter vector $\bthet \in \TH$ that generates a specified target state which is provided as an input data. There are two main paradigms for numerically solving the OCP \eqref{eq:pseudOCP}: \emph{first optimize then discretize} and \emph{first discretize then optimize}. Both techniques produce the same outcome in case of a pure Galerkin approximation, however the former not only results in a strongly consistent scheme (in general), but also offers superior asymptotic convergence properties \cite{Collis2002, Becker2007}. Therefore, we adopt the former approach and use the optimality system \eqref{eq:optSys} as the starting point. Letting $W$ to be a suitable Hilbert space, the state equation \eqref{eq:optSysState} can be represented in the weak form in the following way: 
\begin{align*}
    (\partial_t u_1, \varphi) + (\sigma \nabla u_1 + u_1 \grad \kappa , \nabla \varphi) + (\delta u_1 \nabla u_2, \nabla \varphi) &= (\mu f_1, \varphi), \quad \forall \varphi \in W \\
    (\partial_t u_2, \psi) + (\nabla u_2 , \nabla \psi) + (\alpha u_2, \psi) &= (\beta f_2, \psi), \quad \forall \psi \in W \\
    u_1(0) &= \u_{1,0},\ u_2(0) = \u_{2,0}.
\end{align*}
Similarly, letting $O$  represent the final (target) value of $u_1$, i.e. $u_1(T) = O$, with $O$ being the specified data, the adjoint equation \eqref{eq:optSysAdjoint} reads
\begin{align*}
    -(\partial_t p_1, \varphi) + (\sigma \nabla p_1, \nabla \varphi) + ((\grad \kappa + \delta \nabla u_1) \cdot \nabla p_1, \varphi) &= -(\delta u_1 \nabla p_2, \varphi) - (\mu \partial_{u_1} f_1 p_1, \varphi) \\
    &\quad -(\nu \partial_{u_2} f_1 p_2, \varphi), \\
    -(\partial_t p_2, \psi) + (\nabla p_2 , \nabla \psi) + (\alpha p_2, \psi) = -(\beta \partial_{u_2} f_2 p_2, \psi) &- (\partial_{u_1} f_2 p_1, \psi), \\
    p_1(T) = u_1(T) - O, \quad p_2(0) &= 0
\end{align*}
Together, the above two equations can be compactly written as:
\begin{align}
\label{eq:weakNumForm}
\begin{split}
    \blp \partial_t \u, \bphi \brp + \blp \Aop(\u;\th) \u, \bphi \brp &= \blp \f(\u;\th),\bphi \brp, \quad \forall \bphi \in W \\
    - \blp \partial_t \p, \bpsi \brp + \blp \Aop^*(\p;\th,\u) \p ,\bpsi\brp &= \blp \g(\u;\th)\p,\bpsi \brp, \quad \forall \bpsi \in W \\
    \u(0) &= \u_0, \quad \p(T) = \p_T
\end{split}
\end{align}
System \eqref{eq:weakNumForm} is numerically solved by discretizing it both spatially and temporally. For the spatial discretization we use the finite element method. Consequently, we replace the space $W$ by a finite-dimensional subspace $W_h \subset W$ which consists of continuous piecewise polynomial functions of degree 1, spanned by a nodal basis $\{\varphi_j \}^{N_h}_{j=1}$, with $\text{dim}(W_h) = N_h.$ 
The time interval
$I := [0, T ]$, $T \in \R^{+}$ is divided into $N_{\tau}$ subintervals, each having width $\tau := {|I|\over N_{\tau}}$. Based on this, the temporal grid
points are denoted by $ I_{\tau} := (t_n)_{n \in \{0,\dots,N_{\tau}\}}$ with $t_n := n \tau$. Finally, let $\u^h_n := \u(t_n)^h, \bthet^h_n := \bthet(t_n)^h, \p^h_n := \p(t_n)^h$ denote the finite dimensional approximations of $\u, \p, \th$ at time point $t_n$, respectively. Then for all $\bphi^h, \bpsi^h \in W^h$, the discrete version of \eqref{eq:weakNumForm} is given as:
\begin{subequations}
\label{eq:discreteWeakForm}
\begin{align}
\label{eq:discreteWeakstate}
    \blp \u^h_{n+1}, \bphi^h \brp + \tau \blp \Aop(\u^h_n;\bthet^h_n) \u^h_n, \bphi^h \brp &= \blp \u^h_{n}, \bphi^h \brp + \blp \f(\u^h_n;\bthet^h_n), \tau \bphi^h \brp,  \\
\label{eq:discreteWeakadjoint}
    \blp \p^h_n, \bpsi^h \brp + \tau \blp \Aop^*(\p^h_n;\bthet^h_n,\u^h_n) \p^h_n ,\bpsi^h \brp &= \blp \p^h_{n+1}, \bpsi^h \brp + \blp \g(\u^h_n;\bthet^h_n)\p^h_n, \tau \bpsi^h \brp, \\
    \u^h_0 &= \u(0), \quad \p^h_{N_\tau} = \p(T) \notag
\end{align}
\end{subequations}
Given $\th_{h,\tau} := (\bthet^h_n)_{n \in \N, n < N_{\tau}}$, the finite-element scheme \eqref{eq:discreteWeakstate} can be used to obtain an approximate solution $\u_{h,\tau} := (\u^h_n)_{n \in \N, n < N_{\tau}}$ of the state equation. Analogously, given $\th_{h,\tau}$ and $\u_{h,\tau}$, the finite element scheme \eqref{eq:discreteWeakadjoint} 
can be used to generate an approximate solution $\p_{h,\tau} := (\p^h_n)_{n \in \N, n < N_{\tau}}$ of the adjoint equation. Subsequently, the approximates $\u_{h,\tau}$ and $\p_{h,\tau}$ can be used to compute a new $\th_{h,\tau}$ based on the optimality relation \eqref{eq:optCond}. 
This basically leads to the following iterative method, commonly known as the projected gradient descent method, for computing the optimal parameter function $\bthet^*$. The sequential steps of the procedure are described in Algorithm \ref{algo:pgd}. The algorithm can be viewed as a mapping $(O,\u_0) \mapsto \opt(O,\u_0) = (\u^*,\bthet^*)$, $\opt: \hott \times \hott \to \U \times \THad$, which takes an initial value $\u_0$ and a final (target) value $O$ and computes the optimal solution $\xhat \in \U$ and optimal parameter $\thhat \in \Th_{ad}$. Since Algorithm $\ref{algo:pgd}$, (i.e. the mapping $\opt$),  is a numerical method, it is clear that $\xhat$ and $\thhat$ are the discrete representatives of the corresponding true optimal functions $\u^*$ and $\bthet^*$. 

\SetKwFor{KwFor}{for}{do}{end}
\begin{algorithm}[!htb]
\caption{PGD \label{algo:pgd}}
 \KwData{$O \in \hott$, $\u_0 \in \hott$, $\eps > 0$, $\tau, h > 0$, $N_{\tau} \in \N$, $T > 0$}
 $\thhtau^0 := \bthet^0$ \\
 \KwFor{$k = 1, \dots$}
 {
   $\uhtau^{k+1}$ = $S^h_{\tau}(\thhtau^k)$  using \eqref{eq:discreteWeakstate} \\[1ex]
   $\phtau^{k+1}$ = $(S^h_{\tau})^*(\uhtau^{k+1},\thhtau^k)$  using \eqref{eq:discreteWeakadjoint} \\[1ex]
   $-\nabla J(\uhtau^{k+1},\phtau^{k+1},\thhtau^k) = -G^*_{\th} (\uhtau^{k+1}, \thhtau^{k}) \phtau^{k+1} - \lam \thhtau^{k}$ \\[1ex]
   $\thhtau^{k+1}(\gamma_k) = \projthad(\thhtau^{k} - \gamma_k \nabla J(\uhtau^{k+1},\phtau^{k+1},\thhtau^k))$, $\gamma^{k} \in \{1, {1\over 2}, {1\over 4}, \dots \}$ \\[1ex]
   if $\hat J(\thhtau^{k+1}) < \eps$ \\
    \hspace*{.2cm} exit
 }
\end{algorithm}

\begin{figure}
    \centering
    \subfloat[\scriptsize{Pattern-B \cite{PatB}}\label{fig:patternB}]{
  \includegraphics[width=3.75cm,height=3.2cm]{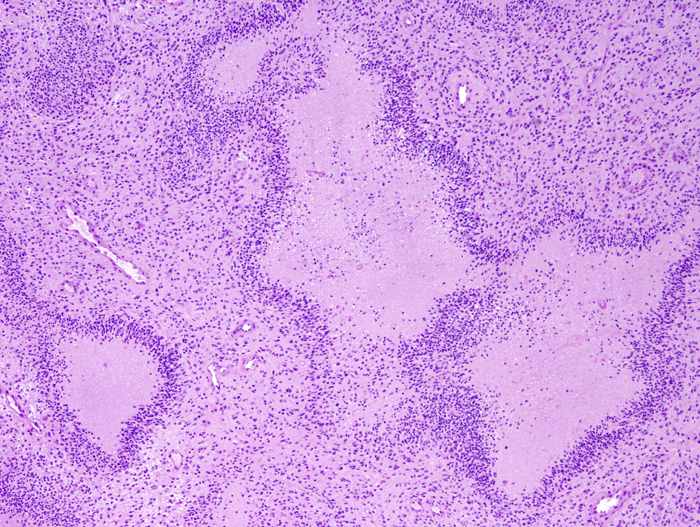}}  \hspace*{.1cm}
    \subfloat[\scriptsize{Pattern-C \cite{PatC}}\label{fig:patternC}]
    {\includegraphics[width=3.75cm,height=3.2cm]{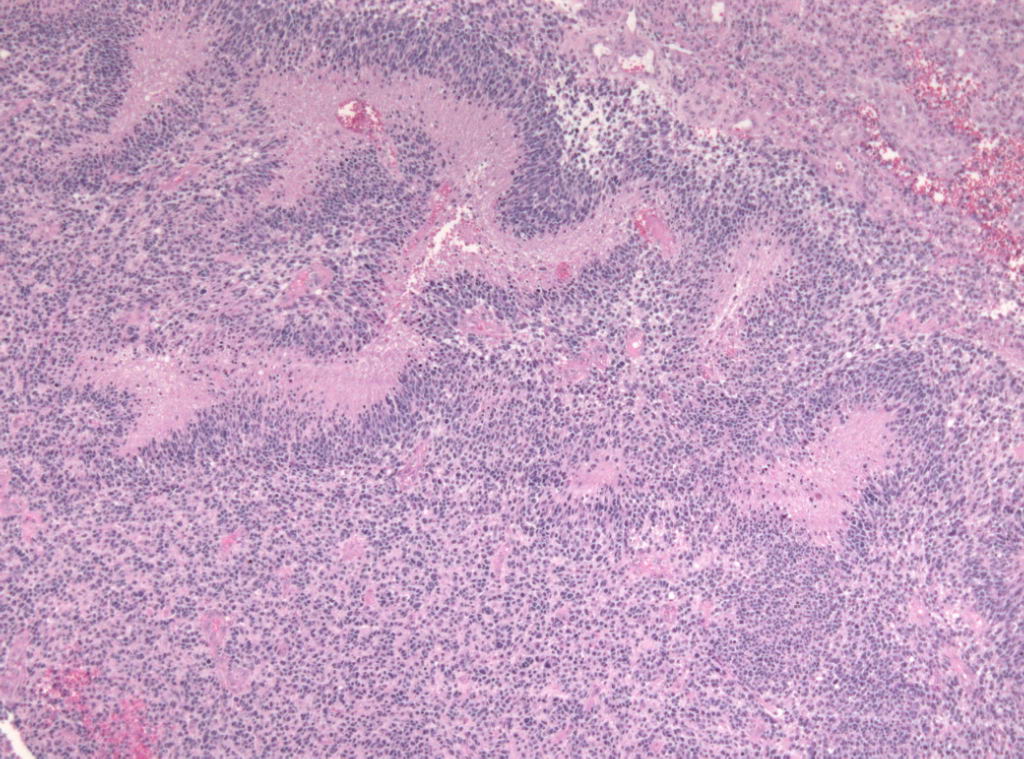}} \hspace*{.1cm}
    \subfloat[\scriptsize{Pattern-E \cite{PatE}}\label{fig:patternE}]
    {\includegraphics[width=3.75cm,height=3.2cm]{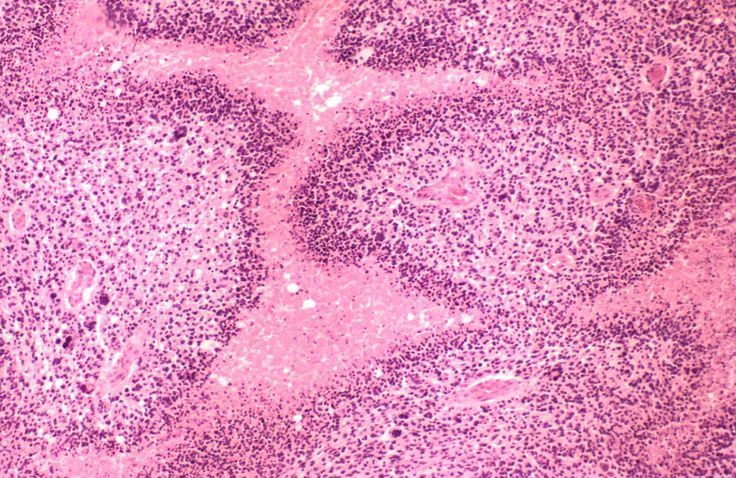}}
    \caption{Target pseudopalisade patterns \label{fig:targetPPDs}.}
\end{figure}

\begin{table}
\caption{Simulation parameters \label{tab:simpar}}
\begin{center}
\begin{tabular}{| c | c |}
\hline
\multicolumn{2}{|c|}{Numerical parameters}\\ 
\hline
T {\small(Total time)}& 10 \\\hline
$\tau$ {\small(Temporal step size)} & 0.1 \\\hline
$h_{x_1}$ {\small(Spatial step size along $x_1$)}& 0.1 \\\hline 
$h_{x_2}$ {\small(Spatial step size along $x_2$)}& 0.1 \\\hline 
$N_{x_1}$ {\small(Grid resolution along $x_1$)}& N {\small(image col size)} \\\hline
$N_{x_2}$ {\small(Grid resolution along $x_2$)}& M {\small(image row size)}\\\hline
$u_{1,0}$ {\small(initial value for $u_1$)} & 0.2 \\\hline
$u_{2,0}$ {\small(initial value for $u_2$)} & 0.5 \\\hline
$\th_{i,0}$ {\small(initial value for $\th_i$, for $i \in \{1,\dots,6\}$)} & 0 \\\hline
$\lam_{i}$ {\small(regularization parameter for $\th_i$ for $i \in \{1,\dots,6\}$ )}& $10^{-4}$\\\hline
\end{tabular}
\end{center}
\end{table}

\begin{table}
\caption{Model parameters \label{tab:modPar}}
\begin{center}
\begin{tabular}{| c | l | c |} 
\hline
\multicolumn{3}{|c|}{Scalar factors for migration coefficients}\\\hline
& \quad \quad \quad Phenomenological relevance & \\\hline
$\gamma_{_{\kappa}}$&{\small speed of pH-taxis for cancer cells}    & .01\\\hline
$\gamma_{_{\delta}}$& {\small speed of advection for cancer cells } & .001\\\hline
$\gamma_{_{pH}}$& {\small constant diffusion coefficient for protons} & .01\\\hline
\end{tabular}
\begin{tabular}{| c |  c |}
\hline
\multicolumn{2}{|c|}{Box constraints}\\\hline
$ {\sigma} \in [.0001, .01]$ & $ {\kappa} \in [-.01, .01]$ \\\hline
$ {\delta} \in [.0001, .01]$ & $ {\mu} \in [.0001, 10]$ \\\hline
$ {\alpha} \in [.0001, 10]$ & $ {\beta} \in [.0001, 10]$ \\\hline
\end{tabular}
\end{center}
\end{table}

\noindent A histological image $\I$ of taken at a specific time $T$ serves as the target state $O \approx \xhat(T)$  for the terminal optimal control problem \eqref{eq:pseudOCP}. It represents the observation data, based on which the above algorithm computes the optimal parameters for the model. The raw images $\I$ cannot be used directly in the optimization algorithm, but instead each need to be transformed into an image that represents the non-dimensionalized tumor density. This pre-processing step is performed using Algorithm \ref{algo:preproc}. Once again, the algorithm can be viewed as a mapping $\I \mapsto \proc(\I) =: \obs$, $\proc: \Loo \to \hott $,  that takes a raw data $\I \in \Loo$ as input and transforms it to an observation variable $\obs \in \hott$.

\begin{algorithm}[!htb]
\caption{Preprocessing steps \label{algo:preproc}}
 \KwData{$\I$: an RGB image of the tissue, with size  $M\times N$}
1. $\I_g = $ gray($\I$). Convert RGB to grayscale image \\
\hspace*{-.17cm}2. $\I_{gs} = \mathfrak{G} * \I_g$. Smoothen the image using a Gaussian filter $\mathfrak{G}$ \\
\hspace*{-.17cm}3. $\I_{gs} = \mathfrak{M} * \I_{gs}$. Remove 'salt and pepper' noise by applying a median filter \\
\hspace*{-.17cm}4. Generate an image mask $\mathfrak{m}$ by applying binary thresholding and performing morphological operation: \\
\qquad    4.1 apply binary to thresholding to extract dominant features \\
\qquad    4.2 perform morphological open operation to remove isolated features. This results in the required mask $\mathfrak{m}$\\
\hspace*{-.17cm}5. $\I_{gsm} = \I_{gs}(\mathfrak{m}) \land \I_{gs}(\mathfrak{m})$. Perform bitwise-and operation of the smoothened gray image with itself using the mask $\mathfrak{m}$. \\ 
\hspace*{-.17cm}6. $\I_{gsmi} = 1 - \I_{gsm}/255$. Normalize the image. \\
\hspace*{-.17cm}7. $O = \I_{gsmi}[::h_x, ::h_y]$. $h_x = {M \over m}$, $h_y = {N \over n}$. Downsample the $M \times N$ image to an $m \times n$ image.
\end{algorithm}
\noindent The final processed image data $O$ represents normalized volumetric concentration of the cancer cells. Thus it serves as a valid measurement for the non-dimensionalized model \eqref{eq:absPseudo}.

\begin{figure}
    \centering \hspace*{-1cm}
    \subfloat[\scriptsize{}\label{fig:pre_proc_target8}]
    {\includegraphics[trim=0cm 1cm 3.6cm 0cm, clip,scale=.105]{figs/results/targets/PatternB.png}} 
    \hspace*{.1cm}
    \subfloat[\scriptsize{}\label{fig:DS16_ker1_target8}]
    {\includegraphics[trim=0cm 1cm 3.6cm 0cm, clip,scale=.1]{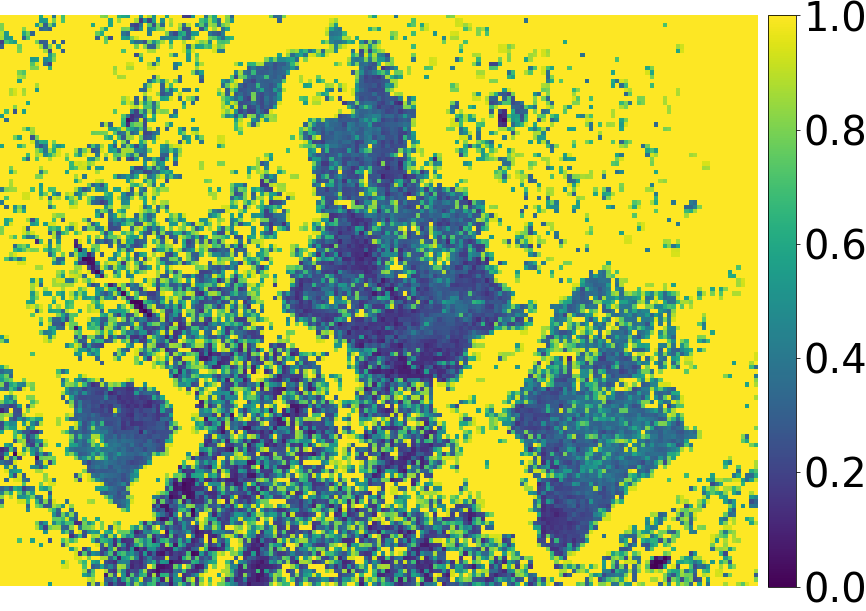}}
    \hspace*{.1cm}
    \subfloat[\scriptsize{}\label{fig:DS16_ker3_target8}]
    {\includegraphics[trim=0cm 0cm 3.6cm 0cm, clip,scale=.1]{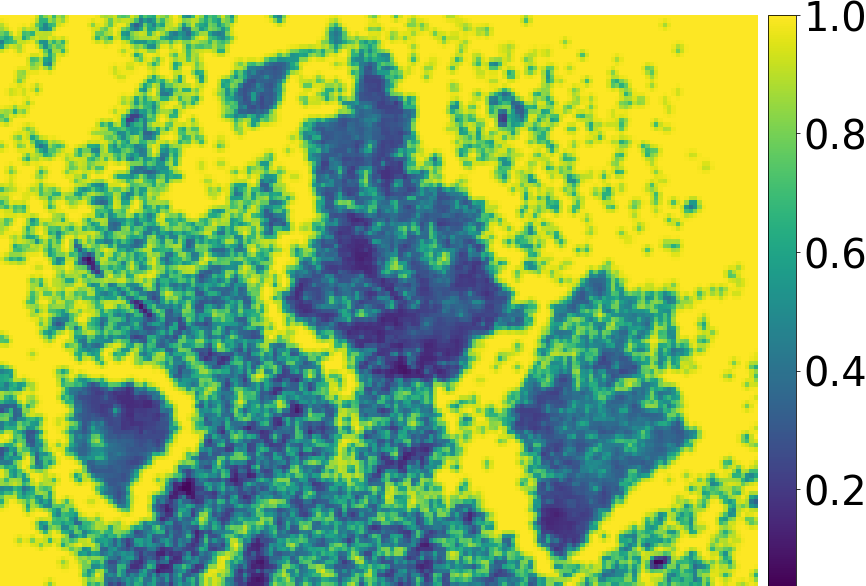}} 
    \hspace*{.1cm}
    \subfloat[\scriptsize{}\label{fig:DS16_ker5_target8}]
    {\includegraphics[trim=0cm 0cm 3.6cm 0cm, clip,scale=.1]{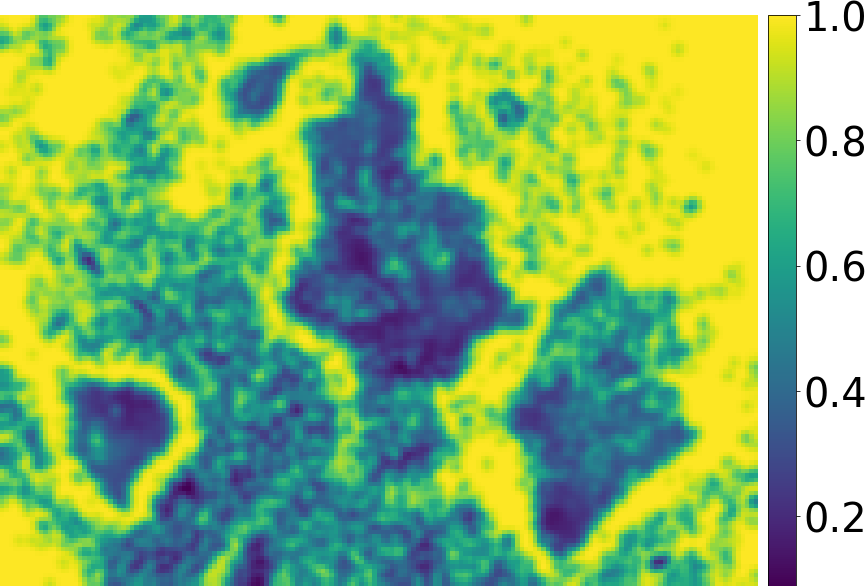}} 
    \hspace*{.1cm}
    \subfloat[\scriptsize{}\label{fig:DS16_ker7_target8}]
    {\includegraphics[trim=0cm 0cm 3.6cm 0cm, clip,scale=.1]{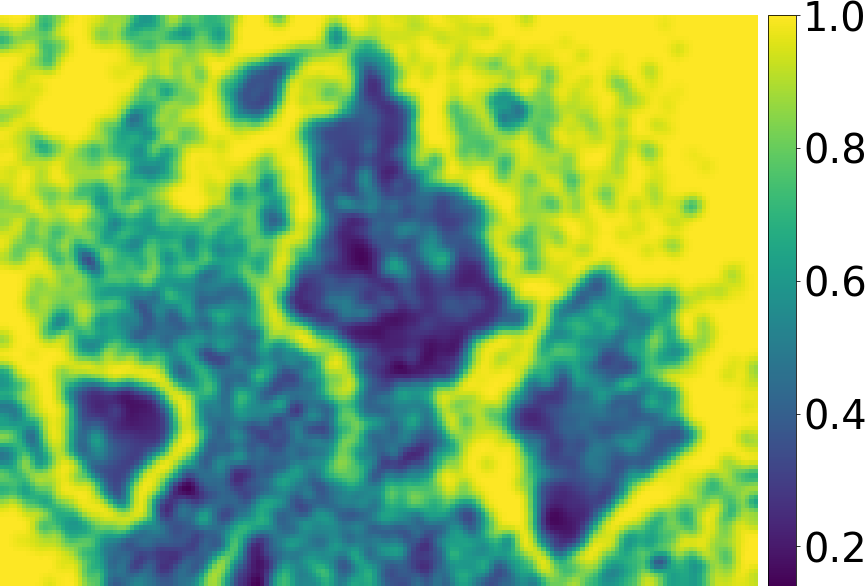}} 
    \caption{Pre-processed version of the raw image $\img$ (cf. Figure \ref{fig:patternB}) and its noisy perturbations (see text).\label{fig:noisy_per_target8}} 
\end{figure}
\begin{figure}
    \centering
    \subfloat[\scriptsize{The top, middle and bottom rows depict the reduction of $e_{\infty}$, $e_{rel}$ and $e_{\dom}$ errors respectively.}\label{fig:noisyOpt1}]
    {\includegraphics[trim=0cm 0cm 0cm 0cm, clip,scale=.14]{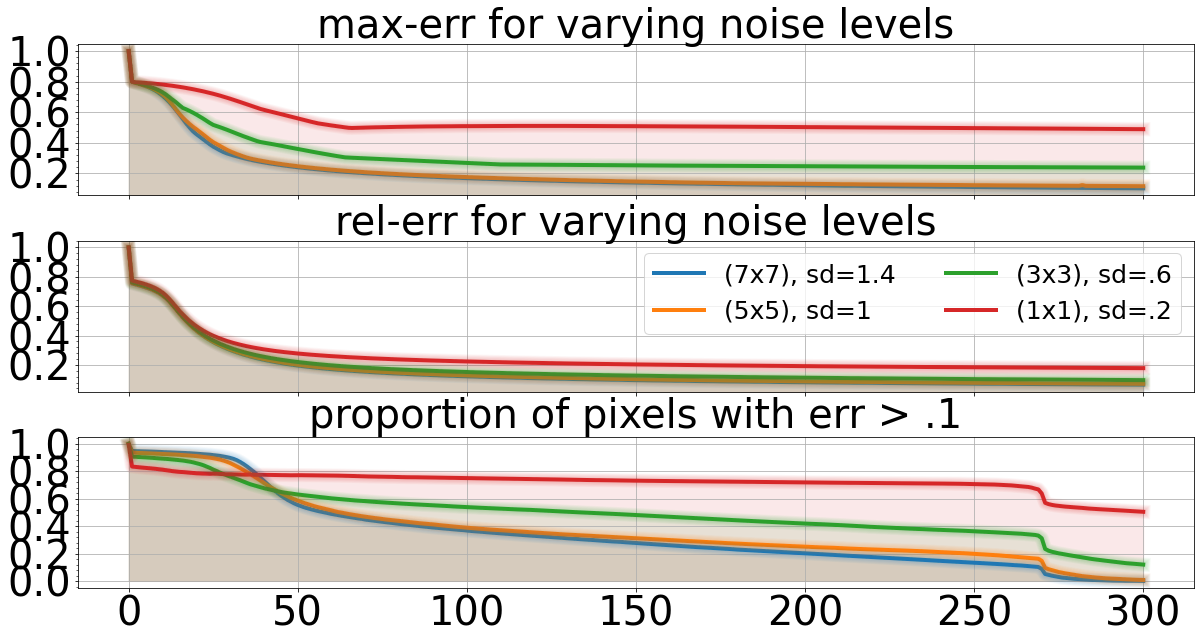}}
    \hspace*{.2cm}
    \subfloat[\scriptsize{Depicts reduction of absolute error $e_2$.}\label{fig:absErr}]
    {\includegraphics[trim=0cm 0cm 0cm 0cm, clip,scale=.14]{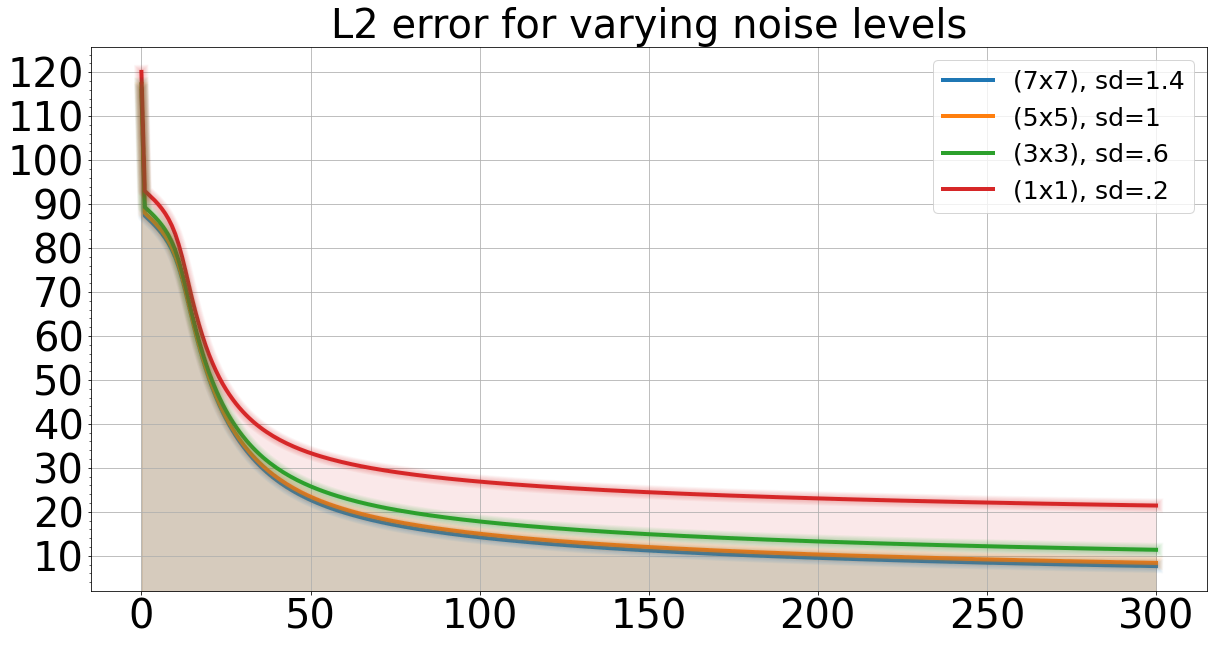}}
    \caption{Error reduction for noisy target images obtained using different standard deviation parameter ${\rm sd} = \sqrt{s} \in \{.2,.6, 1, 1.4\}$ in the smoothing kernel $\Gauss$ and downscalings $\dom_n$ with $n \in \{1,3,5,7\}.$ \label{fig:err_noisy_target8}}
\end{figure}

\paragraph{Evaluation of the optimization algorithm}
In this section we numerically investigate the minimizing properties of the Algorithm \ref{algo:pgd}. To this end we consider different noisy perturbations of a fixed target image and evaluate the obtained outputs of the Algorithm \ref{algo:pgd}. For the target image we consider the processed image $\obs = \proc(\I)$ ( Figure \ref{fig:pre_proc_target8}) obtained after applying the  Algorithm  \ref{algo:preproc} to the raw image $\img$ (Figure \ref{fig:patternB}).  
Let $\Gauss_{k,s}$ be the discrete Gaussian filter with kernel size $k$ and sigma (standard deviation) $s$ and $\DS_n$ to be the $n$-fold down sampling filter, then different perturbed version $\obs_k$ of $\obs$ is obtained applying $\Gauss_{k,s}$ and $\DS_n$ for different values of $k,s$ and $n$. Based on this, Figure \ref{fig:DS16_ker1_target8} is obtained as $\obs_1 = \Gauss{_{1,{1\over5}}} (\DS_4(\obs))$. Similarly, Figures  \ref{fig:DS16_ker3_target8} and \ref{fig:DS16_ker5_target8} are obtained as $\obs_3 = \Gauss{_{3,{3\over5}}} (\DS_4(\obs))$ and $\obs_5 = \Gauss{_{5,1}} (\DS_4(\obs))$ respectively. Due to the smoothing property of the Gaussian filter, increasing the kernel size and sigma results in smoother images, i.e. dampens spatial noise. As a result, we obtain that $\obs_5$ is smoother than $\obs_3$, which is in turn smoother than $\obs_1$.  
Now applying the minimization algorithm $\opt$ (Algorithm \ref{algo:pgd}) to these perturbed inputs we can gauge its performance. To this end, by letting $(\xhat_k,\thhat_k) = \opt(\obs_k)$,
we define the following error metrics:
\begin{align*}
    e^k_2 &:= \|\xhat_k(T) - \obs_k\|_{L^2}, \quad e^k_{\infty} := \|\xhat_k(T)-\obs_k\|_{L^{\infty}} \\
    e^k_{rel} &:= {\|\xhat_k(T) - \obs_k\|_{L^2} \over \| \obs_k\|_{L^2}}, \quad e^k_{\dom}(\eps) := {1\over \vert \dom \vert} \int_{\dom} \mathbbm{1}_{\{e^k_2 > \eps\}} \: dx. 
\end{align*}
Figure \ref{fig:err_noisy_target8} depicts the error reduction profiles corresponding to the noisy target images $\obs_1, \obs_3, \obs_5$ (Figures \ref{fig:DS16_ker1_target8}-\ref{fig:DS16_ker5_target8}). Based on this we can infer the following:
\begin{enumerate}
    \item as can be seen from Figure \ref{fig:absErr}, the absolute error $e^k_2$ tends to a stable low value for each $\obs_k$. It holds that $e^{k_1}_2 < e^{k_2}_2$ when $\obs_{k_1}$ is smoother than $\obs_{k_2}$.
    \item as can be seen from Figure \ref{fig:noisyOpt1}, for smoother target images the error reduction is relatively faster, especially for  $e_{\infty}$ and $e_{\dom}$. 
\end{enumerate}
The deterioration of error reduction for increased noise levels is expected and justified since (based on Theorem \ref{thm:existPseudo}) we expect $\xhat \in \U$, i.e. $\xhat(t) \in \hott$ for all $t \in [0,T]$. Thus for a noisy target pattern, the optimization can only be suboptimal due to the violation of spatial smoothness. 


%

\begin{figure}[!hbtp]
    \centering 
    \subfloat[\tiny{Target Pattern-A}\label{fig:target_pA}]
       {\includegraphics[width=4cm,height=3cm]{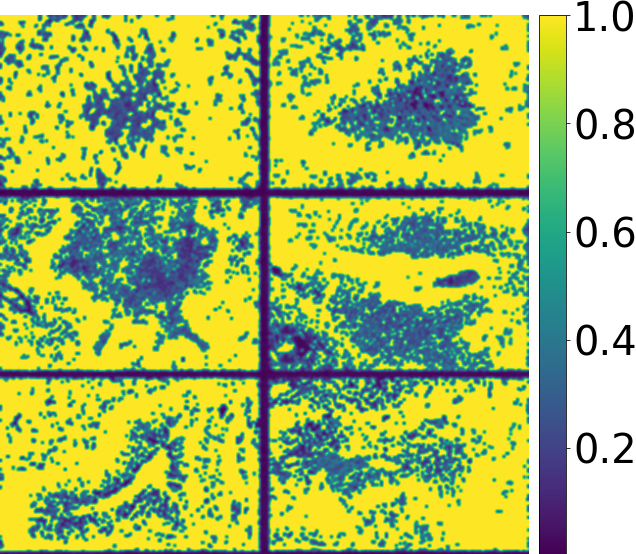}}
    \subfloat[\tiny{Final tumor density.}\label{fig:est_state_pA}]{
  \includegraphics[width=4cm,height=3cm]{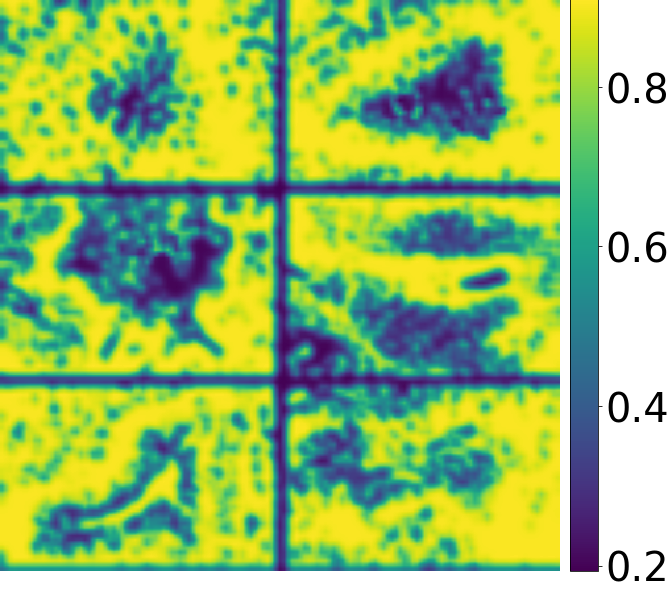}} 
    \subfloat[\tiny{Final acid concentration.}\label{fig:est_in_state_pA}]
    {\includegraphics[width=4cm,height=3cm]{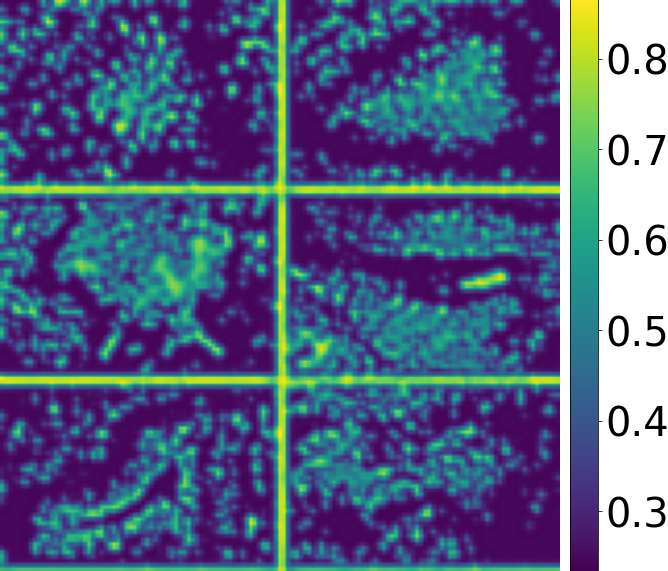}}
    \subfloat[\tiny{Error $\| C(T) - O\|^2$  }\label{fig:advec_final_pA}]{
  \includegraphics[width=4cm,height=3cm]{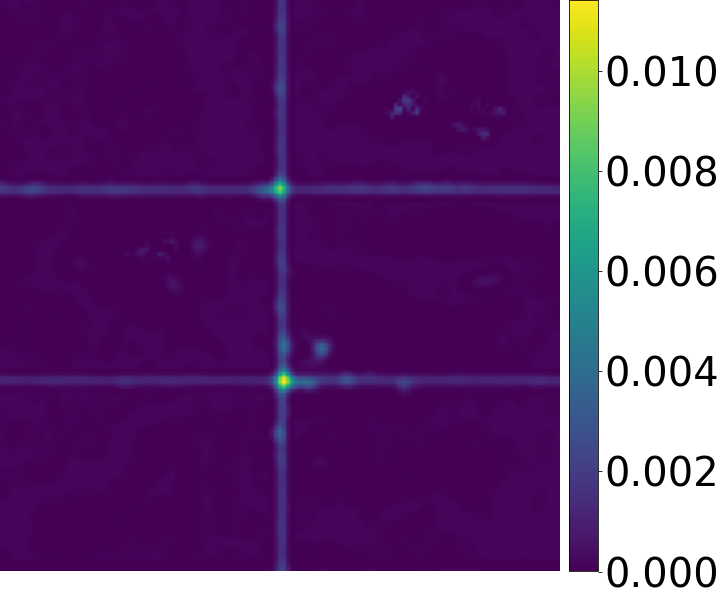}}
    \caption{Optimal cancer and acid distribution for the given target pseudopalisade pattern (see Fig 7. \cite{PatA}). \label{fig:patFormA}}
\end{figure}

\begin{figure}[!hbtp]
    \centering  
    \subfloat[\tiny{Target Pattern-B}\label{fig:target_pB}]
    {\includegraphics[width=4cm,height=3cm]{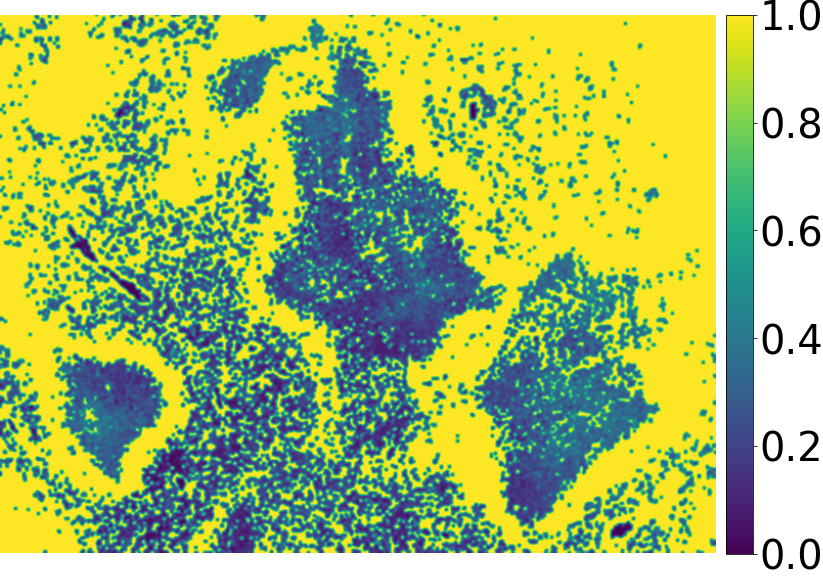}}
    \subfloat[\tiny{Final tumor density.}\label{fig:est_state_pB}]{
  \includegraphics[width=4cm,height=3cm]{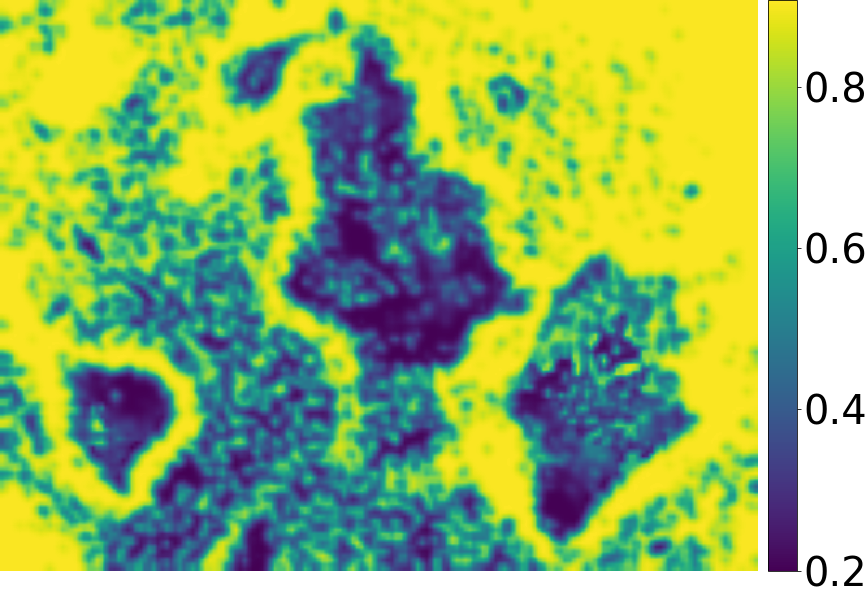}}
    \subfloat[\tiny{Final acid concentration.}\label{fig:est_in_state_pB}]
    {\includegraphics[width=4cm,height=3cm]{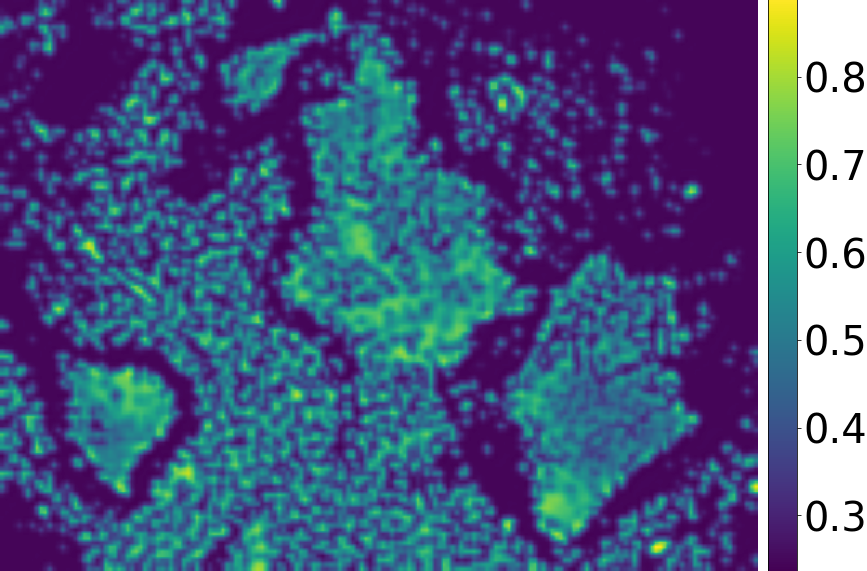}}
    \subfloat[\tiny{Error $\| M(T) - O\|^2$  }\label{fig:advec_final_pB}]{
    \includegraphics[width=4cm,height=3cm]{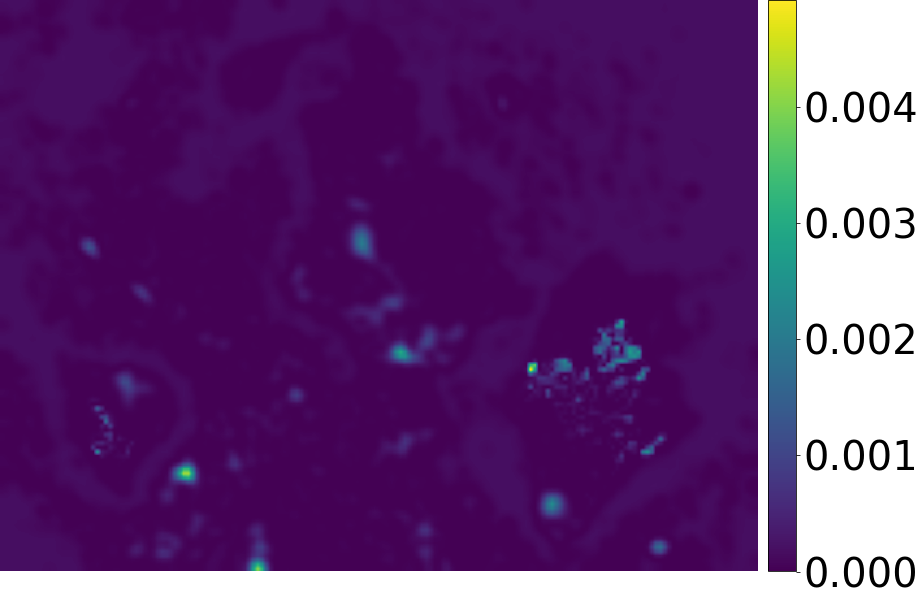}}
    \caption{Optimal cancer and acid distribution for the given target pseudopalisade pattern in Figure \ref{fig:patternB}. \label{fig:patFormB}}
\end{figure}

\begin{figure}[!hbtp]
    \centering 
    \subfloat[\tiny{Target Pattern-C}\label{fig:target_pC}]
    {\includegraphics[width=4cm,height=3cm]{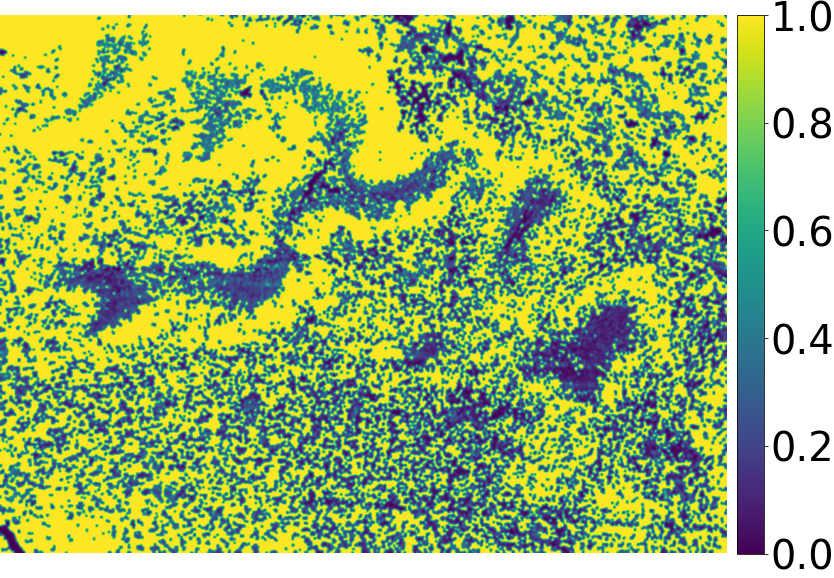}}
    \subfloat[\tiny{Final tumor density.}\label{fig:est_state_pC}]{
  \includegraphics[width=4cm,height=3cm]{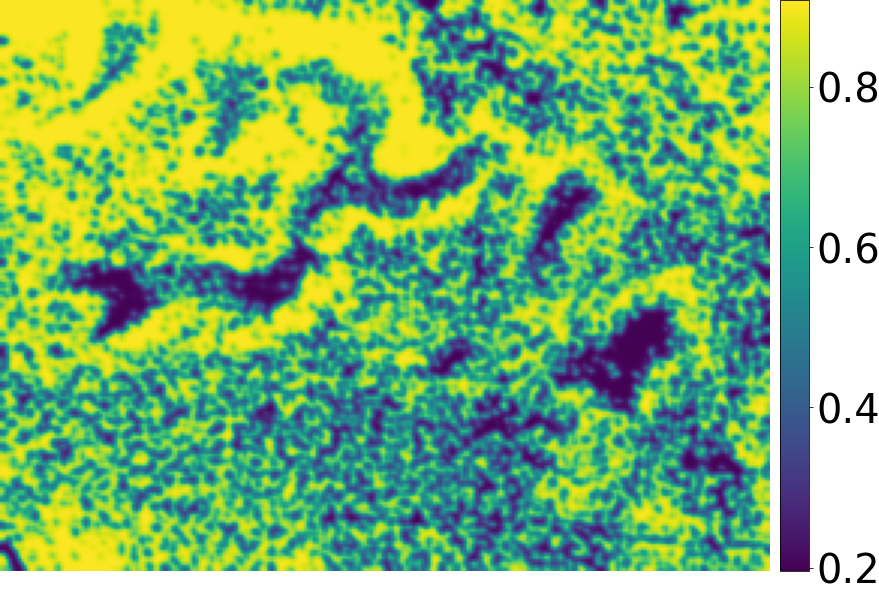}} 
    \subfloat[\tiny{Final acid concentration.}\label{fig:est_in_state_pC}]
    {\includegraphics[width=4cm,height=3cm]{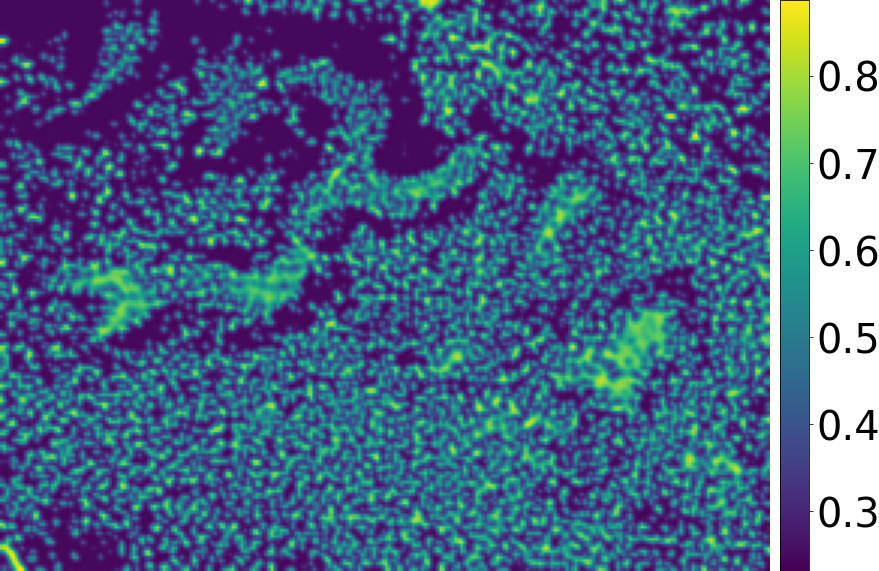}}
    \subfloat[\tiny{Error $\| M(T) - O\|^2$  }\label{fig:advec_final_pC}]{
  \includegraphics[width=4cm,height=3cm]{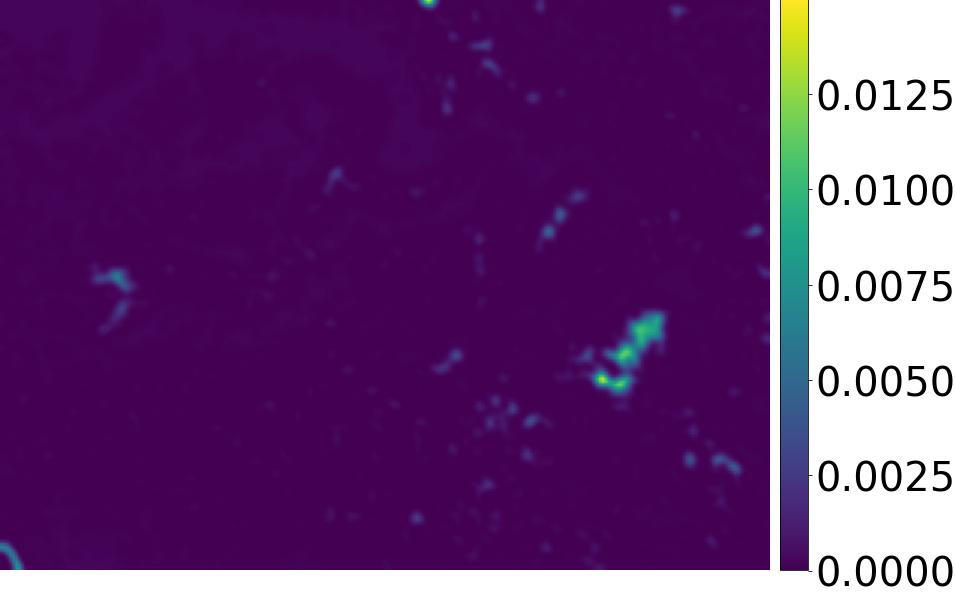}}
    \caption{Optimal cancer and acid distribution for the given target pseudopalisade pattern in Figure \ref{fig:patternC}. \label{fig:patFormC}}
\end{figure}

\begin{figure}[!hbtp]
    \centering 
    \subfloat[\tiny{Target Pattern-E}\label{fig:target_pE}]
    {\includegraphics[width=4cm,height=3cm]{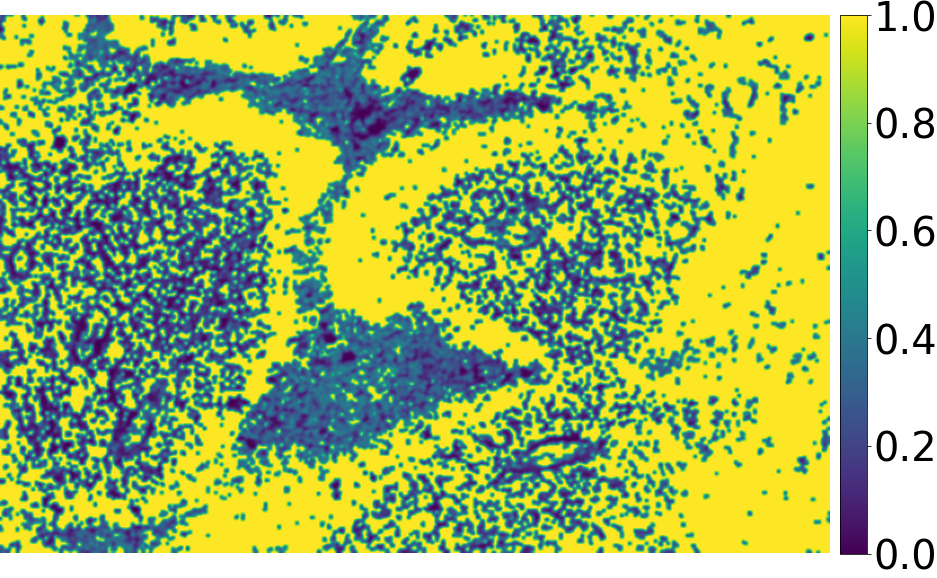}}
    \subfloat[\tiny{Final tumor density.}\label{fig:est_state_pE}]{
  \includegraphics[width=4cm,height=3cm]{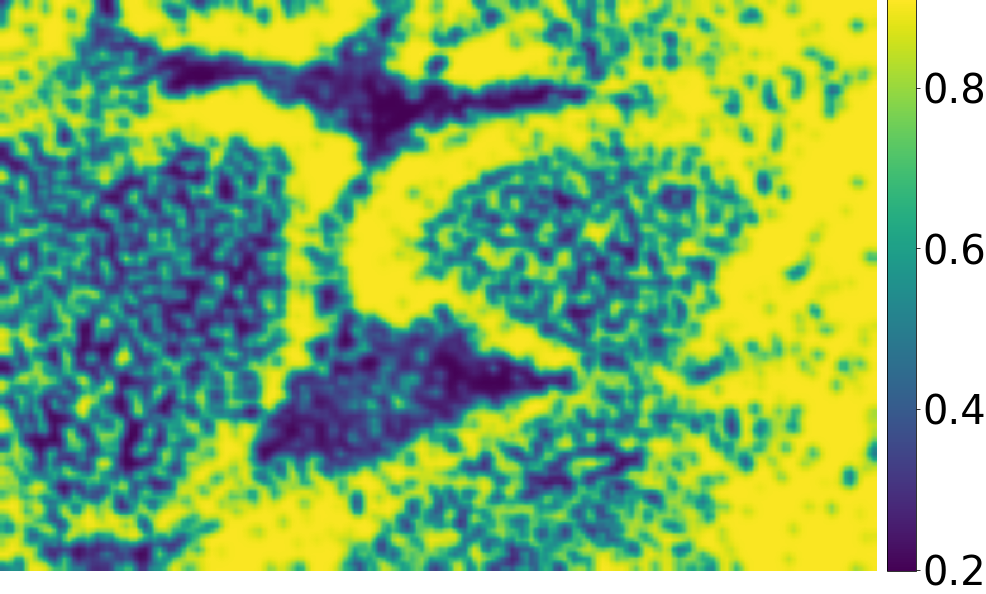}} 
    \subfloat[\tiny{Final acid concentration.}\label{fig:est_in_state_pE}]
    {\includegraphics[width=4cm,height=3cm]{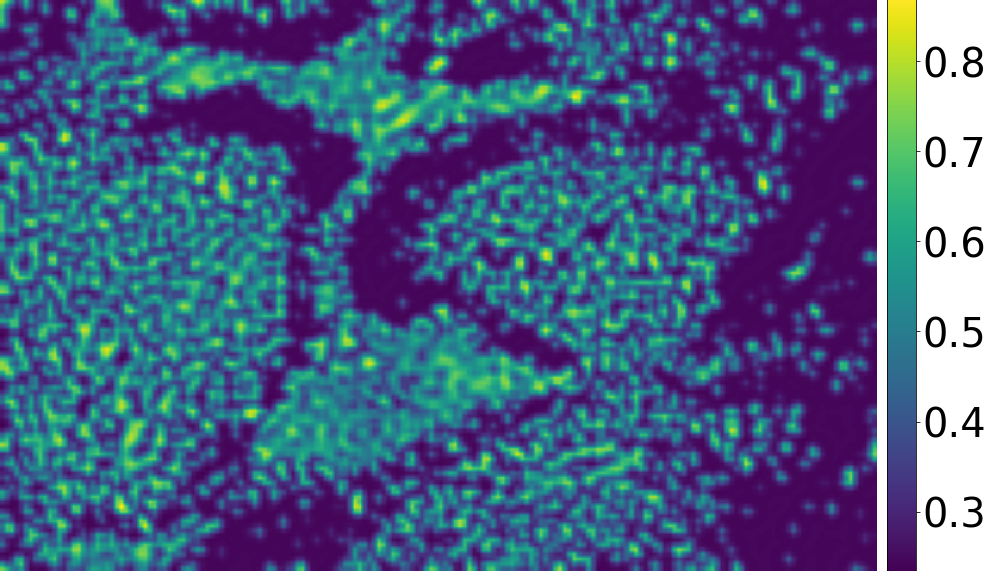}}
    \subfloat[\tiny{Error $\| M(T) - O\|^2$  }\label{fig:advec_final_pE}]{
	  \includegraphics[width=4cm,height=3cm]{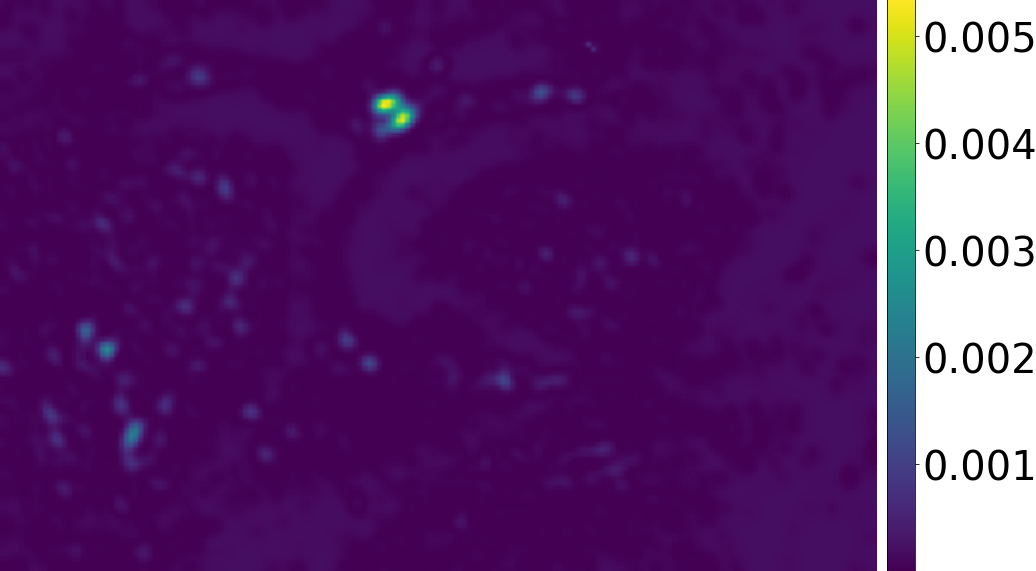}}
    \caption{Optimal cancer and acid distribution for the given target pseudopalisade pattern in Figure \ref{fig:patternE}. \label{fig:patFormE}}
\end{figure}

\begin{figure}[!hbtp]
    \centering  
    \subfloat[\tiny{t = 20 days}\label{fig:can_grow_pB_2}]
    {\includegraphics[width=4cm,height=3cm]{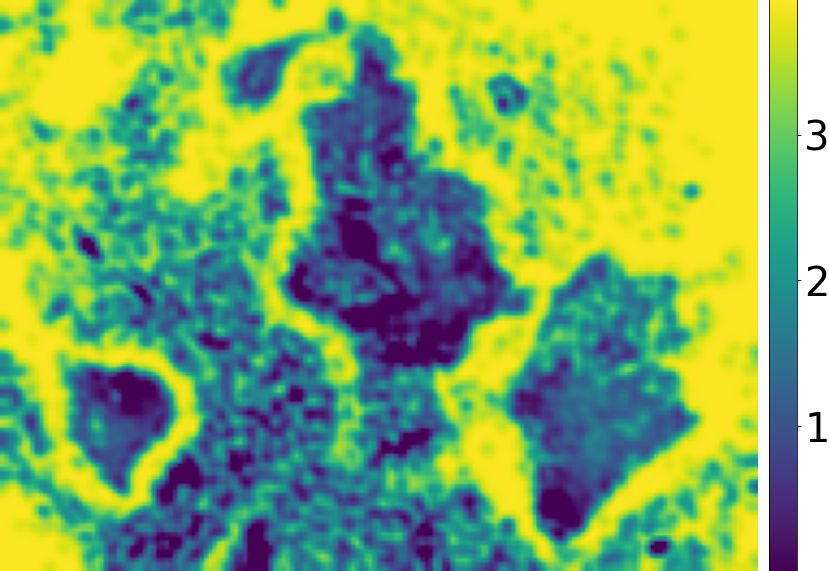}}\hspace*{.1cm}
    \subfloat[\tiny{t = 50 days}\label{fig:can_grow_pB_5}]
    {\includegraphics[width=4cm,height=3cm]{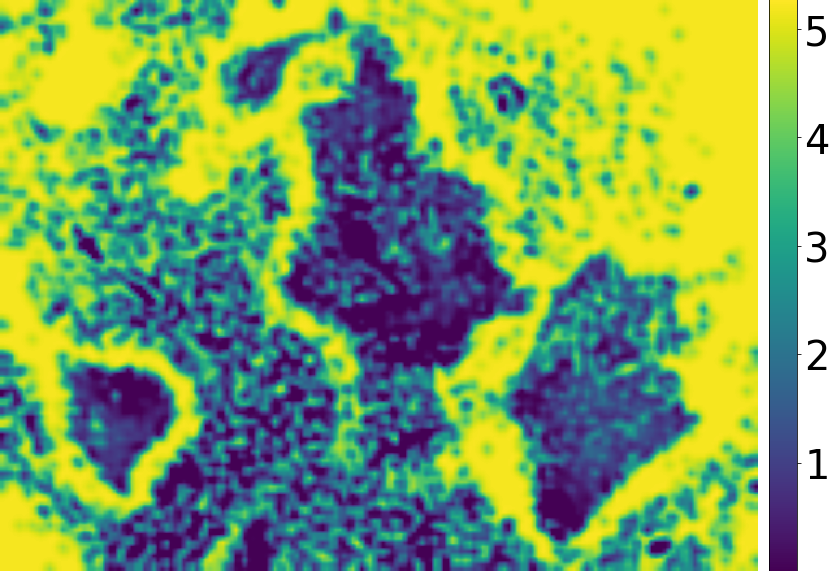}}\hspace*{.1cm}
    \subfloat[\tiny{t = 70 days}\label{fig:can_grow_pB_7}]
    {\includegraphics[width=4cm,height=3cm]{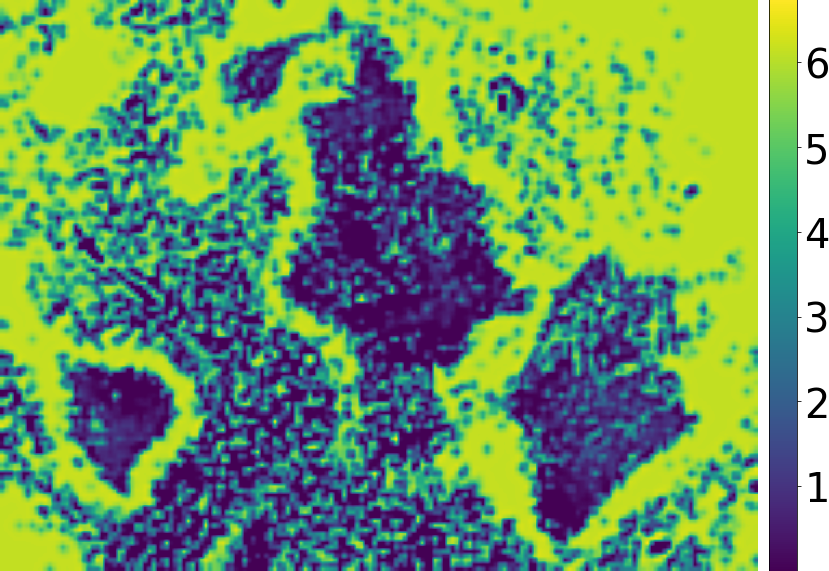}}\hspace*{.1cm}
    \subfloat[\tiny{t = 90 days}\label{fig:can_grow_pB_9}]
    {\includegraphics[width=4cm,height=3cm]{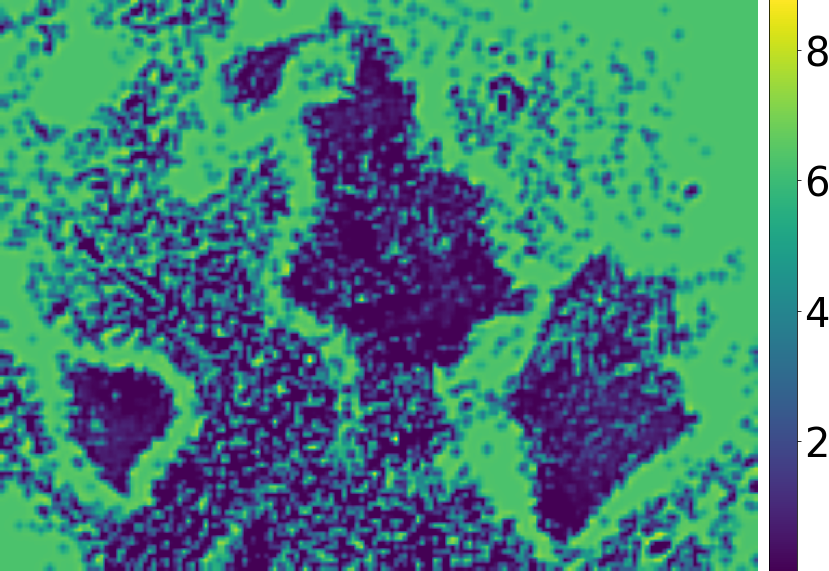}}\\
    \centering  
    \subfloat[\tiny{t = 20 days}\label{fig:prot_dec_pB_2}]
    {\includegraphics[width=4cm,height=3cm]{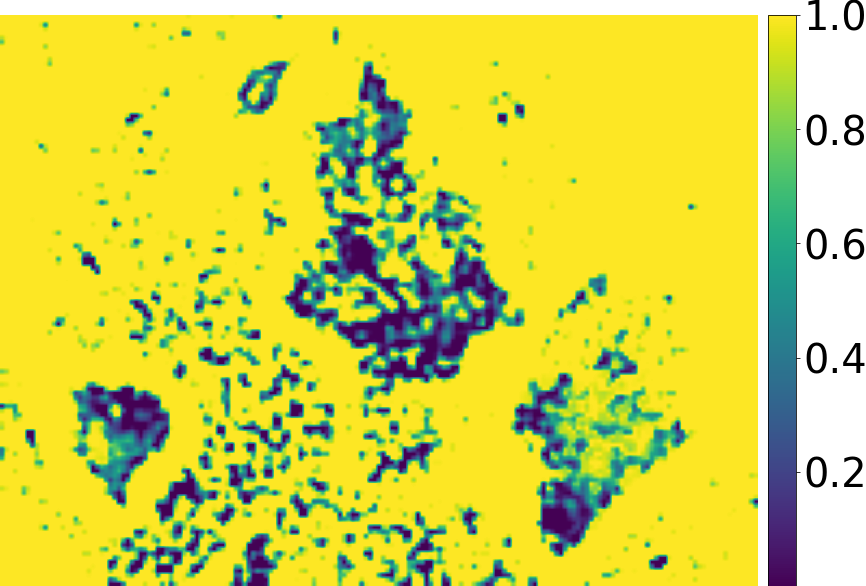}}\hspace*{.1cm}
    \subfloat[\tiny{t = 50 days}\label{fig:prot_dec_pB_5}]
    {\includegraphics[width=4cm,height=3cm]{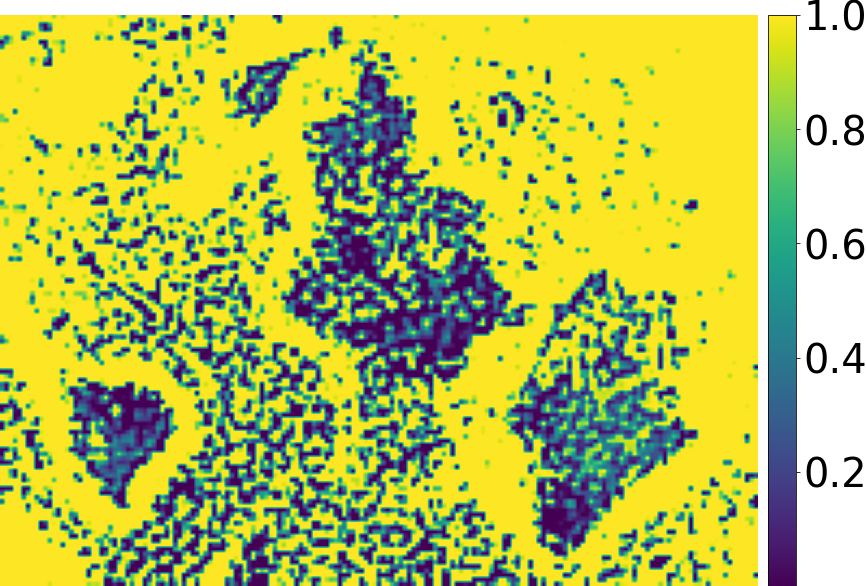}}\hspace*{.1cm}
    \subfloat[\tiny{t = 70 days}\label{fig:prot_dec_pB_7}]
    {\includegraphics[width=4cm,height=3cm]{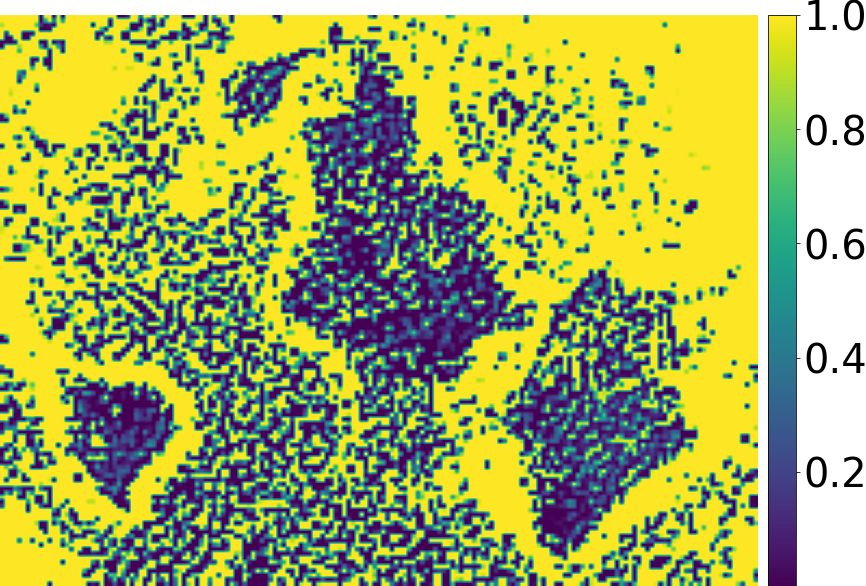}}\hspace*{.1cm}
    \subfloat[\tiny{t = 90 days}\label{fig:prot_dec_pB_9}]
    {\includegraphics[width=4cm,height=3cm]{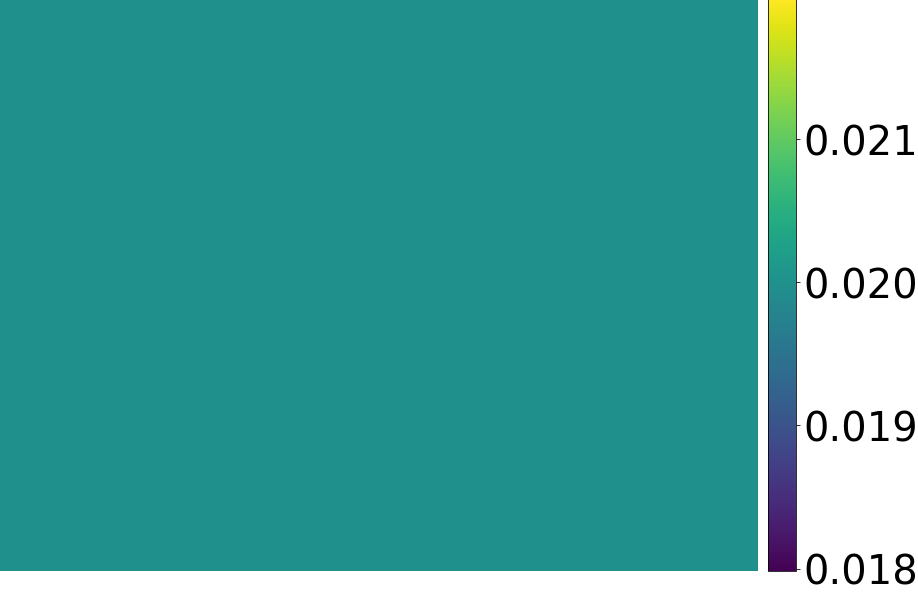}}\\
    \centering 
    \subfloat[\tiny{t = 20 days}\label{fig:prot_grow_pB_2}]
    {\includegraphics[width=4cm,height=3cm]{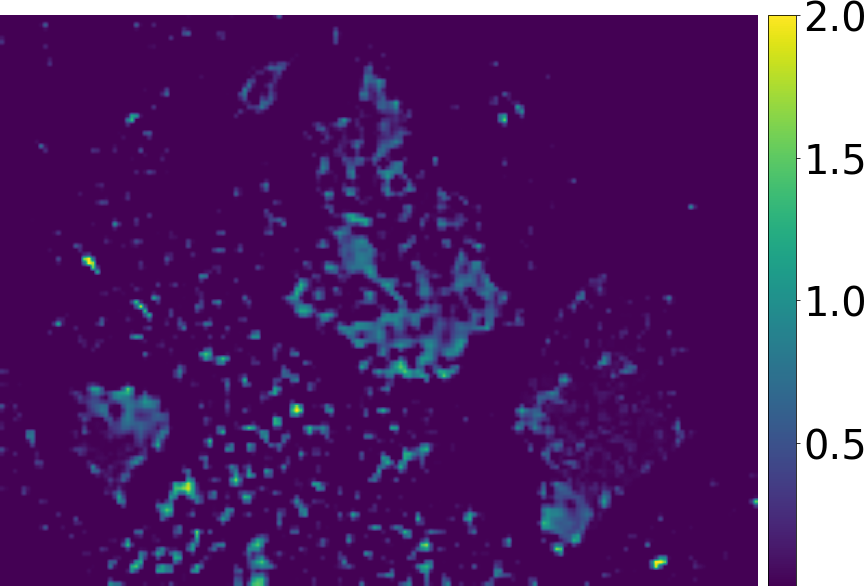}}\hspace*{.1cm}
    \subfloat[\tiny{t = 50 days}\label{fig:prot_grow_pB_5}]
    {\includegraphics[width=4cm,height=3cm]{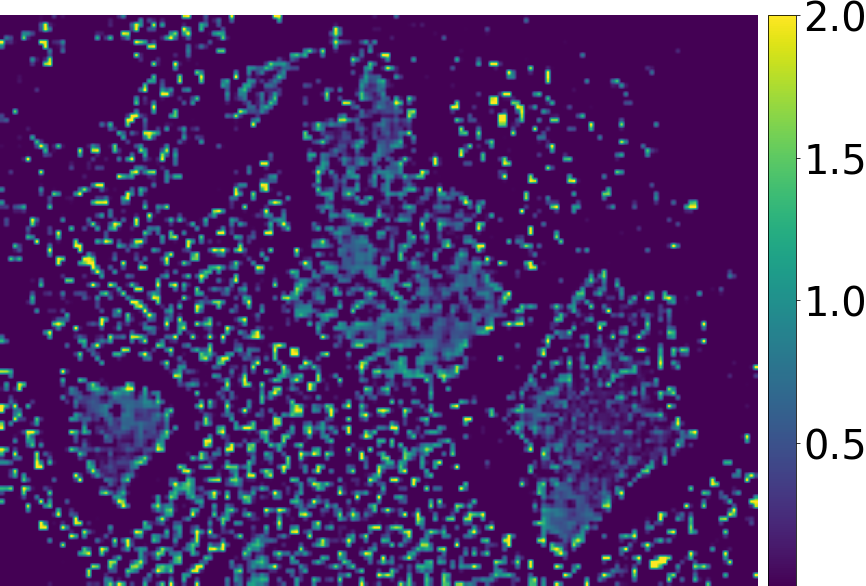}}\hspace*{.1cm}
    \subfloat[\tiny{t = 70 days}\label{fig:prot_grow_pB_7}]
    {\includegraphics[width=4cm,height=3cm]{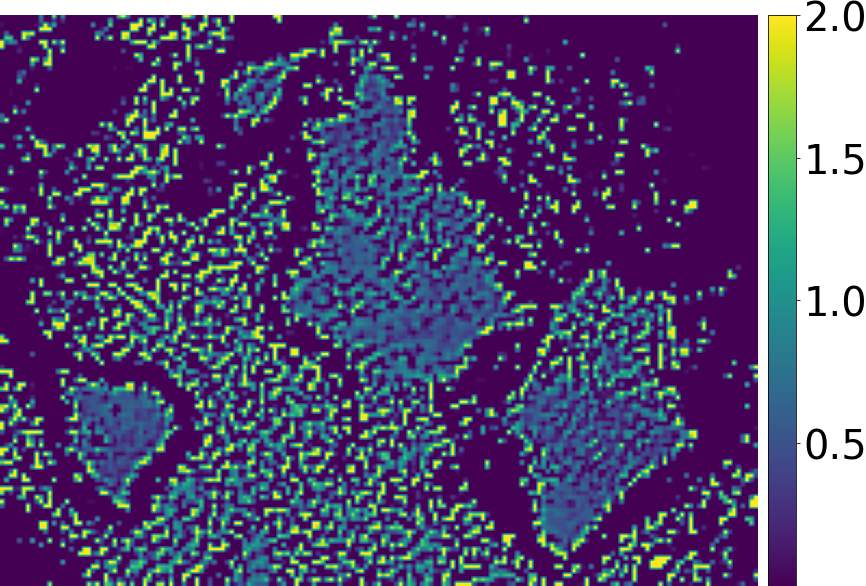}}\hspace*{.1cm}
    \subfloat[\tiny{t = 90 days}\label{fig:prot_grow_pB_9}]
    {\includegraphics[width=4cm,height=3cm]{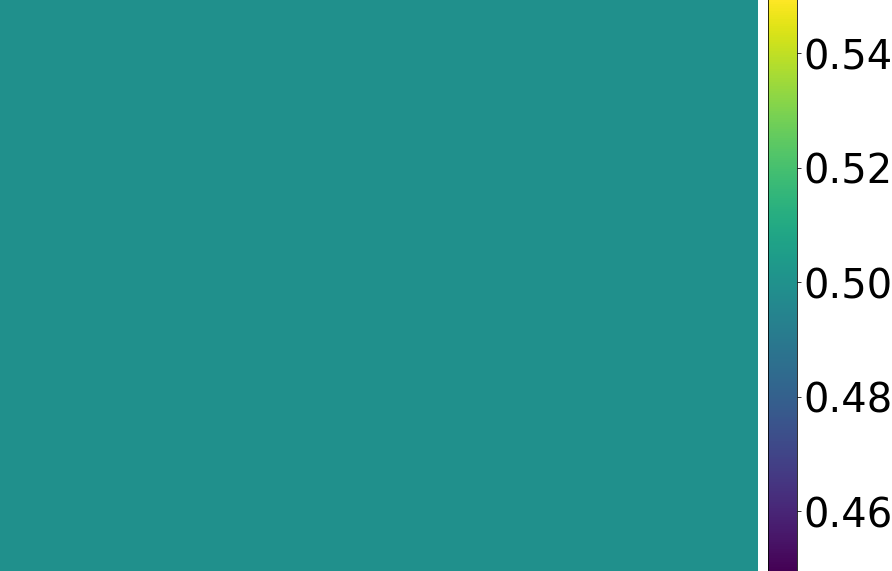}}\\
\caption{Evolution of growth and decay parameter functions for target Pattern-B from Figure \ref{fig:patternB}. The columns are sorted according to time (in days). The leftmost column is the tumor state at t=20, the middle left column is for t=50 days, the middle right one for t=70 days, and the rightmost column for t=90. The rows represent different parameter functions of the model \eqref{eq:absPseudo}. Arranged from top to bottom these are: the growth coefficient $\mu$ (1st row), acid removal rate $\alpha$ (2nd row), acid production rate $\beta$ (3rd row).  \label{fig:can_grow_pB}}
\end{figure}

\begin{figure}[!hbtp]
    \centering 
    \subfloat[\tiny{t = 20 days}\label{fig:can_diff_pB_2}]
    {\includegraphics[width=4cm,height=3cm]{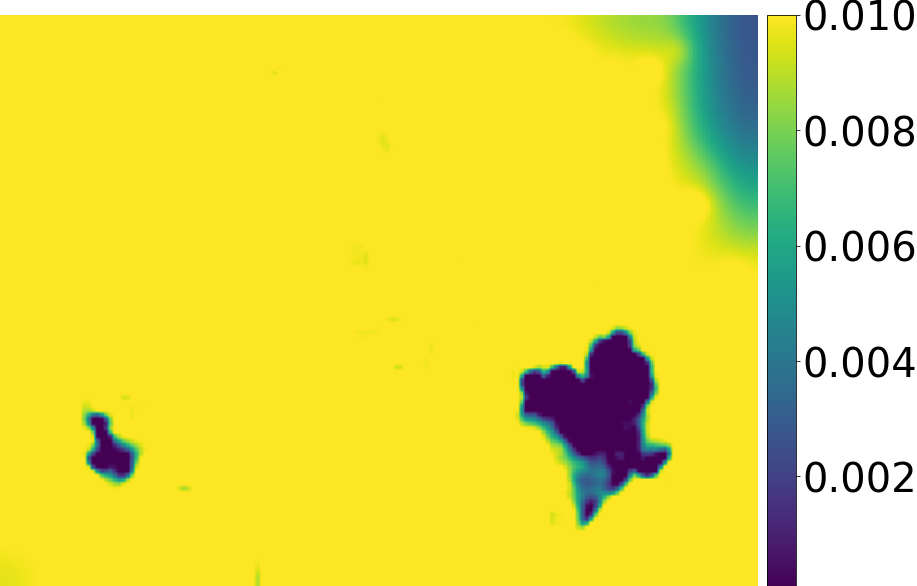}}\hspace*{.1cm}
    \subfloat[\tiny{t = 50 days}\label{fig:can_diff_pB_5}]
    {\includegraphics[width=4cm,height=3cm]{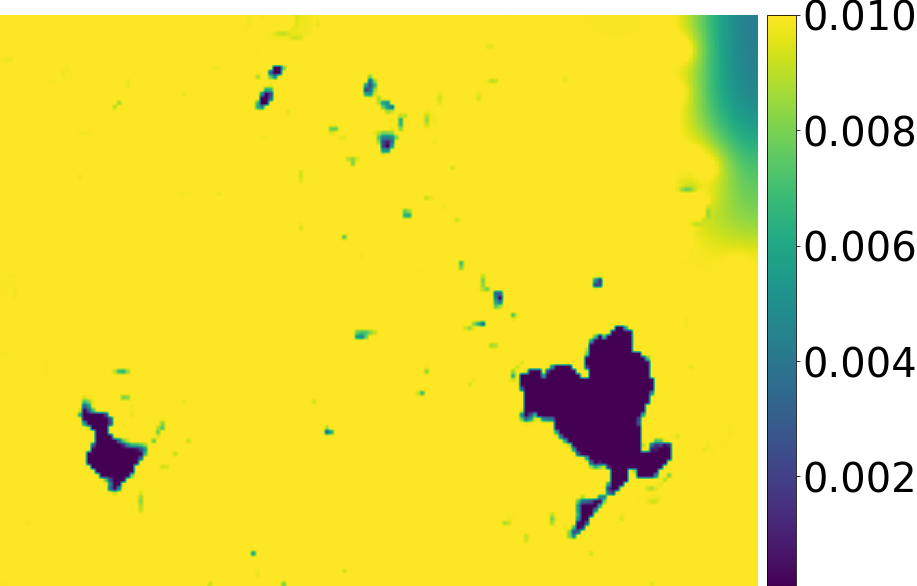}}\hspace*{.1cm}
    \subfloat[\tiny{t = 70 days}\label{fig:can_diff_pB_7}]
    {\includegraphics[width=4cm,height=3cm]{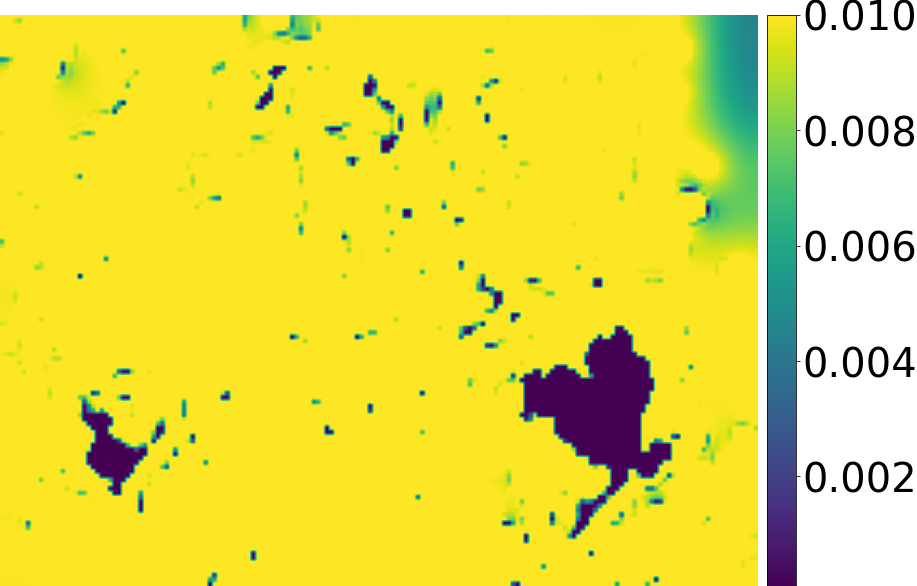}}\hspace*{.1cm}
    \subfloat[\tiny{t = 90 days}\label{fig:can_diff_pB_9}]
    {\includegraphics[width=4cm,height=3cm]{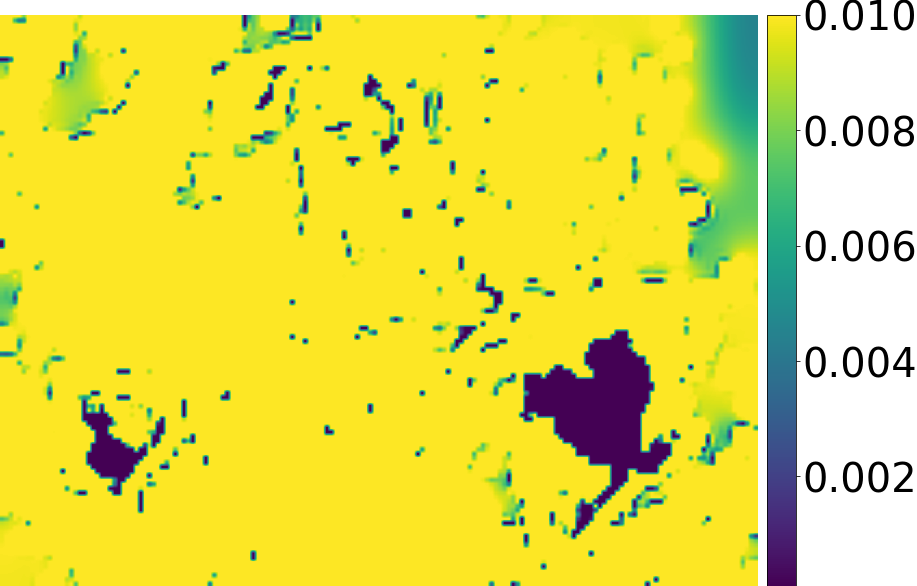}}\\

    \centering 
    \subfloat[\tiny{t = 20 days}\label{fig:can_advec_pB_2}]
    {\includegraphics[width=4cm,height=3cm]{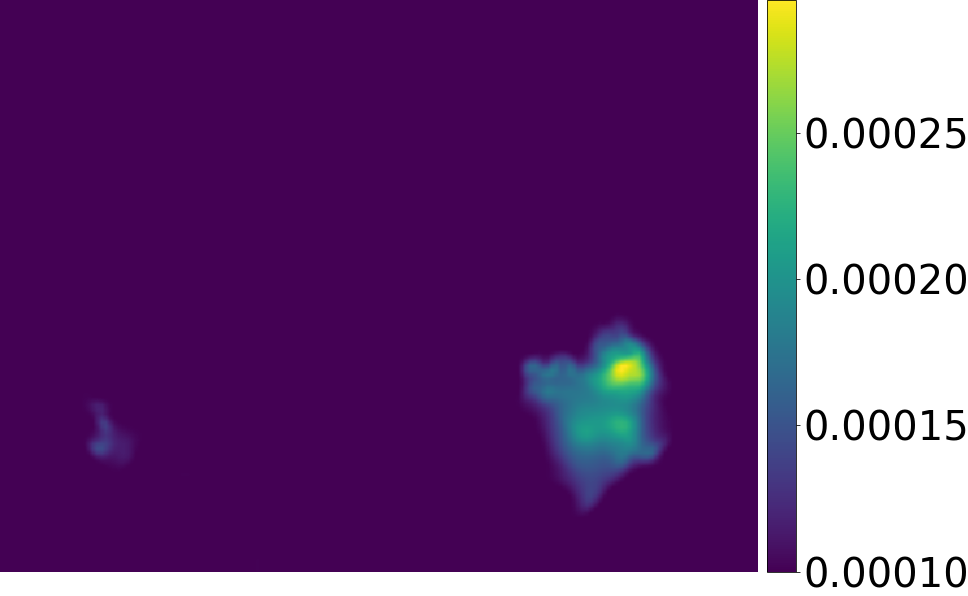}}\hspace*{.1cm}
    \subfloat[\tiny{t = 50 days}\label{fig:can_advec_pB_5}]
    {\includegraphics[width=4cm,height=3cm]{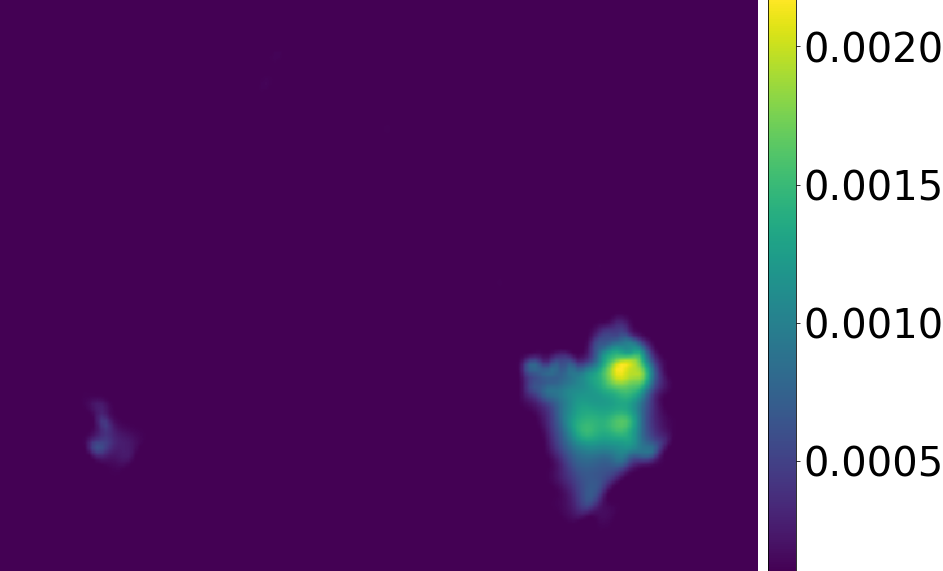}}\hspace*{.1cm}
    \subfloat[\tiny{t = 70 days}\label{fig:can_advec_pB_7}]
    {\includegraphics[width=4cm,height=3cm]{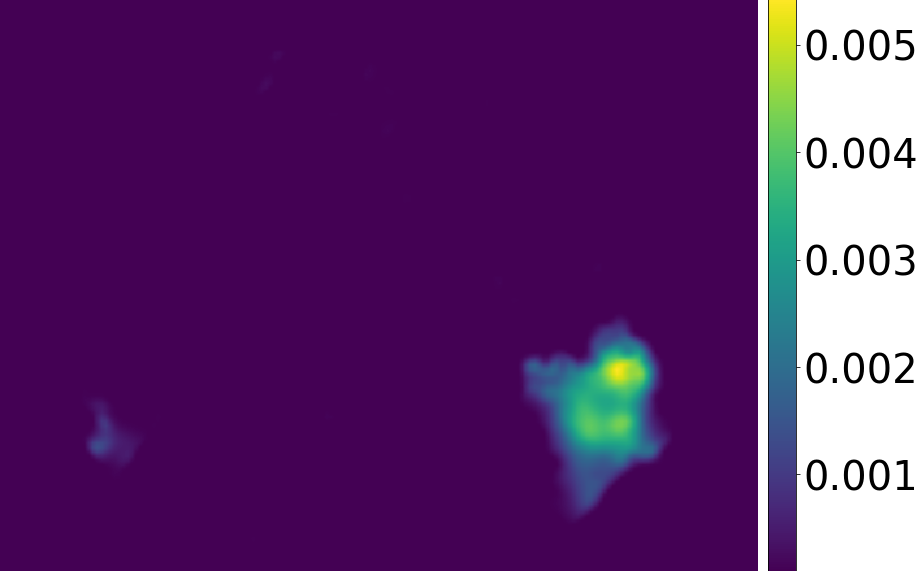}}\hspace*{.1cm}
    \subfloat[\tiny{t = 90 days}\label{fig:can_advec_pB_9}]
    {\includegraphics[width=4cm,height=3cm]{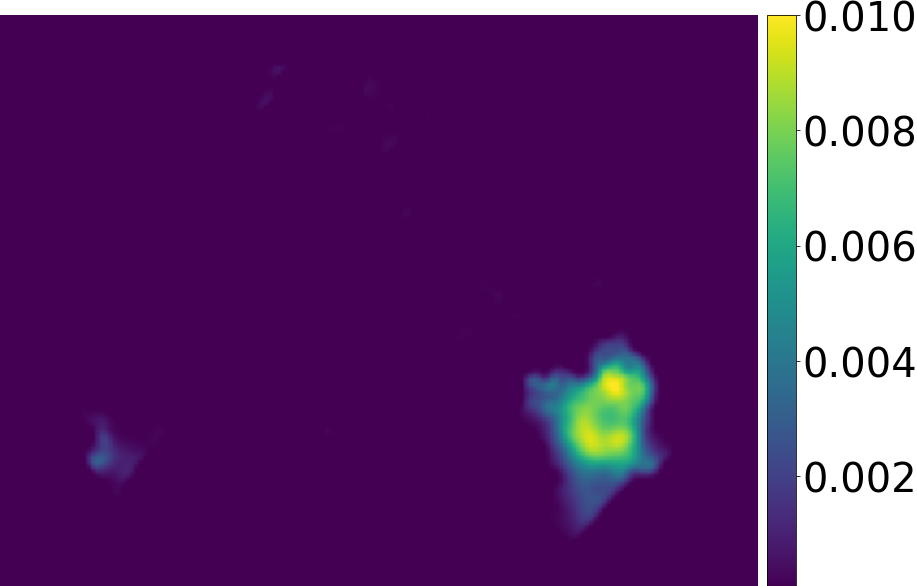}}\\

    \centering 
    \subfloat[\tiny{t = 20 days}\label{fig:can_phtaxis_pB_2}]
    {\includegraphics[width=4cm,height=3cm]{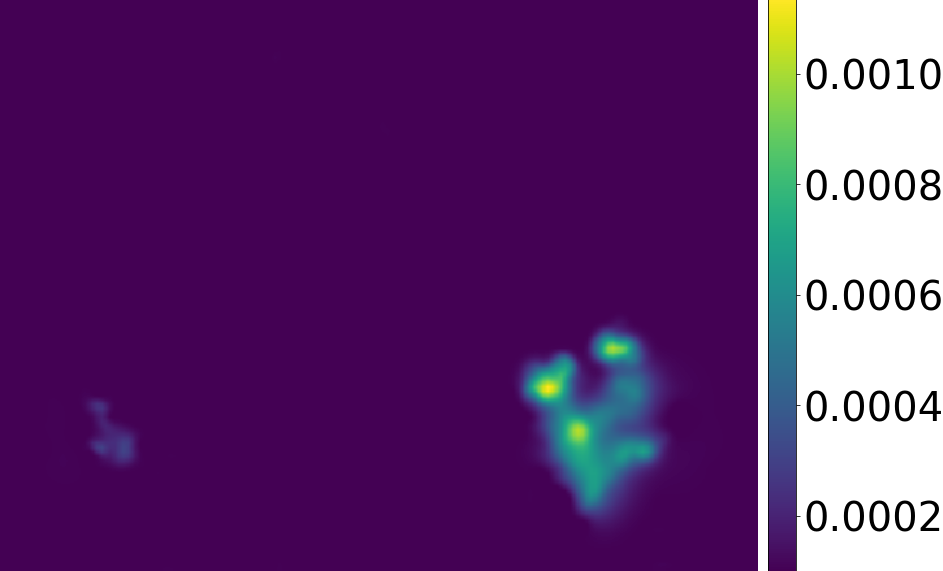}}\hspace*{.1cm}
    \subfloat[\tiny{t = 50 days}\label{fig:can_phtaxis_pB_5}]
    {\includegraphics[width=4cm,height=3cm]{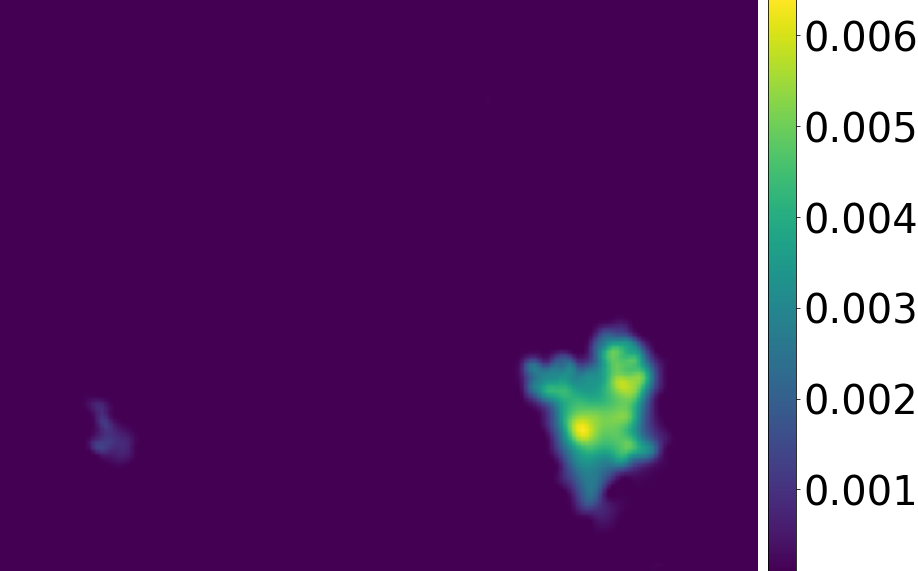}}\hspace*{.1cm}
    \subfloat[\tiny{t = 70 days}\label{fig:can_phtaxis_pB_7}]
    {\includegraphics[width=4cm,height=3cm]{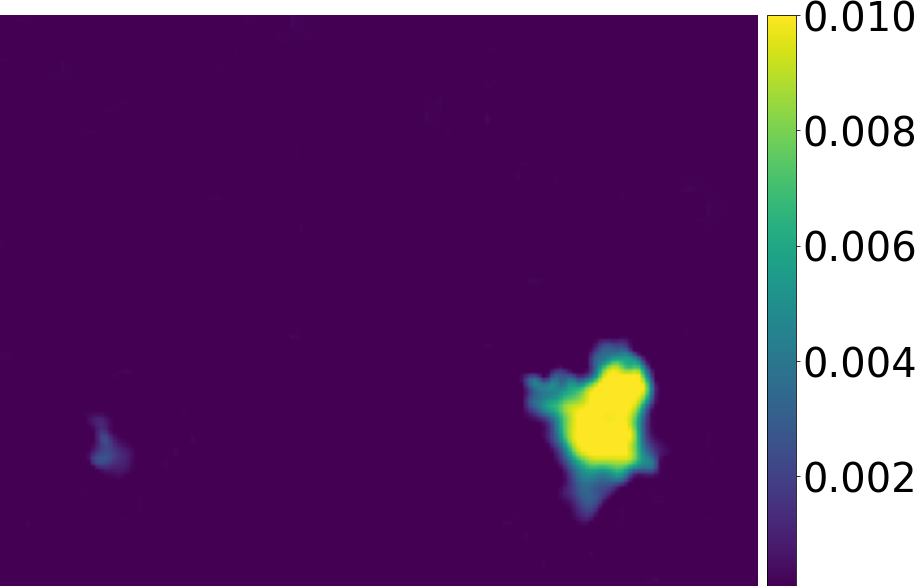}}\hspace*{.1cm}
    \subfloat[\tiny{t = 90 days}\label{fig:can_phtaxis_pB_9}]
    {\includegraphics[width=4cm,height=3cm]{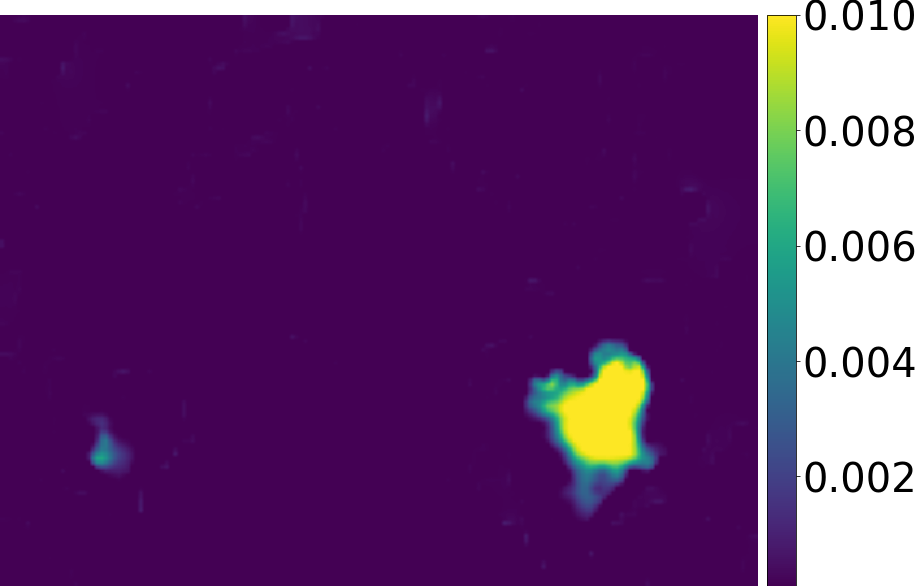}}\\
    \caption{Evolution of tumor motility parameters. The columns are sorted according to time (in days). The leftmost column is the tumor state at t=20, the middle left column is for t=50, the middle right one for t=70, and the rightmost column is for t=90. The rows represent different parameter functions of the model \eqref{eq:absPseudo}. Arranged from top to bottom these are: the diffusion coefficient {$\sigma $} (1st row), advection coefficient $\kappa$ (2nd row) and pH-taxis coefficient $\delta$ (bottom row).  \label{fig:can_phtaxis_pB}}
\end{figure}

\paragraph{Pseudopalisade specific results:} In this section we evaluate the ability of the model to replicate different realistic pseudopalisade patterns $\obs_n$ as shown in Figure \ref{fig:targetPPDs}. We also consider other target patterns which are displayed in pre-processed form.  The generated optimal final density distribution of tumor cells corresponding to $\obs_n$ is as shown in Figures \ref{fig:patFormA}-\ref{fig:patFormE}. Based on these outputs we can observe the following:
\begin{enumerate}
    \item The optimization algorithm is able to accurately generate/recreate the target pseudopalisade pattern. This is evident by looking at the fourth column of the Figures \ref{fig:patFormA}-\ref{fig:patFormE}, where we depicted the $L^2$-norm of the error i.e. the difference between the estimated final tumor density and the required target density. The difference is mainly in the 2nd decimal and only for very small volume fractions. Also, it is important to notice that the estimated final tumor density is much smoother when compared to the required target density. This is a consequence of the wellposedness Theorem \ref{thm:existPseudo} which stipulates that the solution $\xhat $ lies in $\U$ with $\xhat(t) \in \hott$ for every $t \in [0,T]$.
    \item Along with the tumor density, the algorithm also estimates the acid distribution. This is depicted in the third column of Figures \ref{fig:patFormA}-\ref{fig:patFormE}. Based on this we see that, at the regions of higher cancer cell density, in particular at areas of pseudopalisade formation, the acid concentration is relatively low compared to that of the surroundings. This supports the common hypothesis that the center of a pseudopalisade is a necrotic region, with  poor acid removal mechanisms, which results in relatively low pH. 
    \item Another interesting observation is that certain localized areas in the interior region encompassed by a pseudopalisade structure, show relatively high acid concentration. This suggests that these localized regions were the sites of high tumor activity which was likely to be a consequence of high glycolysis  followed by growth and expansion of the tumor front. As a consequence of the excess acid produced and the expanding tumor periphery, acid gets accumulated, primarily in areas of poor vasculature such as the core of the pseudopalisade.
\end{enumerate}
In order to get a deeper understanding about the formation process of pseudopalisade structures, we look at the estimated model parameter function $\th$. We do so for a fixed target pattern, namely for Pattern-B (see Figures \ref{fig:patternB} and \ref{fig:patFormB}). The obtained model parameter functions are depicted in Figures \ref{fig:can_grow_pB}-\ref{fig:can_phtaxis_pB}. Based on the dynamics of the parameters itself we can infer the following:
\begin{enumerate}
    \item The tumor growth rate $\mu$ and acid expulsion rate $\alpha$ resemble structurally very much the target  pattern. Initially the growth and expulsion rates are relatively high and later near the end time, when the tumor distribution is approaching the required pseudopalisade pattern, these rates stabilize. Moreover, it can be observed that $\mu$ and $\alpha$ are positively correlated i.e. higher $\mu$ implies higher $\alpha$, at least at the beginning of tumor evolution. This positive feedback of growth rate and acid buffering indicates the presence of (reminiscent) vasculature and the supporting microenvironment to facilitate removal excess acid. Disruption of this vital supporting element during tumor progression results in the formation of necrotic regions like those appearing toward the end time.
    \item The acid production rate $\beta$ is higher mainly in the regions where there is less tumor density. These regions of higher $\beta$ mainly happen to be the area of tissue necrosis, leading to the eventual accumulation of acid. These above average acid production rates could be attributed to the neoplastic transformation in those regions where excessive glycolysis takes place, in order to fulfill the energy requirements for proliferation. 
    \item Looking at the migration parameters we see that the diffusion coefficient $\sigma$ is lower in the sparsely populated tumor regions which are the main candidate areas for the necrotic core formation. As the tumor progresses, the diffusion coefficient mainly homogenizes and can be approximated to be spatially constant.
    \item The advection coefficient $\kappa$ is initially more pronounced at the outer margin of the necrotic core which later progresses to the inner region of the core. This indicates that harsher region is likely to generate an aggressive stimulus making the tumor cells more mobile.  
    \item The pH-taxis coefficient $\delta$ is mainly  {at the outer edges of areas with acid accumulation, which corresponds to necrosis. Thus, pH-taxis seems to act} mainly at the interfacing/intersecting layer of high-density and low-density regions of tumor. This suggests that acidity facilitates travelling-wave like behavior of tumor-host interface. 
    \item Finally, the clear distinction between regions where the taxis and growth parameters are dominant supports the hypothesis of grow-or-go dichotomy in glioma tumor progression \cite{Hoering2012}.
\end{enumerate}

\newcommand{\bxihat}{\hat\bxi}

\section{Pattern synthesis and disruption} \label{sec:thsynth}
\begin{figure}[!htbp]
    \centering
    \subfloat[\tiny{Disrupted pattern-A}\label{fig:therapyA_state}]
    {\includegraphics[width=4cm,height=3cm]{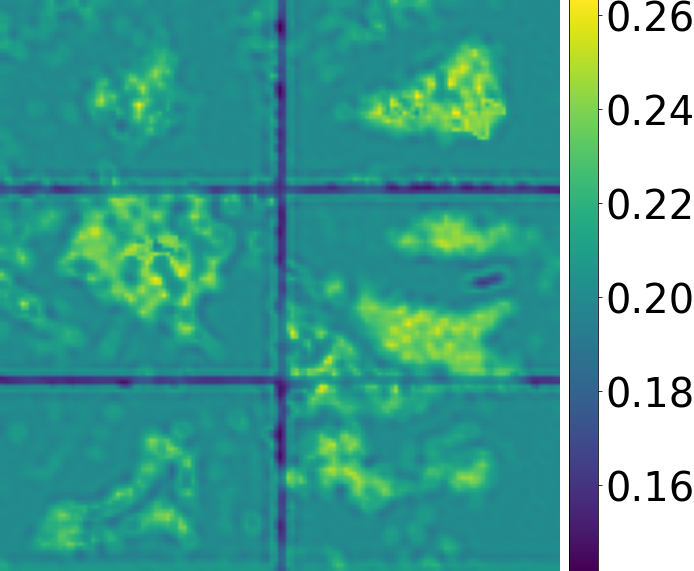}} \hspace*{.01cm}
    \subfloat[\tiny{Disrupted pattern-B}\label{fig:therapyB_state}]{
  \includegraphics[width=4cm,height=3cm]{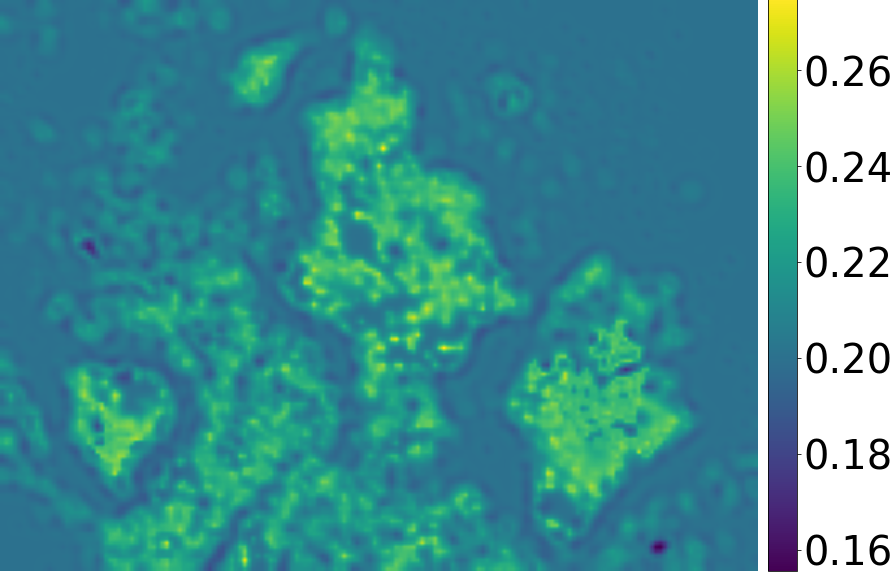}}  \hspace*{.1cm}
    \subfloat[\tiny{Disrupted pattern-C}\label{fig:therapyC_state}]
    {\includegraphics[width=4cm,height=3cm]{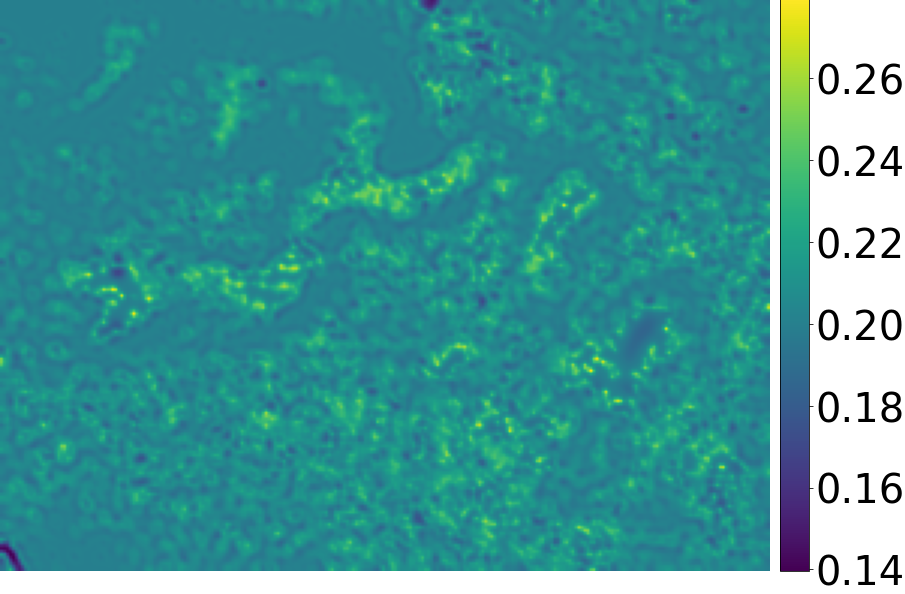}} \hspace*{.1cm}
    \subfloat[\tiny{Disrupted pattern-E}\label{fig:therapyE_state}]
    {\includegraphics[width=4cm,height=3cm]{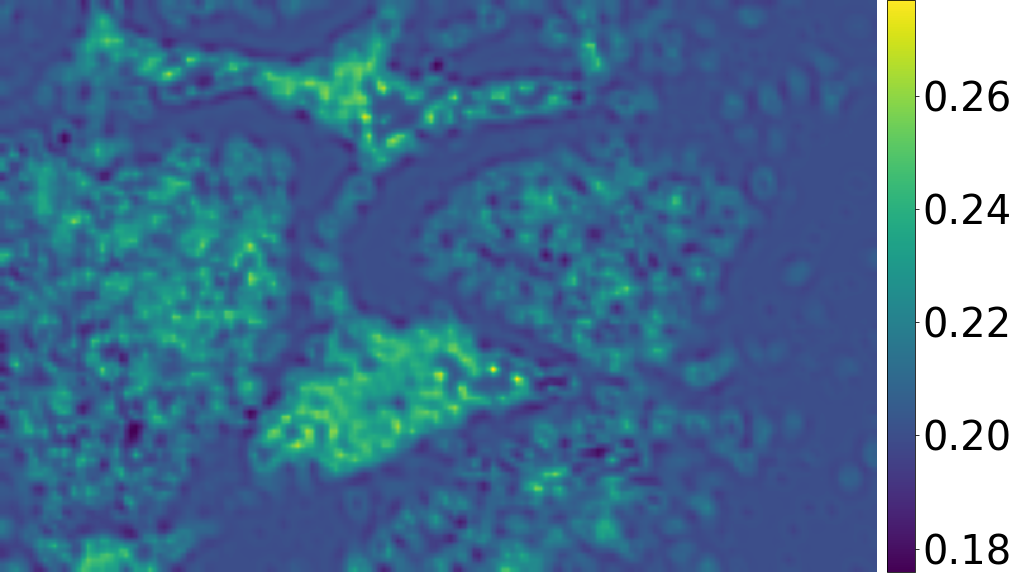}}
    \caption{Final pseudopalisade patterns resulting after applying the pattern neutralizing function $\bxihat$ with the corresponding pattern specific estimated model parameters $\hat \th$. \label{fig:therapyStatePats}}
\end{figure}

\begin{figure}[!htbp]
    \centering
    \subfloat[\tiny{Pattern-A}\label{fig:therapyA_instate}]
    {\includegraphics[width=4cm,height=3cm]{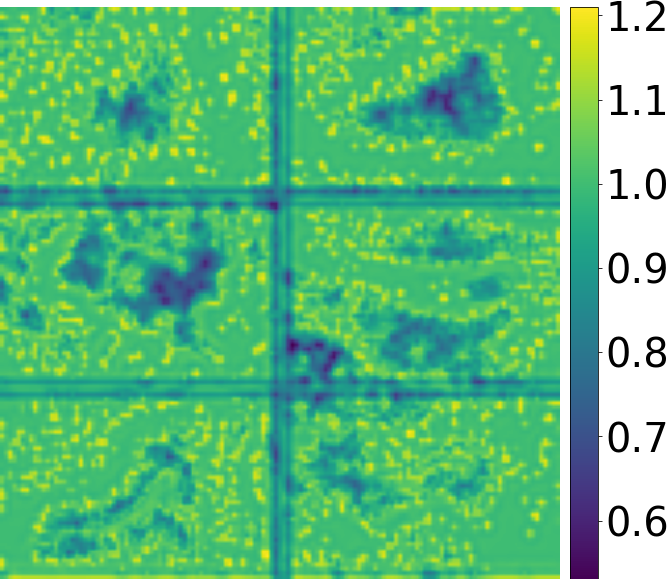}} \hspace*{.01cm}
    \subfloat[\tiny{Pattern-B}\label{fig:therapyB_instate}]{
  \includegraphics[width=4cm,height=3cm]{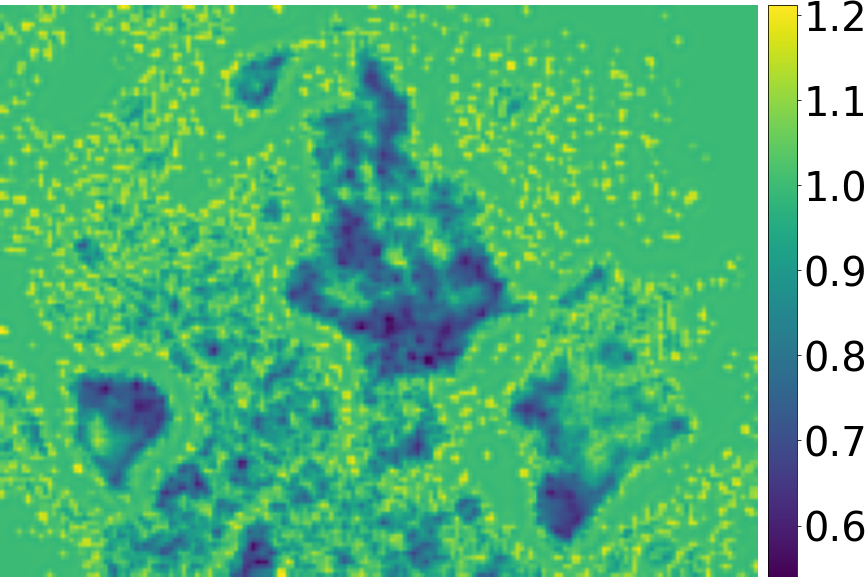}}  \hspace*{.1cm}
    \subfloat[\tiny{Pattern-C}\label{fig:therapyC_instate}]
    {\includegraphics[width=4cm,height=3cm]{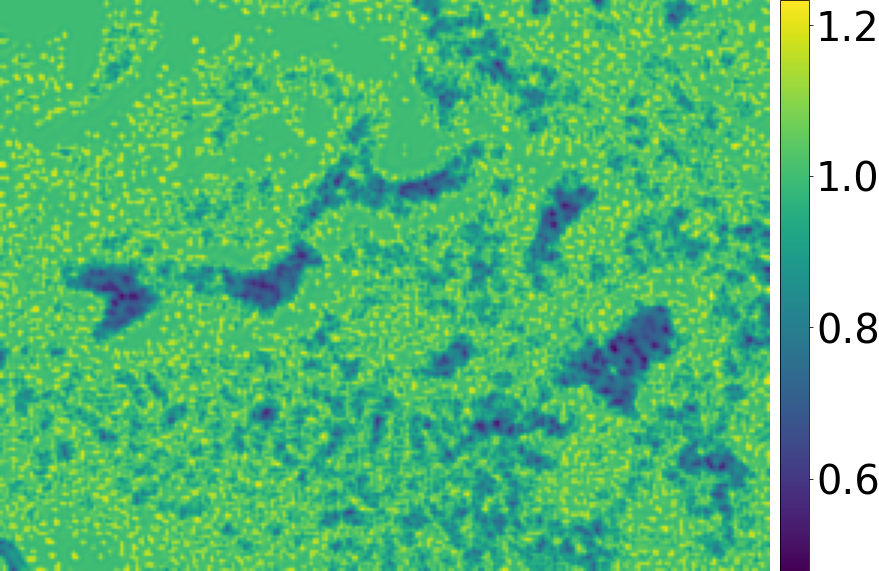}} \hspace*{.1cm}
    \subfloat[\tiny{Pattern-E}\label{fig:therapyE_instate}]
    {\includegraphics[width=4cm,height=3cm]{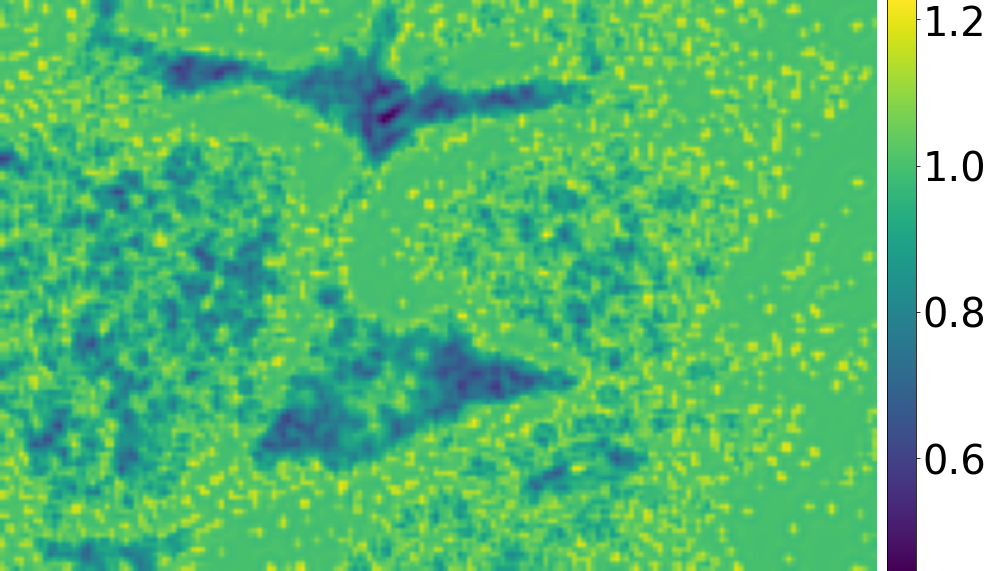}}
    \caption{Final proton distribution corresponding to the pseudopalisade patterns in Figure \ref{fig:therapyStatePats} \label{fig:therapyInStatePats}}
\end{figure}

In this section we consider {an} application of the OCP \eqref{eq:absPseudo}, {aiming to shed light on pattern dynamics under various influences.} One of the main advantages of data-based parameter estimation is its {potential} use in patient-specific therapy design. The diagnostic {histological samples} obtained from a patient can be used to estimate the model parameters which approximately characterize the {(microscopic)} dynamics of GBM progression for that specific patient. Based on this, one can then design or hypothesize different intervention mechanisms that {impair} the development of GBM. One such plausible {way would be to understand how pseudopalisade pattern formation can be disrupted or reversed. To incorporate such a feature,} we adapt the above model as follows:
\begin{subequations}
\label{eq:therapyPPD}
\begin{align}
\partial_t u_1 - \nabla \cdot (\sigma \nabla u_1 + u_1 \nabla \kappa ) &= \nabla \cdot (\delta \: u_1  \nabla u_2)  + \mu  f_1(u_1, u_2) - \xi_1 u_1 \label{eq:thCan}\\ 
\partial_t u_2 - \lap u_2 + \alpha u_2 &=  \beta f_2(u_1,u_2)  + \xi_2 u_2 \label{eq:thHe}\\ 
(\sigma \nabla u_1 + u_1 \nabla \kappa) \cdot \hat n = 0,  \quad \nabla \cdot u_2 &= 0, \quad u_1(0) = \u_{1_0}, \quad u_2(0) = \u_{2_0} \notag 
\end{align}
\end{subequations}
Here, $\bxi = (\xi_1,\xi_2) \in L^2(I; V^{\times 2})$ with $\xi_1 \ge 0$ is a disturbance term that models the disruption of pattern formation mechanisms. Consequently, $\bxi$ is responsible for neutralization or renormalization of the tumor microenvironment, which altogether impedes tumor development. Motivated by this we refer to $\bxi$ as a pattern neutralizing function. Having estimated the model parameters $\th$ for different target patterns (see pictures in Figure \ref{fig:targetPPDs}) we can now ask what kind of external signal is needed to revert or neutralize the cancerous microenvironment. This entails solving the following modified optimization problem: 
\begin{align}
\label{eq:therapyOCP}
    \begin{split}
        \u^*, \bxi^* &= \argmin{\u, \bxi} J(\u,\th,\bxi) \quad {\rm{ s.t. }} \quad G(\u; \th, \bxi) = 0,
    \end{split} 
\end{align}
where $J$ is the quadratic cost function analogous to \eqref{eq:cost} and $G(\u;\th,\bxi)=0$ represents the state equation \eqref{eq:therapyPPD} in abstract form. 
Now, given a non-cancerous or a neutral tissue pattern as the target state, we can solve for the optimal pattern neutralizing function $\bxi^*$. To be more precise, let $\th_{O}$ be the optimal model parameter vector associated to the target pattern Pattern-$O$ where $O \in \{A, B, C, E\}$ i.e. one of the target patterns depicted in Figure \ref{fig:targetPPDs}. Then given some starting value $\u_0$ (representing a non-cancerous initial state) and fixing the model parameters $\th_{O}$, we apply Algorithm \ref{algo:pgd} to determine $\bxihat$  that can counteract the effects of $\th_{O}$ and result in a target state which corresponds to a neutral non-cancerous tissue. Figures \ref{fig:therapyStatePats} - \ref{fig:therapyInStatePats} depict the result of the application of the neutralizing function $\bxi^*$ that is able to counteract the effects of the parameters responsible for generation of the specific pseudopalisade patterns PatternB-PatternE (see Figure \ref{fig:targetPPDs}) and PatternA (see \cite{PatA} Fig. 7.). 

\begin{figure}[!htbp]
    \centering
    \subfloat[\tiny{Disrupted pattern-A}\label{fig:pHtherapyA_state}]
    {\includegraphics[width=4cm,height=3cm]{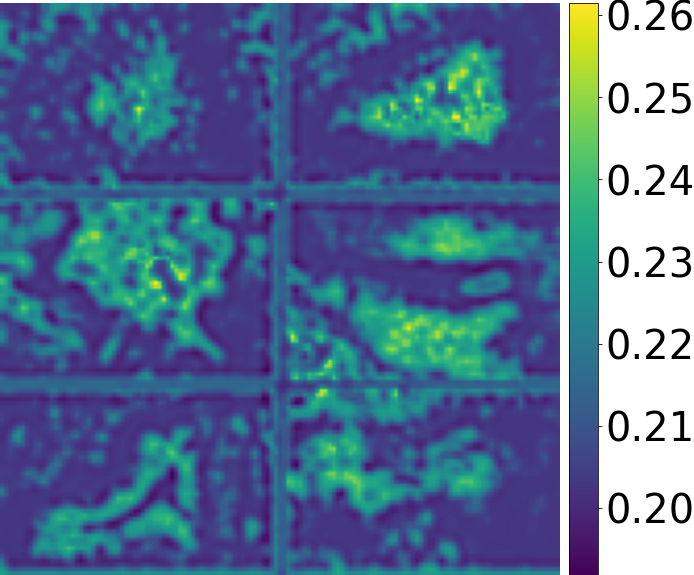}} \hspace*{.01cm}
    \subfloat[\tiny{Disrupted pattern-B}\label{fig:pHtherapyB_state}]{
  \includegraphics[width=4cm,height=3cm]{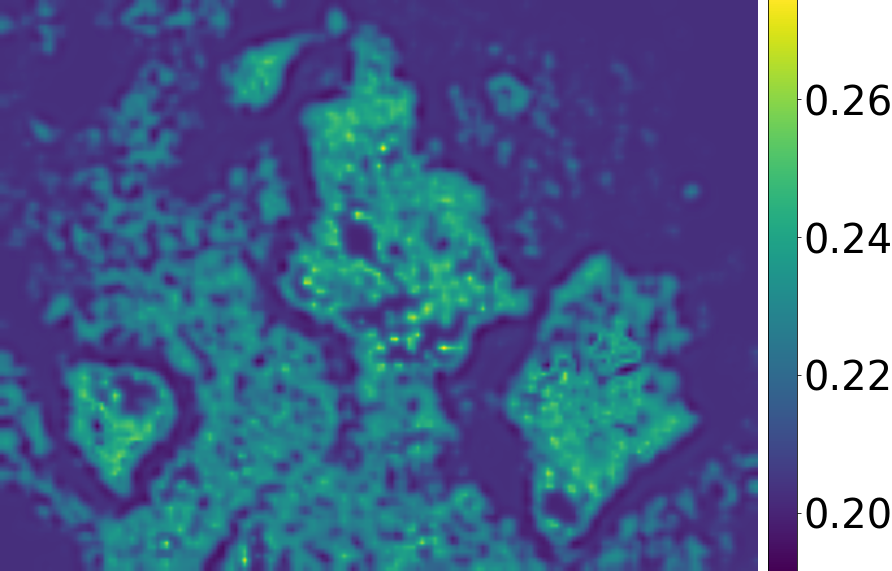}}  \hspace*{.1cm}
    \subfloat[\tiny{Disrupted pattern-C}\label{fig:pHtherapyC_state}]
    {\includegraphics[width=4cm,height=3cm]{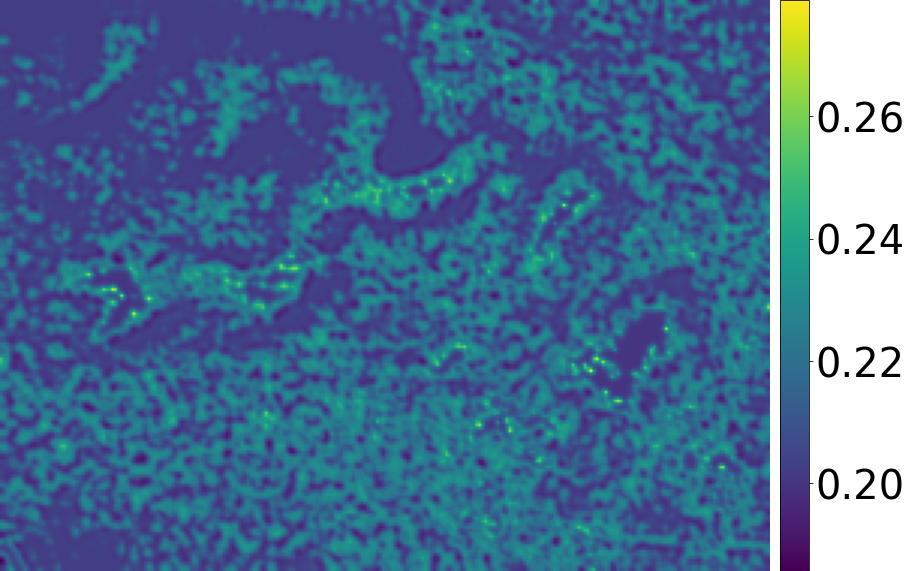}} \hspace*{.1cm}
    \subfloat[\tiny{Disrupted pattern-E}\label{fig:pHtherapyE_state}]
    {\includegraphics[width=4cm,height=3cm]{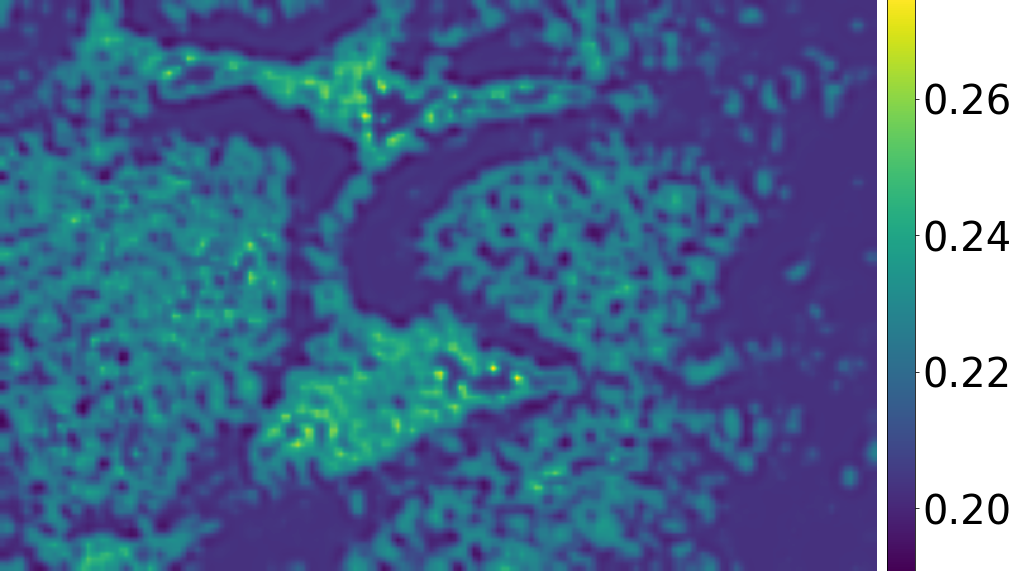}}
    \caption{Final pseudopalisade patterns resulting after applying the pH neutralizing function $\hat \xi_2$ with the corresponding pattern specific estimated model parameters $\hat \th$.  \label{fig:pHtherapyStatePats}}
\end{figure}

\begin{figure}[!htbp]
    \centering 
    \subfloat[\tiny{ Pattern-A}\label{fig:pHtherapyA_instate}]
    {\includegraphics[width=4cm,height=3cm]{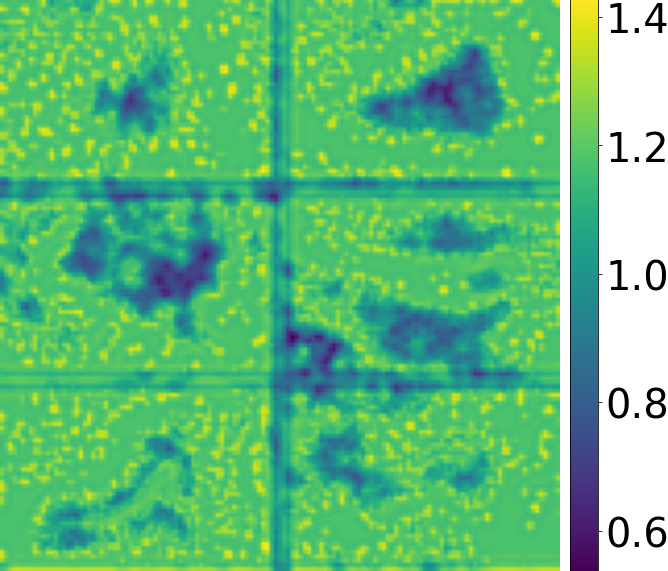}} \hspace*{.01cm}
    \subfloat[\tiny{Pattern-B}\label{fig:pHtherapyB_instate}]{
  \includegraphics[width=4cm,height=3cm]{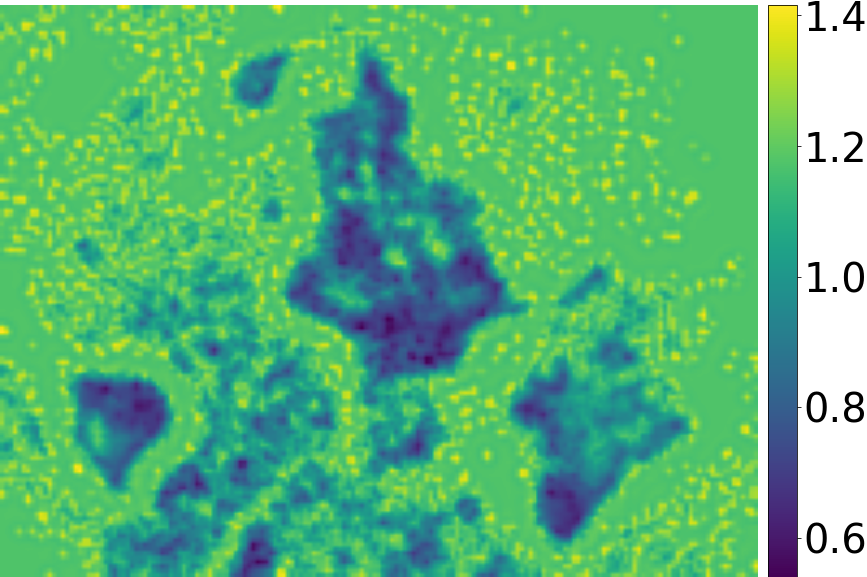}}  \hspace*{.1cm}
    \subfloat[\tiny{Pattern-C}\label{fig:pHtherapyC_instate}]
    {\includegraphics[width=4cm,height=3cm]{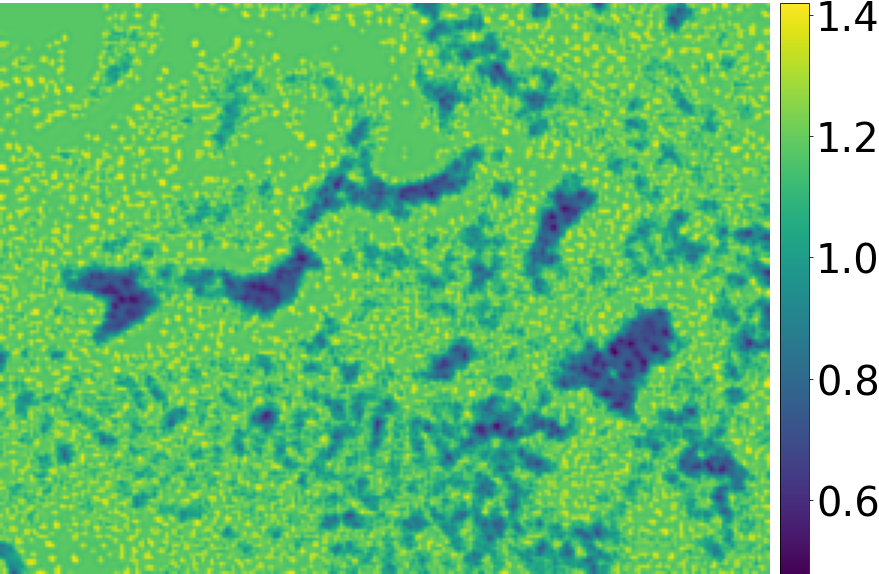}} \hspace*{.1cm}
    \subfloat[\tiny{pattern-E}\label{fig:pHtherapyE_instate}]
    {\includegraphics[width=4cm,height=3cm]{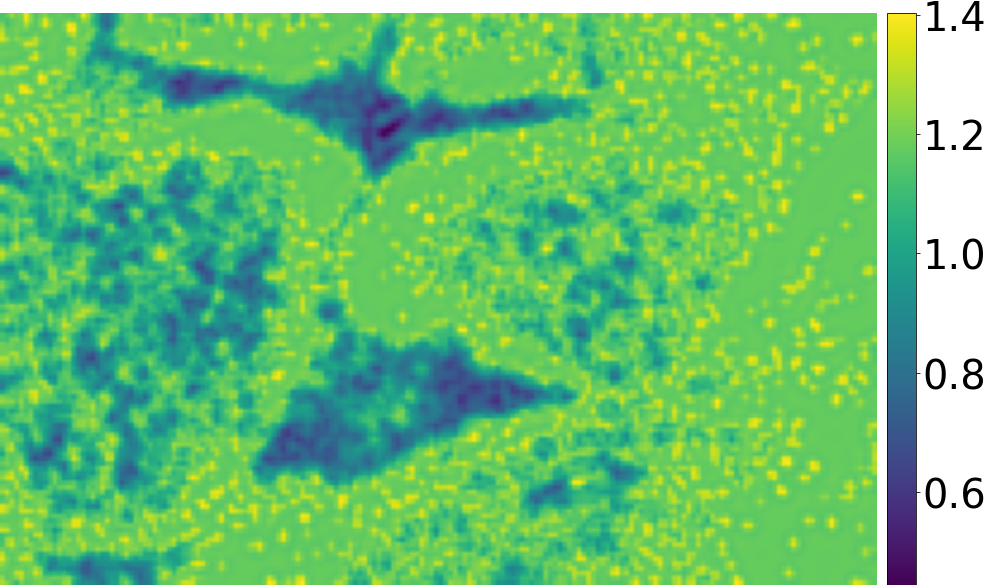}}
    \caption{Final proton distribution corresponding to the pseudopalisade patterns in Figure \ref{fig:pHtherapyStatePats}  \label{fig:pHtherapyInStatePats}}
\end{figure}

\noindent
Alternatively, instead of applying the pattern neutralizing function that directly acts on the cancer cells, one could also look for an indirect method which aims at  manipulating the micro-environment as a means to hinder tumor progression. This means that instead of having two control variables $\xi_1$ and $\xi_2$ in \eqref{eq:therapyPPD} we have only one control variable $\xi_2$ that aims to disrupt the pseudopalisade pattern by appropriately regulating the tissue acidity. Consequently, $\xi_2$ is referred to as pH neutralizing function. Following the above steps we can find the optimal $\hat \xi_2$ that alone can counteract the effects of pseudopalisade forming parameters $\thhat_0$.
Based on the results depicted in Figures \ref{fig:pHtherapyStatePats}-\ref{fig:pHtherapyInStatePats}, we can see that neutralizing the tissue acidity can also serve as an effective pseudopalisade disruption mechanism. {This seems to be in line with therapeutic approaches aiming at tumor alkalinization, see e.g., \cite{Amiri2016,YZY2020}.}

\paragraph{Synthesis of new patterns}
In this section we illustrate how new patterns can be synthesized by combining, linearly or nonlinearly, the optimal parameters $\thhat_{O}$ corresponding to different target pseudopalisade patterns $O$. We let $\thhat_{KL} := (\thhat_K + \thhat_L)/2$, where $\thhat_K$ and $\thhat_L$ are the estimated optimal model parameters for target patterns Pattern-K and Pattern-L, respectively. The patterns generated by different such combinations are depicted in Figures \ref{fig:synthStatePats}-\ref{fig:synthInStatePats}. These indicate that different complex patterns can arise by linear combination of processes responsible for generating simpler patterns.
This  {is of particular importance in view of developing possible therapy concepts,} where one can ask the question if a similar linear combination of pattern neutralizing functions can still be effective for disrupting the new pattern. That is, if $\bxihat_K$ and $\bxihat_L$ are the pattern neutralizing functions that can counteract the effects of $\thhat_K$ and $\thhat_L$ respectively, would their combination $\bxihat_{KL} := (\bxihat_K + \bxihat_L)/2$ be able to counteract the effects of $\thhat_{KL}$? Based on the obtained numerical results, as shown in Figures \ref{fig:threapySynthState}-\ref{fig:threapySynthInState}, we can conclude that the same linear combination of {pattern neutralizing functions} does indeed prove to be effective in counteracting the combined effects of the processes responsible for generating simpler pseudopalisade patterns. 

\begin{figure}[!ht]
    \centering 
    \hspace*{-.5cm}
    \subfloat[\tiny{Pattern-K}\label{fig:tempPat_K}]
    {\includegraphics[width=3cm,height=2cm]{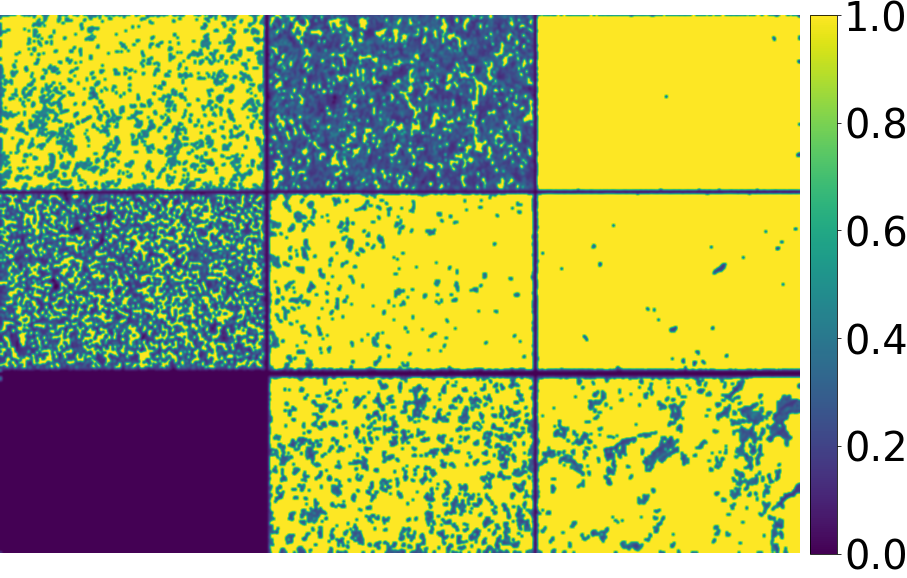}} \hspace*{.01cm}
    \subfloat[\tiny{Pattern-L}\label{fig:tempPat_L}]{
  \includegraphics[width=3cm,height=2cm]{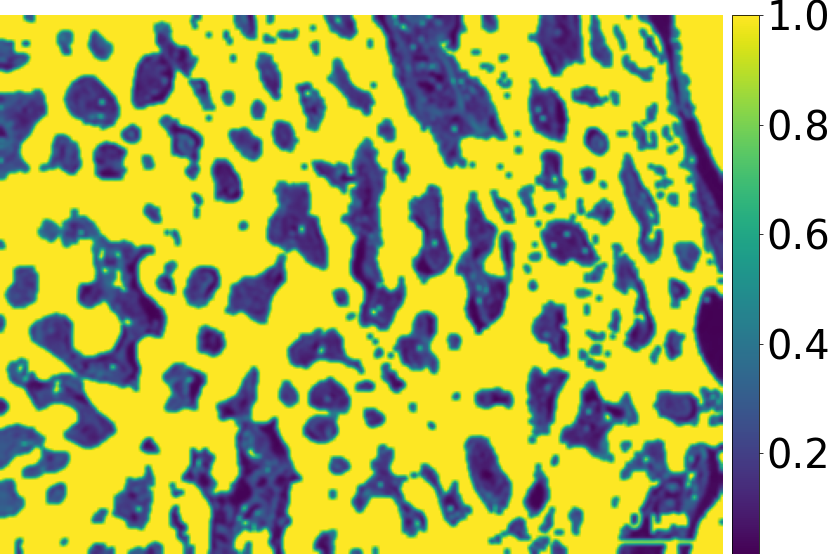}}
  \hspace*{.01cm}
  \subfloat[\tiny{Pattern-I}\label{fig:tempPat_I}]{
  \includegraphics[width=3cm,height=2cm]{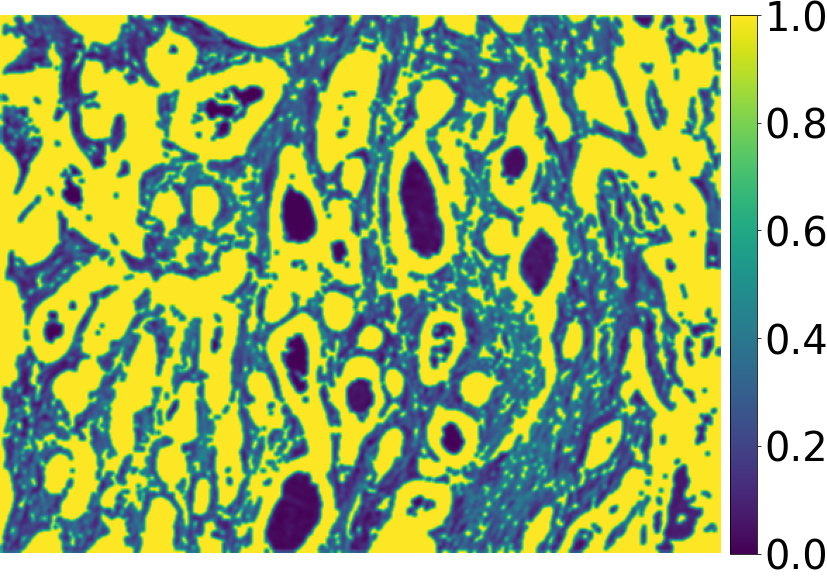}}
  \hspace*{.1cm}
  \subfloat[\tiny{Pattern-G}\label{fig:tempPat_G}]
    {\includegraphics[width=3cm,height=2cm]{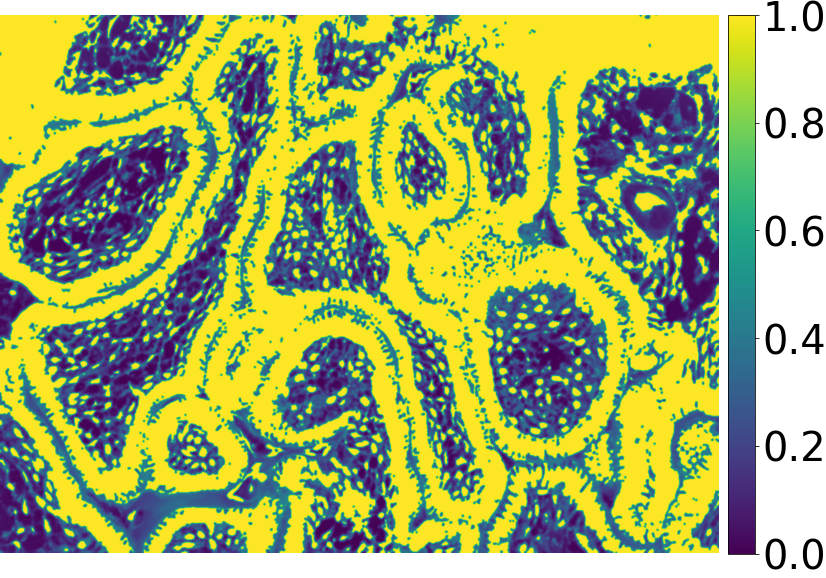}} \hspace*{.1cm}
    \subfloat[\tiny{Pattern-H}\label{fig:tempPat_H}]
    {\includegraphics[width=3cm,height=2cm]{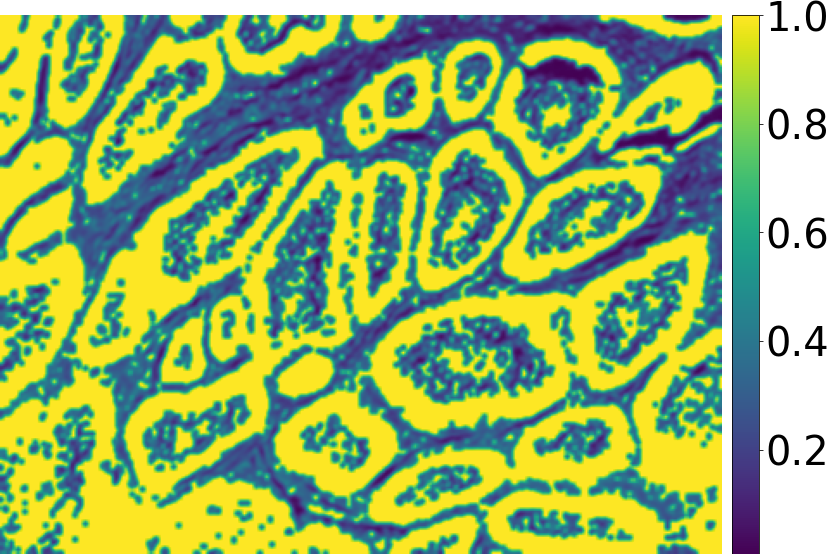}}
    \caption{Template patterns for synthesizing new ones \label{fig:templatePats}}
\end{figure}

\begin{figure}[!ht]
    \centering 
    \subfloat[\tiny{Synthesized pattern-KL}\label{fig:synthKL_state}]
    {\includegraphics[width=4cm,height=3cm]{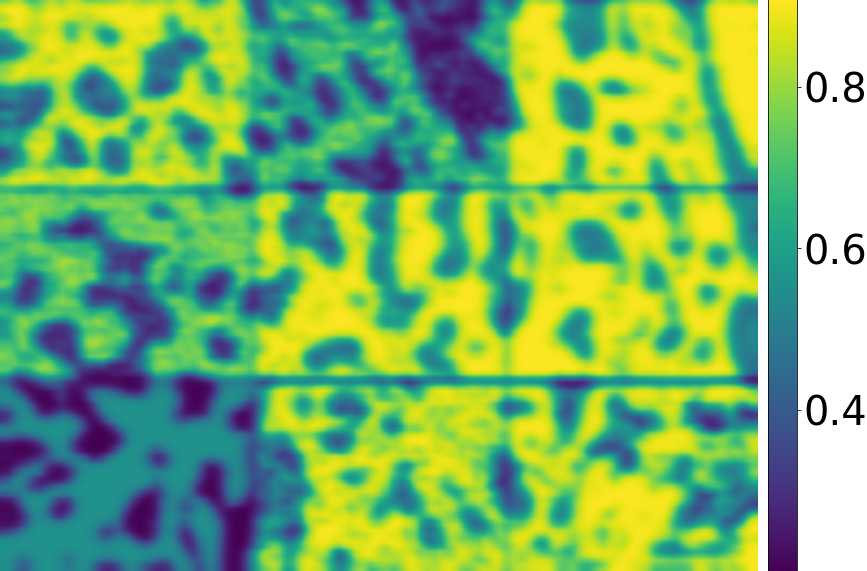}} \hspace*{.01cm}
    \subfloat[\tiny{Synthesized pattern-KI}\label{fig:synthKI_state}]{
  \includegraphics[width=4cm,height=3cm]{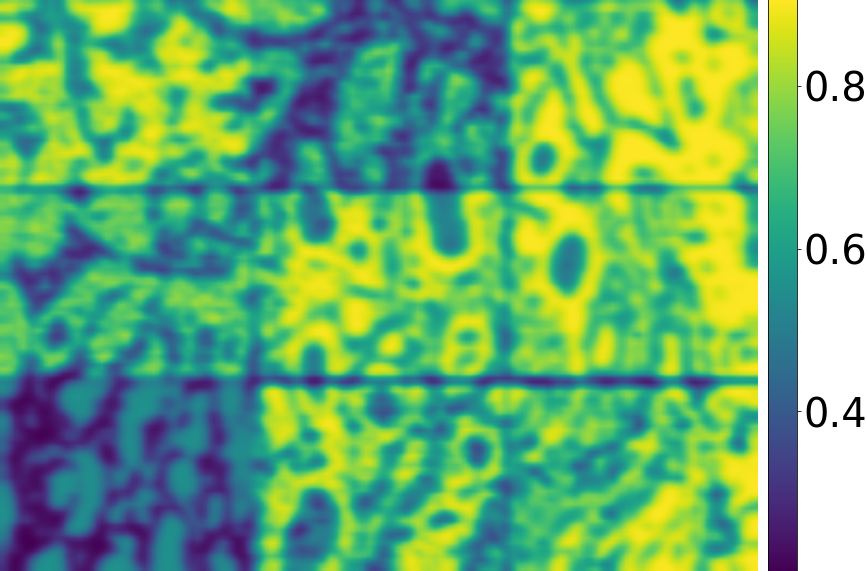}}  \hspace*{.1cm}
  \subfloat[\tiny{Synthesized pattern-KG}\label{fig:synthKG_state}]
    {\includegraphics[width=4cm,height=3cm]{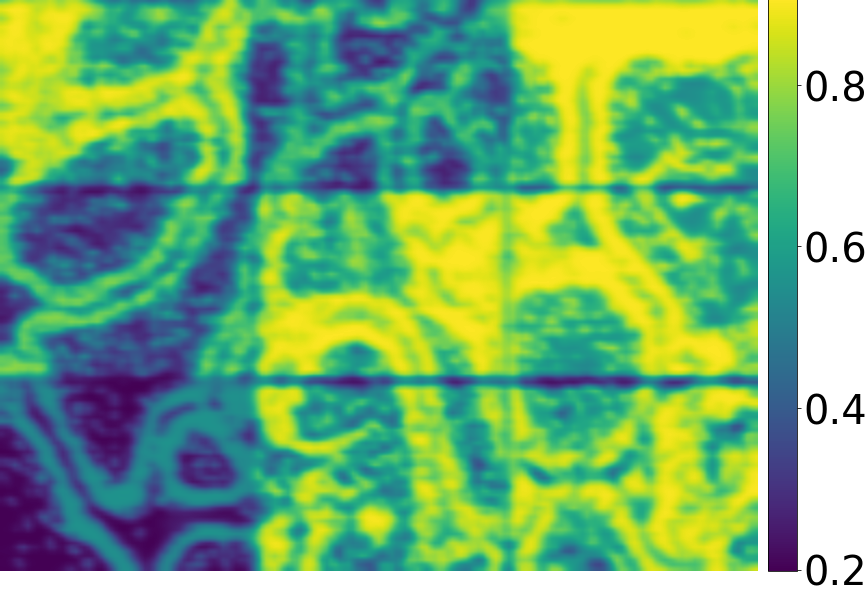}} \hspace*{.1cm}
    \subfloat[\tiny{Synthesized pattern-KH}\label{fig:synthKH_state}]
    {\includegraphics[width=4cm,height=3cm]{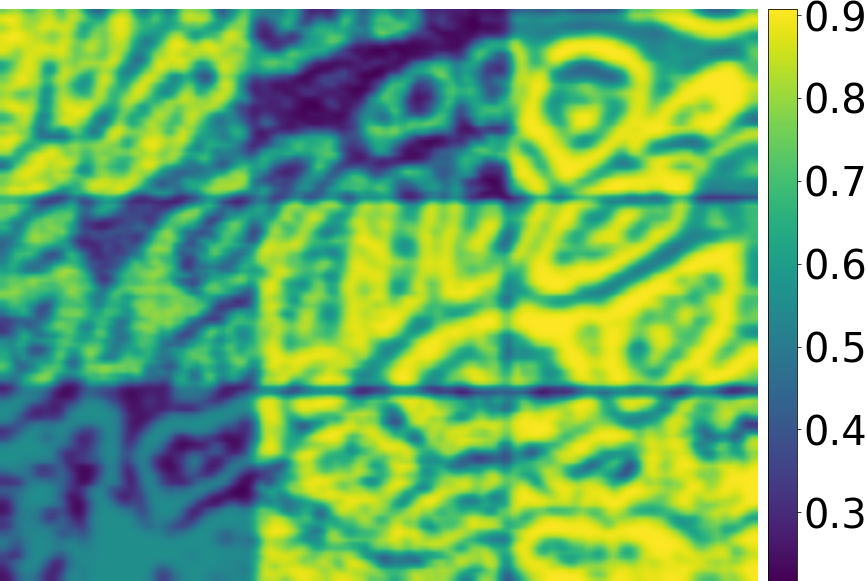}} 
    \caption{Synthesized patterns for different  combinations of the model parameters \label{fig:synthStatePats}}
\end{figure}
\begin{figure}[!ht]
    \centering 
    \subfloat[\tiny{Pattern-KL acid profile}\label{fig:synthKL_instate}]
    {\includegraphics[width=4cm,height=3cm]{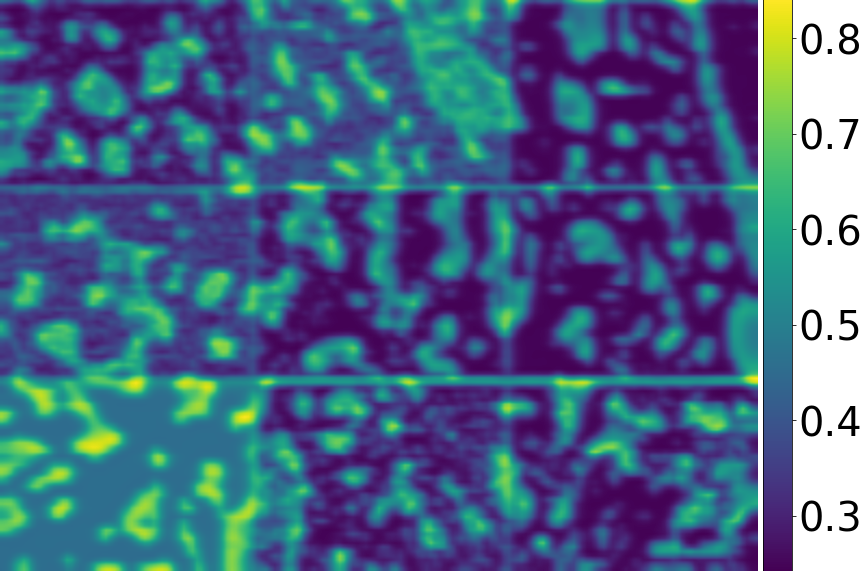}} \hspace*{.01cm}
    \subfloat[\tiny{Pattern-KI acid profile}\label{fig:synthKI_instate}]{
  \includegraphics[width=4cm,height=3cm]{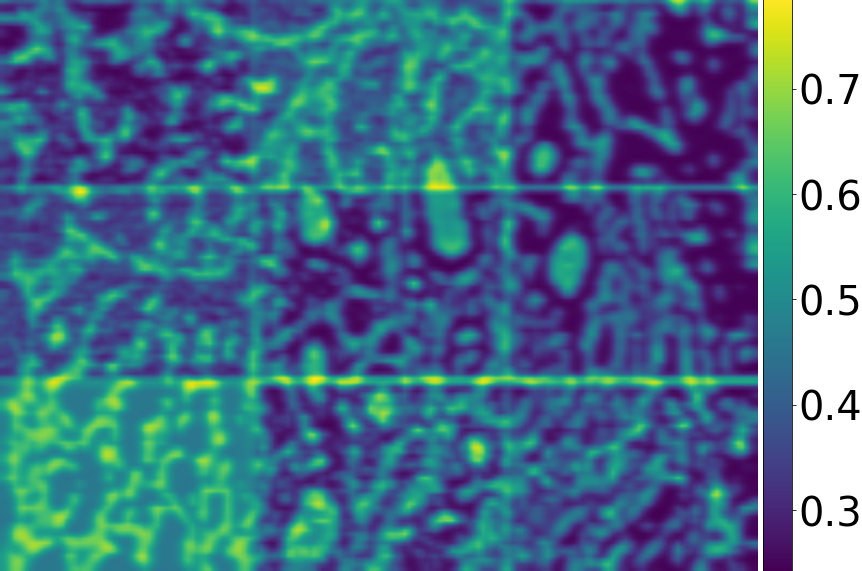}}  \hspace*{.1cm}
  \subfloat[\tiny{Pattern-KG acid profile}\label{fig:synthKG_instate}]
    {\includegraphics[width=4cm,height=3cm]{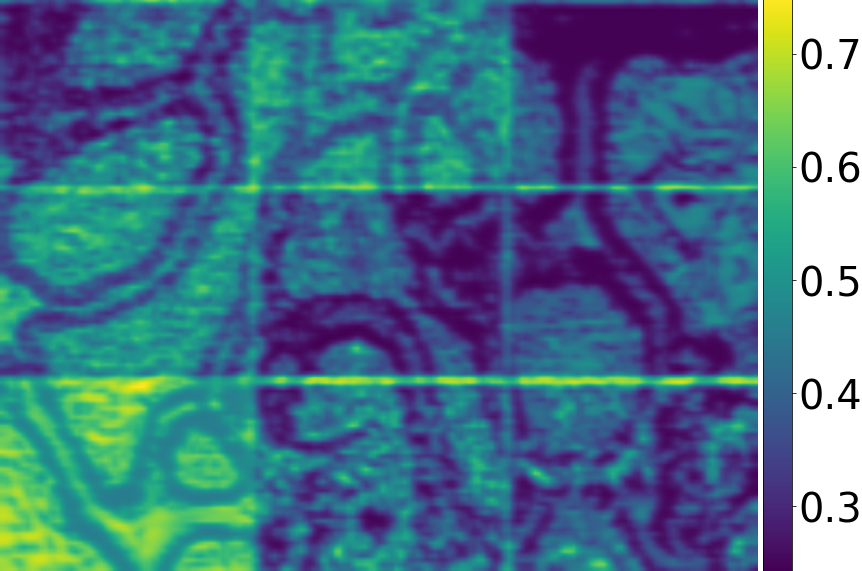}} \hspace*{.1cm}
    \subfloat[\tiny{Pattern-KH acid profile}\label{fig:synthKH_instate}]
    {\includegraphics[width=4cm,height=3cm]{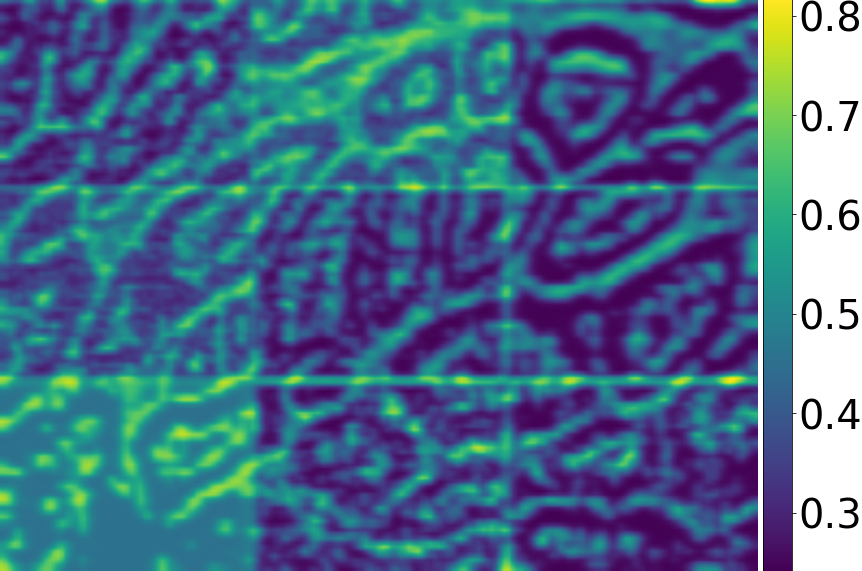}} 
    \caption{Tissue acid profiles for the synthesized patterns in Figure \ref{fig:synthStatePats} \label{fig:synthInStatePats}}
\end{figure}
 \begin{figure}[!ht]
    \centering 
    \subfloat[\tiny{Disrupted pattern-KL}\label{fig:threapy_synthKL_state}]
    {\includegraphics[width=4cm,height=3cm]{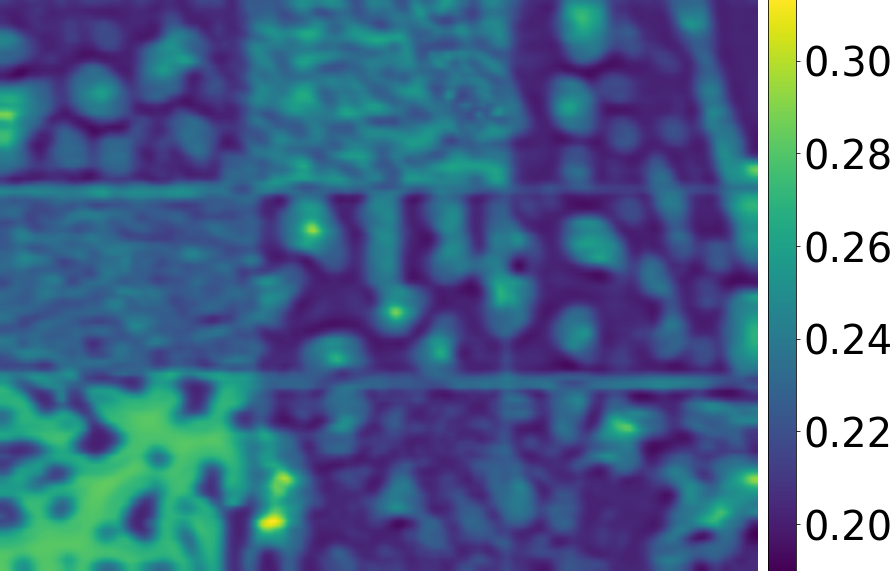}} \hspace*{.01cm}
    \subfloat[\tiny{Disrupted pattern-KI}\label{fig:threapy_synthKI_instate}]
    {\includegraphics[width=4cm,height=3cm]{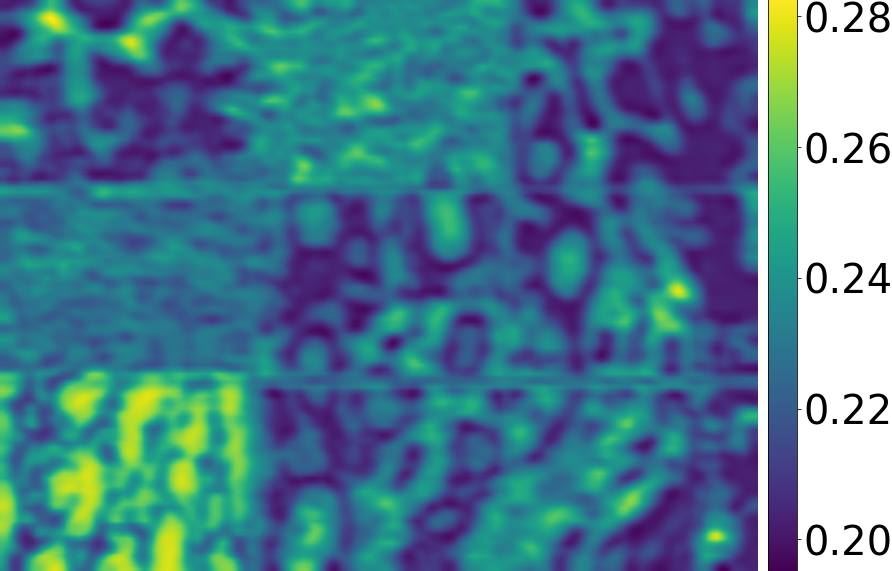}}  \hspace*{.1cm}
    \subfloat[\tiny{Disrupted pattern-KG}\label{fig:threapy_synthKG_instate}]
    {\includegraphics[width=4cm,height=3cm]{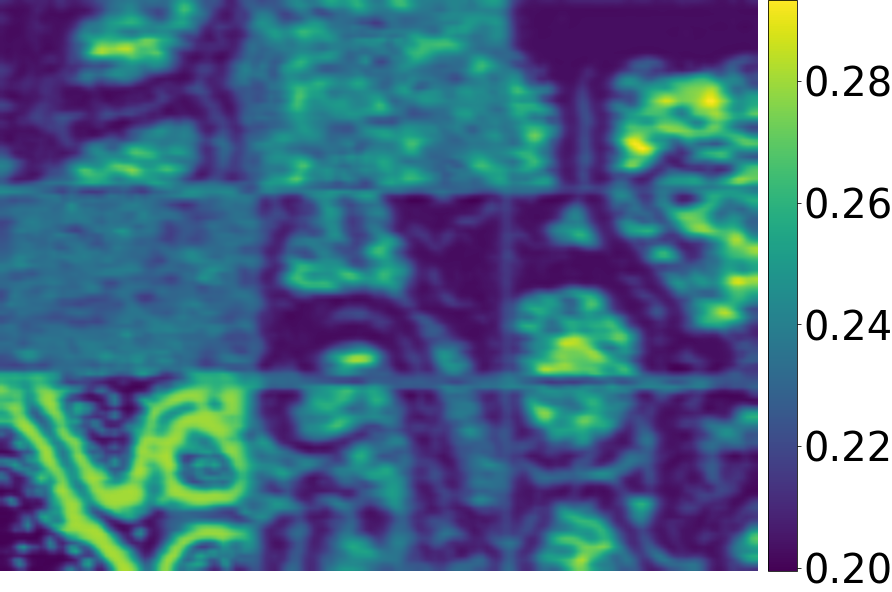}}  \hspace*{.1cm}
    \subfloat[\tiny{Disrupted pattern-KH}\label{fig:threapy_synthKH_instate}]
    {\includegraphics[width=4cm,height=3cm]{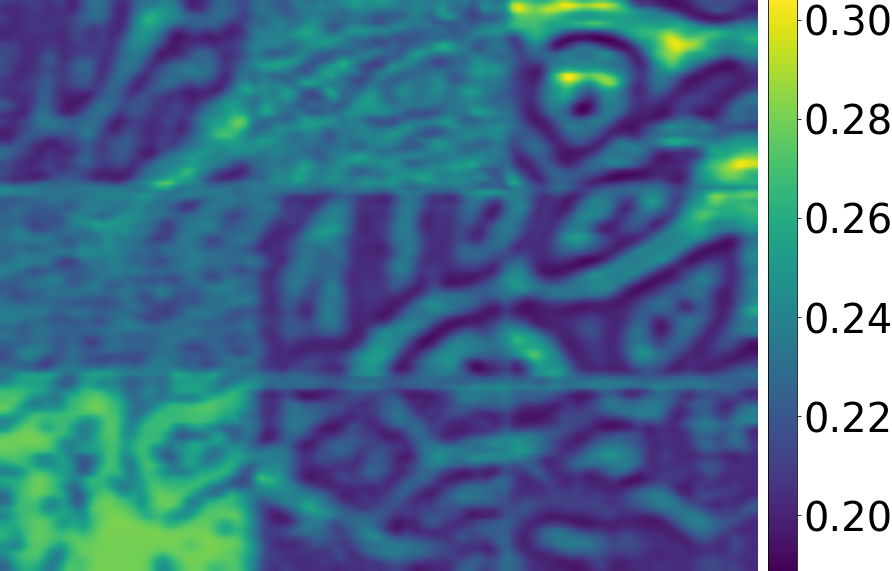}} 
    \caption{Final tumor pattern after applying optimal pH neutralizing function for the synthesized patterns in Figure \ref{fig:synthStatePats} \label{fig:threapySynthState}}
\end{figure}
\begin{figure}[!ht]
    \centering 
    \subfloat[\tiny{Pattern-KL acid profile}\label{fig:threapy_synthKL_instate}]
    {\includegraphics[width=4cm,height=3cm]{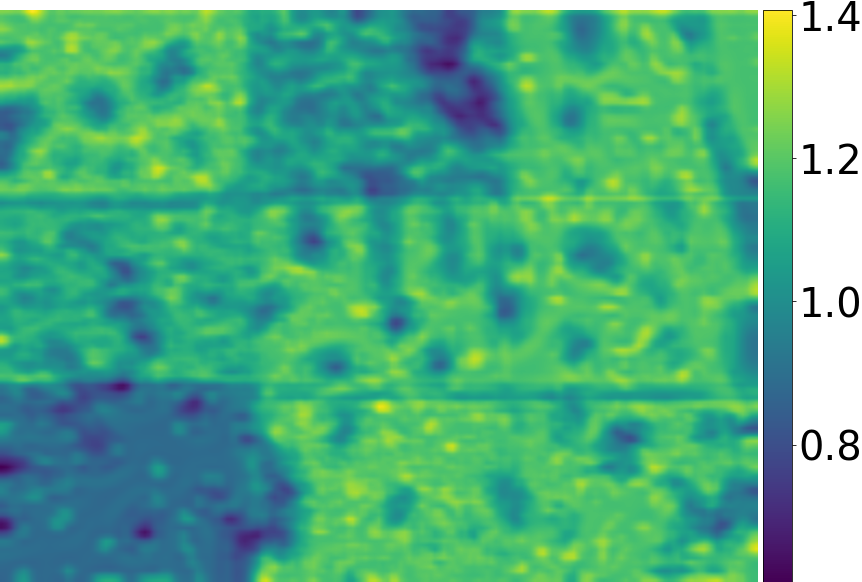}} \hspace*{.01cm}
    \subfloat[\tiny{Pattern-KI acid profile}\label{fig:threapy_synthKI_instate}]
    {\includegraphics[width=4cm,height=3cm]{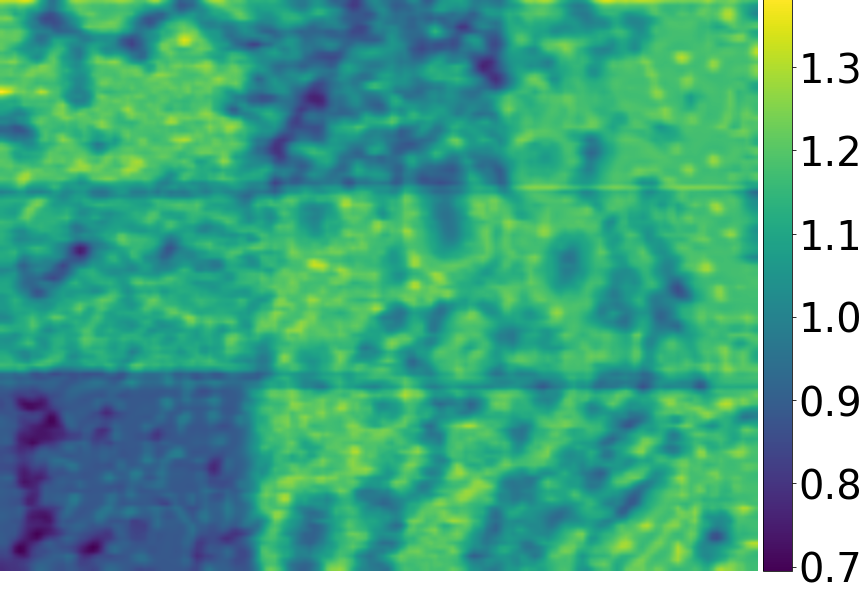}}  \hspace*{.1cm}
    \subfloat[\tiny{Pattern-KG acid profile}\label{fig:threapy_synthKG_instate}]
    {\includegraphics[width=4cm,height=3cm]{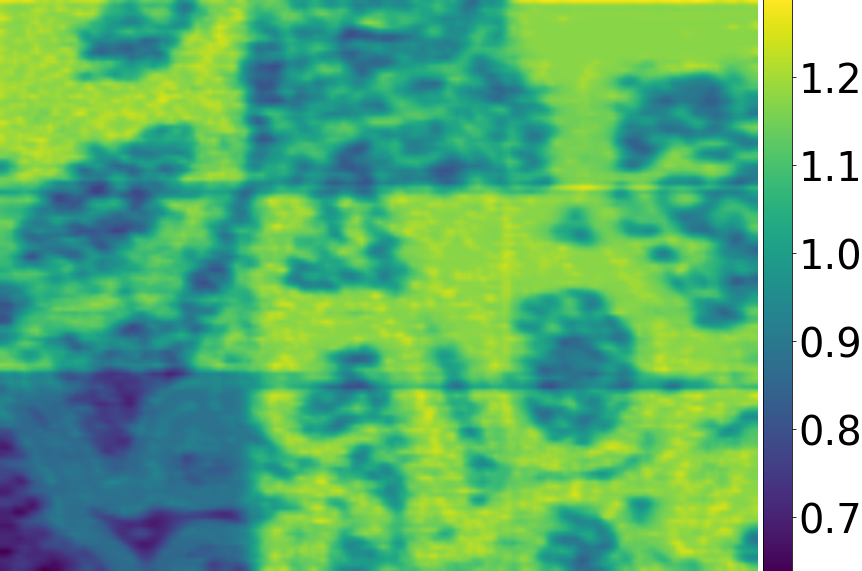}}  \hspace*{.1cm}
    \subfloat[\tiny{Pattern-KH acid profile}\label{fig:threapy_synthKH_instate}]
    {\includegraphics[width=4cm,height=3cm]{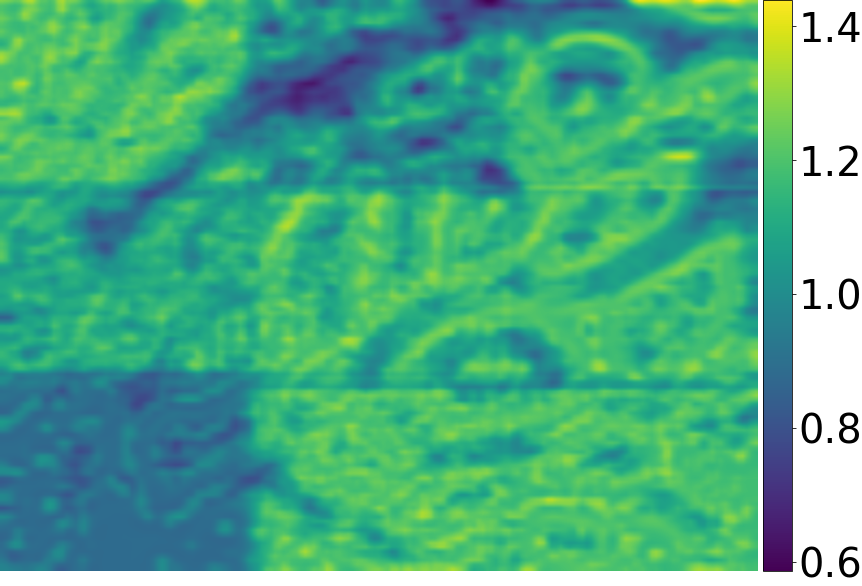}} 
    \caption{Tissue acid profiles for the tumor patterns  in Figure \ref{fig:threapySynthState}  \label{fig:threapySynthInState}}
\end{figure}

\section{Summary and conclusion} \label{sec:conclude}
In this paper we have formulated a terminal valued optimal control problem to understand the process involved in the formation of pseudopalisade structures during the progression of GBMs based on observed data.  Starting from the state of the art multiscale model \cite{Kumar2021} we proposed a modified model for the dynamics of pseudopalisade structures under the influence of tissue acidity which served as a dynamical state equation for the optimization problem. 
We then performed a wellposedness study based on which we are able to establish the existence of an optimal solution to the TOCP. This paved the way to performing numerical simulations for which we used the well known projected gradient descent method. For solving the TOCP problem the required data was obtained by taking experimentally observed images from literature. Here we also proposed an algorithm to convert raw experimental images to model-specific non-dimensionlized volumetric/concentration data. By using these processed images, we were able to successfully recreate the target patterns by estimating the optimal model parameter functions.  Because most of the developed mathematical models only rely on simulation results to  {reproduce experimentally observed qualitative behavior}, the proposed procedure provides an effective alternative approach to validate the model using real data. Based on the target-specific optimal parameters we were able to shed light not only on the dynamical interplay between reaction and migration terms, but also on the {relationship} between tumor progression and acidity. This type of data-specific analysis of dynamics {could be of particular interest to} medical professionals to perform patient-specific diagnosis and in turn {design} patient-specific treatment. From this perspective, we also showed the {feasibility} of different methods 
to normalize the tissue structure and obstruct tumor progression. The computed pattern neutralizing function achieves this not only by modifying the tissue acidity, but also by directly acting on the cancer cell population. Additionally, we highlighted the strength of this data-based approach by its ability to synthesize different unobserved pseudopalisade patterns, by simply combining already known optimal parameter functions computed from specific observed data. This can be further used to determine probable pattern-neutralization functions for a new unobserved pseudopalisade pattern by combining, in a similar way, the pattern-neutralizing functions of the simpler ones. Further direction of research would be to devise an active control strategy and a corresponding control problem for obstructing GBM progression. 

\newpage

\bibliographystyle{abbrv}

\bibliography{sn-bibliography}

\begin{thebibliography}{10}

\bibitem{PatB}
D.~P. Agamanolis.
\newblock Neuropathology: An illustrated interactive course for medical
  students and residents, 2017.
\newblock Available at
  \url{https://neuropathology-web.org/chapter7/chapter7bGliomas.html}.

\bibitem{Alfonso2016}
J.~C. Alfonso, A.~K{\"{o}}hn-Luque, T.~Stylianopoulos, F.~Feuerhake,
  A.~Deutsch, and H.~Hatzikirou.
\newblock {Why one-size-fits-all vaso-modulatory interventions fail to control
  glioma invasion: In silico insights}.
\newblock {\em Scientific Reports}, 6(October):1--15, 2016.

\bibitem{Alfonso2017}
J.~C. Alfonso, K.~Talkenberger, M.~Seifert, B.~Klink, A.~Hawkins-Daarud, K.~R.
  Swanson, H.~Hatzikirou, and A.~Deutsch.
\newblock {The biology and mathematical modelling of glioma invasion: A
  review}.
\newblock {\em Journal of the Royal Society Interface}, 14(136), 2017.

\bibitem{Amiri2016}
A.~Amiri, P.~U. Le, A.~Moquin, G.~Machkalyan, K.~Petrecca, J.~W. Gillard,
  N.~Yoganathan, and D.~Maysinger.
\newblock Inhibition of carbonic anhydrase ix in glioblastoma multiforme.
\newblock {\em European Journal of Pharmaceutics and Biopharmaceutics},
  109:81--92, 2016.

\bibitem{Becker2007}
R.~Becker and B.~Vexler.
\newblock {Optimal control of the convection-diffusion equation using
  stabilized finite element methods}.
\newblock {\em Numerische Mathematik}, 106(3):349--367, 2007.

\bibitem{Bottger2012}
K.~B{\"{o}}ttger, H.~Hatzikirou, A.~Chauviere, and A.~Deutsch.
\newblock {Investigation of the migration/proliferation dichotomy and its
  impact on avascular glioma invasion}.
\newblock {\em Mathematical Modelling of Natural Phenomena}, 7(1):105--135,
  2012.

\bibitem{Horn2007}
M.~C. Brahimi-Horn and J.~Pouyss{\'{e}}gur.
\newblock {Hypoxia in cancer cell metabolism and pH regulation}.
\newblock {\em Essays in Biochemistry}, 43:165--178, aug 2007.

\bibitem{Brat2003}
D.~Brat and T.~Mapstone.
\newblock {Malignant Glioma Physiology: Cellular Response to Hypoxia and Its
  Role in Tumor Progression}.
\newblock {\em Annals of internal medicine}, 138:659--668, may 2003.

\bibitem{Brat2004}
D.~J. Brat, A.~A. Castellano-Sanchez, S.~B. Hunter, M.~Pecot, C.~Cohen, E.~H.
  Hammond, S.~N. Devi, B.~Kaur, and E.~G. {Van Meir}.
\newblock {Pseudopalisades in Glioblastoma Are Hypoxic, Express Extracellular
  Matrix Proteases, and Are Formed by an Actively Migrating Cell Population}.
\newblock {\em Cancer Research}, 64(3):920--927, 2004.

\bibitem{PatA}
D.~J. Brat, A.~A. Castellano-Sanchez, S.~B. Hunter, M.~Pecot, C.~Cohen, E.~H.
  Hammond, S.~N. Devi, B.~Kaur, and E.~G. Van~Meir.
\newblock {Pseudopalisades in Glioblastoma Are Hypoxic, Express Extracellular
  Matrix Proteases, and Are Formed by an Actively Migrating Cell Population}.
\newblock {\em Cancer Research}, 64(3):920--927, 02 2004.

\bibitem{Caiazzo2015}
A.~Caiazzo and I.~Ramis-Conde.
\newblock {Multiscale modelling of palisade formation in gliobastoma
  multiforme}.
\newblock {\em Journal of Theoretical Biology}, 383:145--156, 2015.

\bibitem{Chiche2010}
J.~Chiche, M.~C. Brahimi-Horn, and J.~Pouyss{\'e}gur.
\newblock {Tumour hypoxia induces a metabolic shift causing acidosis: A common
  feature in cancer}.
\newblock {\em Journal of Cellular and Molecular Medicine}, 14(4):771--794,
  2010.

\bibitem{Collis2002}
S.~S. Collis and M.~Heinkenschloss.
\newblock {Analysis of SUPG Method Applied to the Solution of Optimal Control
  Problems}, March 2002.
\newblock Accessible at \url{https://hdl.handle.net/1911/101983}.

\bibitem{Colombo2015}
M.~C. Colombo, C.~Giverso, E.~Faggiano, C.~Boffano, F.~Acerbi, and
  P.~Ciarletta.
\newblock {Towards the personalized treatment of glioblastoma: Integrating
  patient-specific clinical data in a continuous mechanical model}.
\newblock {\em PLoS ONE}, 10(7):1--23, 2015.

\bibitem{Conte2022}
M.~Conte, Y.~Dzierma, S.~Knobe, and C.~Surulescu.
\newblock Mathematical modeling of glioma invasion and therapy approaches.
\newblock {\em arXiv:2203.11578}, 2022.

\bibitem{Conte2020}
M.~Conte, L.~Gerardo-Giorda, and M.~Groppi.
\newblock {Glioma invasion and its interplay with nervous tissue and therapy: A
  multiscale model}.
\newblock {\em Journal of Theoretical Biology}, 486:110088, 2020.

\bibitem{Conte2021}
M.~Conte and C.~Surulescu.
\newblock Mathematical modeling of glioma invasion: acid- and vasculature
  mediated go-or-grow dichotomy and the influence of tissue anisotropy.
\newblock {\em Applied Mathematics and Computation}, 407:126305, 2021.

\bibitem{Corbin2018}
G.~Corbin, A.~Hunt, F.~Schneider, A.~Klar, and C.~Surulescu.
\newblock {Higher-Order Models for Glioma Invasion: From a Two-Scale
  Description to Effective Equations for Mass density and Momentum}.
\newblock {\em Mathematical Models and Methods in Applied Sciences}, 28, jan
  2018.

\bibitem{Corbin2021}
G.~Corbin, A.~Klar, C.~Surulescu, C.~Engwer, M.~Wenske, J.~Nieto, and J.~Soler.
\newblock Modeling glioma invasion with anisotropy- and hypoxia-triggered
  motility enhancement: From subcellular dynamics to macroscopic pdes with
  multiple taxis.
\newblock {\em Mathematical Models and Methods in Applied Sciences},
  31(01):177--222, 2021.

\bibitem{Dietrich2020}
A.~Dietrich, N.~Kolbe, N.~Sfakianakis, and C.~Surulescu.
\newblock Multiscale modeling of glioma invasion: from receptor binding to
  flux-limited macroscopic pdes.
\newblock {\em arXiv: Tissues and Organs}, 2020.
\newblock To appear in SIAM Multiscale Modeling and Simulation.

\bibitem{Dolecek2012}
T.~A. Dolecek, J.~M. Propp, N.~E. Stroup, and C.~Kruchko.
\newblock {CBTRUS statistical report: Primary brain and central nervous system
  tumors diagnosed in the United States in 2005-2009}.
\newblock {\em Neuro-Oncology}, 14(SUPPL.5), 2012.

\bibitem{Engwer2015}
C.~Engwer, T.~Hillen, M.~Knappitsch, and C.~Surulescu.
\newblock {Glioma follow white matter tracts: a multiscale DTI-based model}.
\newblock {\em Journal of Mathematical Biology}, 71(3):551--582, 2015.

\bibitem{Engwer2016a}
C.~Engwer, A.~Hunt, and C.~Surulescu.
\newblock {Effective equations for anisotropic glioma spread with
  proliferation: A multiscale approach and comparisons with previous settings}.
\newblock {\em Mathematical Medicine and Biology}, 33(4):435--459, 2016.

\bibitem{Engwer2016}
C.~Engwer, M.~Knappitsch, and C.~Surulescu.
\newblock {A multiscale model for glioma spread including cell-tissue
  interactions and proliferation}.
\newblock {\em Mathematical Biosciences and Engineering}, 13(2):443--460, 2016.

\bibitem{Estrella2013}
V.~Estrella, T.~Chen, M.~Lloyd, J.~Wojtkowiak, H.~H. Cornnell,
  A.~Ibrahim-Hashim, K.~Bailey, Y.~Balagurunathan, J.~M. Rothberg, B.~F.
  Sloane, J.~Johnson, R.~A. Gatenby, and R.~J. Gillies.
\newblock {Acidity generated by the tumor microenvironment drives local
  invasion}.
\newblock {\em Cancer Research}, 73(5):1524--1535, 2013.

\bibitem{PatC}
J.~Florian.
\newblock Glioblastoma showing areas of pseudopalisading necrosis, 2010.
\newblock Available at
  \url{https://commons.wikimedia.org/wiki/File:GBM_pseudopalisading_necrosis.jpg}.

\bibitem{Gatenby03}
R.~Gatenby and E.~T. Gawlinski.
\newblock The glycolytic phenotype in carcinogenesis and tumor invasion:
  insights through mathematical models.
\newblock {\em Cancer research}, 63 14:3847--54, 2003.

\bibitem{Gholami2016}
A.~Gholami, A.~Mang, and G.~Biros.
\newblock {An inverse problem formulation for parameter estimation of a
  reaction–diffusion model of low grade gliomas}.
\newblock {\em Journal of Mathematical Biology}, 72(1-2):409--433, 2016.

\bibitem{Harpold2007}
H.~L. Harpold, E.~C. Alvord, and K.~R. Swanson.
\newblock {The evolution of mathematical modeling of glioma proliferation and
  invasion}.
\newblock {\em Journal of Neuropathology and Experimental Neurology},
  66(1):1--9, 2007.

\bibitem{Hatzikirou2012}
H.~Hatzikirou, D.~Basanta, M.~Simon, K.~Schaller, and A.~Deutsch.
\newblock {'Go or grow': The key to the emergence of invasion in tumour
  progression?}
\newblock {\em Mathematical Medicine and Biology}, 29(1):49--65, 2012.

\bibitem{Hatzikirou2005}
H.~Hatzikirou, A.~Deutsch, C.~Schaller, M.~Simon, and K.~Swanson.
\newblock {Mathematical modelling of glioblastoma tumour development: A
  review}.
\newblock {\em Mathematical Models and Methods in Applied Sciences},
  15(11):1779--1794, 2005.

\bibitem{Finotti12}
T.~V.~P. {Heather Finotti Suzanne Lenhart}.
\newblock {Optimal control of advective direction in reaction-diffusion
  population models}.
\newblock {\em Evolution Equations \& Control Theory}, 1(1):81--107, 2012.

\bibitem{Hinow2008}
{Hinow, P., P. Gerlee, et al.}
\newblock A spatial model of tumor-host interaction: Application of
  chemotherapy.
\newblock {\em Math Biosci Eng}, 6(3):521--546, 2009.

\bibitem{Hinze08}
M.~Hinze, R.~Pinnau, M.~Ulbrich, and S.~Ulbrich.
\newblock {\em {Optimization with PDE Constraints}}.
\newblock Mathematical Modelling: Theory and Applications. Springer
  Netherlands, 2008.

\bibitem{Hiremath2015}
S.~Hiremath and C.~Surulescu.
\newblock {A stochastic multiscale model for acid mediated cancer invasion}.
\newblock {\em Nonlinear Analysis: Real World Applications}, 22:176--205, 2015.

\bibitem{Hiremath2016}
S.~A. Hiremath and C.~Surulescu.
\newblock {A stochastic model featuring acid-induced gaps during tumor
  progression}.
\newblock {\em Nonlinearity}, 29(3):851--914, 2016.

\bibitem{Hiremath2017}
S.~A. Hiremath and C.~Surulescu.
\newblock {Mathematical Models for Acid-Mediated Tumor Invasion: From
  Deterministic to Stochastic Approaches BT - Multiscale Models in Mechano and
  Tumor Biology}.
\newblock pages 45--71, Cham, 2017. Springer International Publishing.

\bibitem{Hiremath2018}
S.~A. Hiremath, C.~Surulescu, A.~Zhigun, and S.~Sonner.
\newblock {On a coupled SDE-PDE system modeling acid-mediated tumor invasion}.
\newblock {\em Discrete and Continuous Dynamical Systems - Series B},
  23(9):3685--3715, 2018.

\bibitem{Hogea2008}
C.~Hogea, C.~Davatzikos, and G.~Biros.
\newblock {An image-driven parameter estimation problem for a
  reaction-diffusion glioma growth model with mass effects}.
\newblock {\em Journal of Mathematical Biology}, 56(6):793--825, 2008.

\bibitem{Hoering2012}
E.~H\"{o}ring, P.~Harter, J.~Seznec, J.~Schittenhelm, H.-J. B\"{u}hring,
  S.~Bhattacharyya, E.~von Hattingen, C.~Zachskorn, M.~Mittelbronn, and
  U.~Naumann.
\newblock The go or grow potential of gliomas is linked to the neuropeptide
  processing enzyme carboxypeptidase e and mediated by metabolic stress.
\newblock {\em Acta Neuropathologica}, 124(1):83--97, 2012.

\bibitem{Hunt2016}
A.~Hunt and C.~Surulescu.
\newblock {A Multiscale Modeling Approach to Glioma Invasion with Therapy}.
\newblock {\em Vietnam Journal of Mathematics}, 45, jul 2016.

\bibitem{Jbabdi2005}
S.~Jbabdi, E.~Mandonnet, H.~Duffau, L.~Capelle, K.~R. Swanson,
  M.~P{\'e}l{\'e}grini-Issac, R.~Guillevin, and H.~Benali.
\newblock {Simulation of anisotropic growth of low-grade gliomas using
  diffusion tensor imaging}.
\newblock {\em Magnetic Resonance in Medicine}, 54(3):616--624, 2005.

\bibitem{Jing2019}
X.~Jing, F.~Yang, C.~Shao, K.~Wei, M.~Xie, H.~Shen, and Y.~Shu.
\newblock {Role of hypoxia in cancer therapy by regulating the tumor
  microenvironment}.
\newblock {\em Molecular Cancer}, 18(1):1--15, 2019.

\bibitem{Khain2011}
E.~Khain, M.~Katakowski, S.~Hopkins, A.~Szalad, X.~Zheng, F.~Jiang, and
  M.~Chopp.
\newblock {Collective behavior of brain tumor cells: the role of hypoxia.}
\newblock {\em Physical review. E, Statistical, nonlinear, and soft matter
  physics}, 83(3 Pt 1):31920, mar 2011.

\bibitem{Kim2009}
Y.~Kim, S.~Lawler, M.~O. Nowicki, E.~A. Chiocca, and A.~Friedman.
\newblock {A mathematical model for pattern formation of glioma cells outside
  the tumor spheroid core}.
\newblock {\em Journal of Theoretical Biology}, 260(3):359--371, 2009.

\bibitem{Kim2013}
Y.~Kim and S.~Roh.
\newblock {A hybrid model for cell proliferation and migration in
  glioblastoma}.
\newblock {\em Discrete and Continuous Dynamical Systems - Series B},
  18(4):969--1015, 2013.

\bibitem{Kleihues1995}
P.~Kleihues, F.~Soylemezoglu, B.~Sch{\"{a}}uble, B.~W. Scheithauer, and P.~C.
  Burger.
\newblock {Histopathology, classification, and grading of gliomas.}
\newblock {\em Glia}, 15(3):211--221, nov 1995.

\bibitem{Konukoglu2010}
E.~Konukoglu, O.~Clatz, P.~Y. Bondiau, H.~Delingette, and N.~Ayache.
\newblock {Extrapolating glioma invasion margin in brain magnetic resonance
  images: Suggesting new irradiation margins}.
\newblock {\em Medical Image Analysis}, 14(2):111--125, 2010.

\bibitem{Kumar2021}
P.~Kumar, J.~Li, and C.~Surulescu.
\newblock {Multiscale modeling of glioma pseudopalisades: contributions from
  the tumor microenvironment}.
\newblock {\em Journal of Mathematical Biology}, 82(6):1--45, 2021.

\bibitem{KS2020}
P.~Kumar and C.~Surulescu.
\newblock A flux-limited model for glioma patterning with hypoxia-induced
  angiogenesis.
\newblock {\em Symmetry}, 12(11), 2020.

\bibitem{KSZ21}
P.~Kumar, C.~Surulescu, and A.~Zhigun.
\newblock Multiphase modelling of glioma pseudopalisading under acidosis.
\newblock {\em Mathematics in Engineering}, 4(6):1--28, 2022.

\bibitem{Gonzalez2012}
A.~Mart{\'{i}}nez-Gonz{\'{a}}lez, G.~F. Calvo, L.~A. {P{\'{e}}rez Romasanta},
  and V.~M. P{\'{e}}rez-Garc{\'{i}}a.
\newblock {Hypoxic Cell Waves Around Necrotic Cores in Glioblastoma: A
  Biomathematical Model and Its Therapeutic Implications}.
\newblock {\em Bulletin of Mathematical Biology}, 74(12):2875--2896, 2012.

\bibitem{Martirosyan2015}
N.~L. Martirosyan, E.~M. Rutter, W.~L. Ramey, E.~J. Kostelich, Y.~Kuang, and
  M.~C. Preul.
\newblock {Mathematically modeling the biological properties of gliomas: A
  review}.
\newblock {\em Mathematical Biosciences and Engineering}, 12(4):879--905, 2015.

\bibitem{PatE}
C.~McKinney.
\newblock Glioblastoma multiforme pseudopalisading necrosis.
\newblock Available at \url{https://www.pinterest.com/pin/471048442246524817/}.

\bibitem{Murray2002}
J.~D. Murray.
\newblock {\em Mathematical Biology I. An Introduction}, volume~17 of {\em
  Interdisciplinary Applied Mathematics}.
\newblock Springer, New York, 3 edition, 2002.

\bibitem{Necas19}
J.~Necas, J.~Malek, M.~Rokyta, and M.~Ruzicka.
\newblock {\em {Weak and Measure-Valued Solutions to Evolutionary PDEs}}.
\newblock Chapman and Hall/CRC Press, 1996.

\bibitem{Painter2013}
K.~J. Painter and T.~Hillen.
\newblock {Mathematical modelling of glioma growth: The use of Diffusion Tensor
  Imaging (DTI) data to predict the anisotropic pathways of cancer invasion}.
\newblock {\em Journal of Theoretical Biology}, 323:25--39, 2013.

\bibitem{Piasentin2020}
N.~Piasentin, E.~Milotti, and R.~Chignola.
\newblock {The control of acidity in tumor cells: a biophysical model}.
\newblock {\em Scientific Reports}, 10(1):1--14, 2020.

\bibitem{Plate1992}
K.~H. Plate, G.~Breier, H.~A. Weich, and W.~Risau.
\newblock {Vascular endothelial growth factor is a potential tumour
  angiogenesis factor in human gliomas in vivo}.
\newblock {\em Nature}, 359(6398):845--848, 1992.

\bibitem{Rong2006}
Y.~Rong, D.~L. Durden, E.~G. {Van Meir}, and D.~J. Brat.
\newblock {'Pseudopalisading' necrosis in glioblastoma: A familiar morphologic
  feature that links vascular pathology, hypoxia, and angiogenesis}.
\newblock {\em Journal of Neuropathology and Experimental Neurology},
  65(6):529--539, 2006.

\bibitem{Sander2002}
L.~M. Sander and T.~S. Deisboeck.
\newblock {Growth patterns of microscopic brain tumors.}
\newblock {\em Physical review. E, Statistical, nonlinear, and soft matter
  physics}, 66(5 Pt 1):51901, nov 2002.

\bibitem{Sturrock2015}
M.~Sturrock, W.~Hao, J.~Schwartzbaum, and G.~A. Rempala.
\newblock {A mathematical model of pre-diagnostic glioma growth}.
\newblock {\em Journal of theoretical biology}, 380:299--308, sep 2015.

\bibitem{Swan2018}
A.~Swan, T.~Hillen, J.~C. Bowman, and A.~D. Murtha.
\newblock {A Patient-Specific Anisotropic Diffusion Model for Brain Tumour
  Spread}.
\newblock {\em Bulletin of Mathematical Biology}, 80(5):1259--1291, 2018.

\bibitem{Swanson2011}
K.~R. Swanson, R.~C. Rockne, J.~Claridge, M.~A. Chaplain, E.~C. Alvord, and
  A.~R. Anderson.
\newblock {Quantifying the role of angiogenesis in malignant progression of
  gliomas: In Silico modeling integrates imaging and histology}.
\newblock {\em Cancer Research}, 71(24):7366--7375, 2011.

\bibitem{Webb99}
S.~Webb, J.~Sherratt, and R.~Fish.
\newblock Mathematical modelling of tumor acidity: Regulation of intracellular
  ph.
\newblock {\em Journal of Theoretical Biology}, 196(2):237--250, 1999.

\bibitem{Webb04}
S.~Webb, J.~Sherratt, and R.~Fish.
\newblock Alterations in proteolytic activity at low ph and its association
  with invasion: A theoretical model.
\newblock {\em Clinical \& Experimental Metastasis}, 17:397--407, 2004.

\bibitem{Wippold2006}
F.~J.~n. Wippold, M.~L{\"{a}}mmle, F.~Anatelli, J.~Lennerz, and A.~Perry.
\newblock {Neuropathology for the neuroradiologist: palisades and
  pseudopalisades.}
\newblock {\em AJNR. American journal of neuroradiology}, 27(10):2037--2041,
  2006.

\bibitem{Yagi09}
A.~Yagi.
\newblock {\em Abstract Parabolic Evolution Equations and Their Applications}.
\newblock Springer Monographs in Mathematics. Springer, 2009.

\bibitem{YZY2020}
M.~Yang, X.~Zhong, and Y.~Yuan.
\newblock Does baking soda function as a magic bullet for patients with cancer?
  a mini review.
\newblock {\em Integrative Cancer Therapies}, 19:1534735420922579, 2020.
\newblock PMID: 32448009.

\bibitem{Yu2021}
V.~Y. Yu, D.~Nguyen, D.~OConnor, D.~Ruan, T.~Kaprealian, R.~Chin, and K.~Sheng.
\newblock {Treating Glioblastoma Multiforme (GBM) with super hyperfractionated
  radiation therapy: Implication of temporal dose fractionation optimization
  including cancer stem cell dynamics}.
\newblock {\em PLoS ONE}, 16(2 February):1--16, 2021.

\bibitem{Zacher2010}
R.~Zacher.
\newblock {De Giorgi-Nash-Moser} estimates for evolutionary partial
  integro-differential equations, 2010.
\newblock Habilitationsschrift from Univ. Halle-Wittenberg, accessible at
  \url{https://dx.doi.org/10.25673/387}.

\bibitem{Zagzag2000}
D.~Zagzag, R.~Amirnovin, M.~A. Greco, H.~Yee, J.~Holash, S.~J. Wiegand,
  S.~Zabski, G.~D. Yancopoulos, and M.~Grumet.
\newblock {Vascular apoptosis and involution in gliomas precede
  neovascularization: A novel concept for glioma growth and angiogenesis}.
\newblock {\em Laboratory Investigation}, 80(6):837--849, 2000.

\end{thebibliography}

\end{document}